\documentclass[authoryear, 11pt]{imsart}
\RequirePackage[OT1]{fontenc}
\RequirePackage{amsthm,amsmath,natbib}
\RequirePackage{hypernat}
\startlocaldefs
\numberwithin{equation}{section}
\theoremstyle{plain}

\endlocaldefs
\usepackage{psfrag, epsfig, graphicx}
\usepackage{amsfonts,amsmath,latexsym,amssymb}
\usepackage[active]{srcltx}
\usepackage{fullpage}

\newtheorem{lemma}{Lemma}
\newtheorem{theorem}{Theorem}
\newtheorem{proposition}{Proposition}

\newtheorem{assumption}{Assumption}
\newtheorem{remark}{Remark}
\newtheorem{definition}{Definition}
\newtheorem{corollary}{Corollary}

\renewcommand{\kappa}{\varkappa}

\newcommand{\rd}{{\rm d}}
\newcommand{\ve}{\varepsilon}
\newcommand{\e}{\varepsilon}

\newcommand{\cA}{{\cal A}}
\newcommand{\cB}{{\cal B}}

\newcommand{\cE}{{\cal E}}
\newcommand{\cF}{{\cal F}}

\newcommand{\cJ}{{\cal J}}
\newcommand{\cK}{{\cal K}}
\newcommand{\cL}{{\cal L}}
\newcommand{\cM}{{\cal M}}

\newcommand{\cP}{{\cal P}}

\newcommand{\cR}{{\cal R}}
\newcommand{\cS}{{\cal S}}
\newcommand{\cT}{{\cal T}}

\newcommand{\cY}{{\cal Y}}

\newcommand{\bga}{\boldsymbol{\gamma}}
\newcommand{\bde}{\boldsymbol{\delta}}

\newcommand{\bmu}{\boldsymbol{\mu}}
\newcommand{\bleta}{\boldsymbol{\eta}}
\newcommand{\bB}{\mathbb B}
\newcommand{\bC}{\mathbb C}

\newcommand{\bE}{\mathbb E}
\newcommand{\bF}{\mathbb F}

\newcommand{\bH}{\mathbb H}

\newcommand{\bK}{{\mathbb K}}
\newcommand{\bL}{{\mathbb L}}

\newcommand{\bN}{{\mathbb N}}

\newcommand{\bP}{{\mathbb P}}

\newcommand{\bR}{{\mathbb R}}
\newcommand{\bS}{{\mathbb S}}

\newcommand{\mA}{\mathfrak{A}}
\newcommand{\mB}{\mathfrak{B}}

\newcommand{\mS}{\mathfrak{S}}
\newcommand{\mL}{\mathfrak{L}}

\newcommand{\mmp}{\mathfrak{p}}

\newcommand{\mK}{\mathfrak{K}}

\newcommand{\mh}{\mathfrak{h}}

\newcommand{\mJ}{\mathfrak{J}}
\newcommand{\mm}{\mathfrak{m}}

\newcommand{\mH}{\mathfrak{H}}

\newcommand{\ma}{\mathfrak{a}}

\newcommand{\epr}{\hfill\hbox{\hskip 4pt
                \vrule width 5pt height 6pt depth 1.5pt}\vspace{0.5cm}\par}

\newcommand{\blh}{\boldsymbol{h}}

\begin{document}
\begin{frontmatter}
\title{
Adaptive estimation over anisotropic functional classes via oracle approach}
\runtitle{Global choice}
\begin{aug}
\author{\fnms{Oleg} \snm{Lepski}
\ead[label=e2]{lepski@cmi.univ-mrs.fr}}
\runauthor{O. Lepski}



\address{Institut de Math\'ematique de Marseille\\
Aix-Marseille  Universit\'e   \\
 39, rue F. Joliot-Curie \\
13453 Marseille, France\\
\printead{e2}\\ }
\end{aug}

\maketitle

\begin{abstract}
We address the problem of adaptive minimax  estimation in white gaussian noise model
under $\bL_p$--loss, $1\leq p\leq\infty,$ on the anisotropic Nikolskii classes.
We present  the estimation procedure  based on a  \textsf{new} data-driven selection scheme from the family of
kernel estimators with varying bandwidths. For proposed estimator  we establish so-called  $\bL_p$-norm oracle inequality
and use it  for deriving minimax adaptive results. We prove the existence of rate-adaptive estimators and
 fully characterize behavior of the minimax risk for different
relationships between regularity parameters and
norm indexes in  definitions of the functional class and of the risk.
In particular some new asymptotics of the minimax risk are discovered including necessary and sufficient conditions for existence a uniformly consistent estimator. We provide also with  detailed overview of existing methods and results and formulate open problems in adaptive minimax estimation.


\end{abstract}


\end{frontmatter}

\section{Introduction}
\label{sec:introduction}
Let  $\bR^d,\;d\geq 1,$ be equipped with  Borel $\sigma$-algebra $\mB(\bR^d)$  and Lebesgue measure $\nu_d$. Put $\widetilde{\mB}(\bR^d)=\left\{B\in\mB(\bR^d):\;\;\nu_d(B)<\infty\right\}$ and let $\big(W(B),\; B\in\widetilde{\mB}(\bR^d)\big)$ be the white noise with intensity
$\nu_d$. Set also for any $\cA\in\mB(\bR^d)$ and any $1\leq p<\infty$
\begin{eqnarray*}
\bL_p\big(\cA,\nu_d\big)&=&\Big\{g:\cA\to\bR:\; \|g\|_{p,\cA}^p:=\int_{\cA}|g(t)|^p\nu_d(\rd t)<\infty\Big\};
\\
\bL_\infty\big(\cA\big)&=&\Big\{g:\cA\to\bR:\; \|g\|_{\infty,\cA}:=\sup_{t\in\cA}|g(t)|<\infty\Big\}.
\end{eqnarray*}

\subsection{Statistical model and $\bL_p$-risk}

Consider the sequence of  statistical experiments (called gaussian white noise model) generated by the observation $X^{\e}=\left\{X_\epsilon(g),\;g\in\bL_2\big(\bR^d,\nu_d\big)\right\}_{\e}$ where
\begin{equation}
\label{eq:WGN-model}
X_\e(g)=\int f(t)g(t)\nu_d(\rd t)+\e\int g(t)W(\rd t).
\end{equation}
Here $\e\in (0,1)$ is understood as the noise level  which is usually supposed  sufficiently  small.
\par
The goal is to recover unknown signal $f$ from observation $X^{\e}$ on a given cube $(-b,b)^d,\; b>0$.
The quality of an estimation procedure will be described by $\bL_p$-risk, $1\leq p\leq\infty,$ defined in (\ref{eq:L_p-risk}) below and
as an estimator we understand any $X^{\e}$-measurable Borel function belonging to $\bL_p\big(\bR^d,\nu_d\big)$. Without loss of generality and for ease of the notation we will assume that functions to be estimated
vanish outside $(-b,b)^d$.

Thus, for any estimator $\tilde{f}_\e$
and any $f\in\bL_p\big(\bR^d,\nu_d\big)\cap\bL_2\big(\bR^d,\nu_d\big)$ we define its $\bL_p$-risk  as
\begin{equation}
\label{eq:L_p-risk}
\cR^{(p)}_\e\big[\tilde{f}_\e; f\big]=\left\{\bE^{(\e)}_f\left(\big\|\tilde{f}_\e-f\big\|^q_{p}\right)\right\}^{\frac{1}{q}},\;\; q\geq 1.
\end{equation}
Here and later $\|\cdot\|_p, 1\leq p\leq \infty,$ stands for $\|\cdot\|_{p,(-b,b)^d}$ and $\bE^{(\e)}_f$ denote the mathematical expectation with respect to the probability law of $X^\e$.

\smallskip

Let $\bF$ be  a given subset of $\bL_p\big(\bR^d,\nu_d\big)\cap\bL_2\big(\bR^d,\nu_d\big)$. For any estimator $\tilde{f}_\e$ define its {\it maximal risk}  by
$
\cR^{(p)}_\e\big[\tilde{f}_\e; \bF\big]=\sup_{f\in\bF}\cR^{(p)}_\e\big[\tilde{f}_\e; f\big]
$
and its {\it minimax risk} on $\bF$ is given by
\begin{equation}
\label{eq:minmax-risk}
\phi_\e(\bF):=\inf_{\tilde{f}_\e}\cR^{(p)}_\e\big[\tilde{f}_\e; \bF\big].
\end{equation}
Here infimum is taken over all possible estimators. An estimator whose maximal risk is proportional to $\phi_\e(\bF)$ is called minimax on $\bF$.

\subsection{{Adaptive estimation}}

Let $\big\{\bF_\vartheta,\vartheta\in\Theta\big\}$ be the collection of subsets of $\bL_p\big(\bR^d,\nu_d\big)\cap\bL_2\big(\bR^d,\nu_d\big)$, where $\vartheta$ is a nuisance parameter which may have very complicated structure.

The problem of adaptive estimation can be formulated as follows:
{\it is it possible to construct a single estimator $\hat{f}_\e$
 which would be  simultaneously minimax on each class
 $\bF_\vartheta,\;\vartheta\in\Theta$, i.e.}
$$
 \cR^{(p)}_\e\big[\hat{f}_\e; \bF_\vartheta\big]\sim \phi_\e(\bF_\vartheta),\;\e\to 0,\;\;\forall \vartheta\in\Theta?
$$
We refer to this question as {\it the  problem of adaptive
estimation over  the scale of }
$\{\bF_\vartheta,\;\vartheta\in\Theta \}$.
If such estimator exists we will call it optimally or rate-adaptive.

In the present paper we will be interested in adaptive estimation over the scale
$$
\bF_\vartheta=\bN_{\vec{r},d}\big(\vec{\beta},\vec{L}\big),\;\vartheta=\big(\vec{\beta},\vec{r},\vec{L}\big),
$$
where $\bN_{\vec{r},d}\big(\vec{\beta},\vec{L}\big)$ is an anisotropic Nikolskii class, see Section \ref{sec:nikolski} for formal definition.
Here we only mention that for any $f\in\bN_{\vec{r},d}\big(\vec{\beta},\vec{L}\big)$ the coordinate $\beta_i$ of the vector  $\vec{\beta}=(\beta_1,\ldots,\beta_d)\in (0,\infty)^d$ represents the smoothness of $f$ in the direction $i$ and  the coordinate $r_i$ of the vector
$\vec{r}=(r_1,\ldots,r_d)\in [1,\infty]^d$ represents the index of the norm in which $\beta_i$ is measured. Moreover, $\bN_{\vec{r},d}\big(\vec{\beta},\vec{L}\big)$  is the intersection of the balls in some semi-metric space and the vector $\vec{L}\in (0,\infty)^d$ represents the radii of these balls.

The aforementioned dependence on the direction is usually referred to {\it anisotropy} of the underlying function and the corresponding functional class.
The use of the integral norm in the definition of  the smoothness   is referred to {\it inhomogeneity} of the underlying function. The latter means that the function $f$ can be sufficiently smooth on some part of the  observation domain  and rather irregular on the other part. Thus,  the adaptive  estimation over the scale $\big\{\bN_{\vec{r},d}\big(\vec{\beta},\vec{L}\big),\;\big(\vec{\beta},\vec{r},\vec{L}\big)\in(0,\infty)^d\times[1,\infty]^d\times(0,\infty)^d\big\}$
can be viewed as the adaptation to anisotropy  and inhomogeneity of the function to be estimated.

\subsection{Historical notes} The history of the adaptive estimation over  scales of sets of smooth functions counts nowadays 30 years. During this time the variety of functional classes was introduced in the nonparametric statistics in particular Sobolev, Nikolskii and Besov ones.
The relations between different scales as well as between classes belonging to the same scale can be found, for instance, in \cite{Nikolski}.
It is worth mentioning that although  considered classes are  different,  the same estimation procedure may be minimax on them.
In such situations we will say that the class $\bF_1$ {\it  statistically equivalent} to the class $\bF_2$ and write $\bF_1\bowtie\bF_2$. Also for two sequences $a_\e\to 0$ and $b_\e\to 0$ we will write $a_\e\sim b_\e$  and $a_\e\gtrsim b_\e$ if $0<\lim_{\e\to 0}a_\e b^{-1}_\e<\infty$ and $\lim_{\e\to 0}a_\e b^{-1}_\e\geq 1$ respectively.

\paragraph{Estimation of univariate functions}

The first adaptive results were obtained  in \cite{Pinsker}. The authors studied the problem of adaptive estimation  over the scale of periodic Sobolev classes (Sobolev ellipsoids), $W(\beta,L)$, in the univariate model (\ref{eq:WGN-model}) under $\bL_2$-loss ($p=2$).
The exact asymptotics of minimax risk on $W(\beta,L)$ is given by $P(L)\e^{\frac{2\beta}{2\beta+1}}$, where $P(L)$ is the Pinsker constant.  The authors proposed the  estimation procedure based on blockwise bayesian  construction and showed that it is adaptive {\it efficient} over the scale of considered classes.
Noting that $W(\beta,L)\bowtie\bN_{2,1}(\beta,L)$ one can assert that Efroimovich-Pinsker estimator is rate-adaptive on $\bN_{2,1}(\beta,L)$ as well.

Starting from this pioneering paper a variety of adaptive methods under $\bL_2$-loss were proposed in different  statistical models such as density and spectral density estimation, nonparametric regression, deconvolution model, inverse problems and many others. Let us mention some of them.
\begin{itemize}
\item
Extension of   Efroimovich-Pinsker method, \cite{Efr, Efr2};
\item
Unbiased risk minimization, \cite{golubev92, golubev93};
\item
Model selection, \cite{Barron-Birge, B-M, birge};
\item
Aggregation of estimators, \cite{nemirovski00}, \cite{Nem},   \cite{Weg}, \cite{Tsyb03}, \cite{rigollet-tsybakov}, \cite{Bunea}, \cite{Gol};
\item
Exponential weights, \cite{Barron}, \cite{dal08}, \cite{rigollet-tsybakov11};
\item
Risk hull method, \cite{cav-golubev};
\item
Blockwise Stein method, \cite{C99}, \cite{cav-tsyb},  \cite{rigollet}.
\end{itemize}

Some of aforementioned papers deal with not only adaptation over the scale of functional classes but contain sharp oracle inequalities
(about oracle approach and its relation to adaptive estimation see for instance \cite{GL12} and the references therein). Without any doubts the adaptation  under $\bL_2$-loss is the best developed area of the adaptive estimation.
Rather detailed overview and some new ideas related to this topic  can be found in the recent paper \cite{Baraud}.

The adaptive estimation under $\bL_p$-loss, $1\leq p\leq \infty$ was initiated in \cite{lepski91} over the collection of H\"older classes, i.e.
$\bN_{\infty,1}(\beta,L)$.  The asymptotics of minimax risk  is given by
$$
\phi\big(\bN_{\infty,1}(\beta,L)\big)\sim
\left\{
\begin{array}{ll}
\e^{\frac{2\beta}{2\beta+1}},\quad& p\in[1,\infty);
\\
\big(\e^2|\ln(\e)|\big)^{\frac{\beta}{2\beta+1}},\quad& p=\infty.
\end{array}
\right.
$$
The author constructed the optimally-adaptive estimator which  is obtained by the selection from the family of
piecewise polynomial estimators. Selection rule is based on pairwise comparison of estimators (bias-majorant tradeoff).
Some sharp results were obtain in \cite{lepski92b}, where {\it efficient} adaptive estimator was proposed in the case of $\bL_\infty$-loss, see also \cite{Tsyb}.

Recent development in adaptive univariate density estimation under $\bL_\infty$-loss can be found in \cite{GN-1}, \cite{Gach}. Another "extreme" case, the estimation under $\bL_1$-loss,
was scrutinized in \cite{dev-lug96}, \cite{dev-lug97}.

\smallskip

The consideration of the classes of inhomogeneous functions  in nonparametric statistics has been started in \cite{nemirovski}, where the minimax rates of convergence were established and minimax estimators were constructed in the case of generalized Sobolev classes. The adaptive estimation problem over the scale of Besov classes $\bB^{\beta}_{r,q}(L)$ was studied for the first time in \cite{Donoho} in the framework of the density model. We note  that $\bB^{\beta}_{r,\infty}=\bN_{r,1}(\beta,L)$ and although
$\bB^{\beta}_{r,q}\supset\bN_{r,1}(\beta,L)$ for any $q\geq 1$, see \cite{Nikolski}, one has $\bB^{\beta}_{r,q}\bowtie\bN_{r,1}(\beta,L)$.

The same problem in the  univariate model (\ref{eq:WGN-model}) was studied in  \cite{LepMam}. The asymptotics of minimax risk is given by

$$
\phi\big(\bB^{\beta}_{r,q}(L)\big)\sim
\left\{
\begin{array}{ll}
\e^{\frac{2\beta}{2\beta+1}},\quad& (2\beta+1)r>p;
\\
\big(\e^2|\ln(\e)|\big)^{\frac{\beta-1/r+1/p}{2\beta-2/r+1}},\quad& (2\beta+1)r\leq p.
\end{array}
\right.
$$

The set of parameters satisfying  $r(2\beta+1)>p$  is called in the literature {\em the dense zone} and the case  $r(2\beta+1)\leq p$ is referred to {\em the sparse zone}. As it was shown  in \cite{Donoho} hard threshold wavelet estimator is {\it nearly} adaptive over the scale of Besov classes. The latter means that the maximal risk of the proposed estimator differs from $\phi\big(\bB^{\beta}_{r,q}(L)\big)$ by logarithmic factor on the dense zone and on the boundary $(2\beta+1)r=p$. The similar  result was proved in \cite{LepMam} but for completely different estimation procedure: for the first time local bandwidth selection scheme was used for the estimation of entire function. Moreover, the computations of the  maximal risk of the proposed estimator on
$\bB^{\beta}_{r,q}(L)$ was made  by   integration of the local oracle inequality.

It is important to emphasize that both aforementioned results were proved under additional assumption
\begin{equation}
\label{eq1:class-param}
1 -(\beta r)^{-1}+(\beta p)^{-1}>0.
\end{equation}
Independently, the approach similar to \cite{LepMam} was proposed in \cite{GN97}. The authors  constructed {\it nearly} adaptive estimation over the scale of generalized Sobolev classes.

The optimally adaptive estimator over the scale of Besov classes  was built in \cite{Jud97}. The estimation procedure is the hard threshold wavelet construction with random thresholds those choice are  based on some modification of the  comparison scheme proposed in \cite{lepski91}. Several years later
similar result was obtained in \cite{Silver}. The estimation method is  again hard threshold wavelet estimator but with
empirical bayes selection of thresholds. Both results were obtained under additional condition  $\beta>1/r$ which is slightly stronger than (\ref{eq1:class-param}).
{\it Efficient} adaptive estimator over the scale of Besov classes under $\bL_2$-loss was constructed in \cite{CH05} by use of empirical bayes thresholding.

We finish this part with mentioning  the papers \cite{Juditsky}, \cite{patricia}, where very interesting phenomena related to the adaptive density estimation under $\bL_p$-loss with unbounded support were observed, and   the paper \cite{Gol}, where $\bL_p$-aggregation of estimators
was proposed.


\paragraph{ Multivariate function estimation}

Much less is known when adaptive estimation of multivariate function is considered. The principal difficulty is related to the fact that the methods developed in the univariate case cannot be directly generalized to the multivariate setting.

In the series of papers  starting in the end of  70's Ibragimov and Hasminskii studied the problem of minimax  estimation over $\bN_{\vec{r},d}\big(\vec{\beta},\vec{L}\big)$ under  $\bL_p$--losses in different statistical models, see \cite{Has-Ibr} and references therein.
Note, however, that these authors treated only the case $r_i=p, i=1,\ldots,d$,  that allowed to prove that standard linear estimators (kernel, local polynomial, etc) are minimax. The optimally adaptive estimator corresponding to the latter case was constructed in \cite{GL11} in the density model on $\bR^d$.

The estimation over {\it isotropic} Besov class $\bB^{\vec{\beta}}_{\vec{r},q}\big(\vec{L}\big)$ was studied  in \cite{delyon-iod}, where the authors established the asymptotic of minimax risk under $\bL_p$-loss and constructed minimax estimators. Here the isotropy means that
$\vec{\beta}=(b,\ldots,b)$, $\vec{r}=(r,\ldots,r)$ and $\vec{L}=(L,\ldots,L)$. The asymptotics of minimax risk is given by
$$
\phi\big(\bB^{\vec{\beta}}_{\vec{r},q}\big(\vec{L}\big)\big)\sim
\left\{
\begin{array}{ll}
\e^{\frac{2b}{2b+d}},\quad& (2b+d)r>dp;
\\
\big(\e^2|\ln(\e)|\big)^{\frac{b-d/r+d/p}{2b-2d/r+d}},\quad& (2b+d)r\leq dp.
\end{array}
\right.
$$

Nearly adaptive with respect to  $\bL_p$-risk, $1\leq p<\infty$, estimator  over collection of  {\it isotropic} Nikolskii classes  $\bN_{\vec{r},d}\big(\vec{\beta},\vec{L}\big)\bowtie\bB^{\vec{\beta}}_{\vec{r},q}\big(\vec{L}\big)$  was built in \cite{GL08}. The proposed procedure is based on the special algorithm of  local bandwidth selection from the family of kernel estimators. The corresponding upper bound for maximal risks is proved  under additional assumption $b>d/r$.

\cite{bertin2} considered the problem of adaptive estimation over the scale of anisotropic H\"older classes, i.e. $\bN_{\vec{r},d}\big(\vec{\beta},\vec{L}\big)$ with
$r_i=\infty$ for any $i=1,\ldots, d$ under $\bL_\infty$-loss. The asymptotics of minimax risk is given here by
$$
\phi\big(\bN_{\vec{\infty},d}\big(\vec{\beta},\vec{L}\big)\big)\sim \big(\e^2|\ln(\e)|\big)^{\frac{\beta}{2\beta+1}},
$$
where $1/\beta=1/\beta_1+\cdots+1/\beta_d$. The construction of the optimally adaptive estimator is based on the selection rule from the family of kernel estimators developed in \cite{LepLev98}.

\cite{akakpo} studied the problem of adaptive estimation over the scale of anisotropic Besov classes $\bB^{\vec{\beta}}_{\vec{r},q}\big(\vec{L}\big)$
under $\bL_2$-loss in multivariate density model on the unit cube.
The construction of the optimally-adaptive estimator
is based on model selection approach and it uses  sophisticated approximation bounds. Note however that all results are proved in the situation  where
 coordinates of the vector $\vec{r}$ are the same ($r_i=r,\; i=1,\ldots d$).

\smallskip

For the first time the minimax and minimax adaptive estimation
over the scale of anisotropic  classes $\bN_{\vec{r},d}(\vec{\beta},\vec{L})$
under $\bL_p$-loss in the multivariate  model (\ref{eq:WGN-model}) was studied in full generality in \cite{lepski-kerk, lepski-kerk-08}.

To describe the results obtained in this paper we will need the following notations used in the sequel as well. Set
$
\omega^{-1}= (\beta_1r_1)^{-1}+\cdots+(\beta_dr_d)^{-1}
$
and define for any $1\leq s\leq\infty$
\begin{eqnarray*}
&& \tau(s)=1-1/\omega+1/(s\beta),\qquad  \kappa(s)=\omega(2+1/\beta)-s.
\end{eqnarray*}
In \cite{lepski-kerk} under assumption
\begin{equation}
\label{eq2:class-param}
\tau(\infty)>0,\qquad \sum_{i=1}^d\big[1/(r_i\beta_i)-1/(p\beta_i)\big]_+<2/p
\end{equation}
(called by the authors {\it the dense zone}) the following asymptotics of minimax risk was found
$$
\phi_\e\big(\bN_{\vec{r},d}\big(\vec{\beta},\vec{L}\big)\big)\sim
\e^{\frac{\beta}{2\beta+1}}.
$$
In \cite{lepski-kerk-08} under assumption
\begin{equation}
\label{eq3:class-param}
\tau(\infty)>0,\qquad \kappa(p)\leq 0,\qquad \vec{r}\in [1,p]^d,
\end{equation}
(called by the authors {\it the sparse zone}) the following asymptotics of minimax risk was found
$$
\phi_\e\big(\bN_{\vec{r},d}\big(\vec{\beta},\vec{L}\big)\big)\sim
\big(\e^2|\ln(\e)|\big)^{\frac{\tau(p)}{2\tau(2)}}.
$$
The authors built {\it nearly} adaptive with respect to  $\bL_p$-risk, $1\leq p<\infty$, estimator.
Its construction is based on the pointwise  bandwidths selection rule  which differs from whose
presented in \cite{LepLev98} as well as from the construction developed several years  later in \cite{GL08, GL13}. It is important to emphasize that
the method developed in the present paper is in some sense a "global" version of the aforementioned procedure.

The existence of an optimally-adaptive estimator as well as the asymptotics of minimax risk in the case, where assumptions (\ref{eq2:class-param}) and
(\ref{eq3:class-param}) are not fulfilled,  remained an open problem. Note also that  the assumption (\ref{eq1:class-param}) appeared in the univariate case can be rewritten as  $\tau(p)>0$. The minimax as well as adaptive estimation in the case $\tau(p)\leq 0$ was not investigated. One can suppose that a uniformly consistent estimator on $\bN_{r,1}\big(\beta,L\big)$ does not exist if  $\tau(p)\leq 0$ since $\tau(p)>0$ is the sufficient condition for the compact embedding of the univariate Nikolskii space into $\bL_p$, see \cite{Nikolski}.

\smallskip

The attempt to shed light on aforementioned problems was recently undertaken in \cite{GL13} in the framework of the density estimation  on $\bR^d$.
The authors are interested in adaptive estimation under $\bL_p$-loss, $p\in [1,\infty)$ over the collection of functional classes
$$
\bF_\vartheta=\bN_{\vec{r},d}\big(\vec{\beta},\vec{L}, M\big) \;:=\; \bN_{\vec{r},d}\big(\vec{\beta},\vec{L}\big)
\;\cap \; \left\{f:\|f\|_\infty\leq M\right\},\quad \vartheta=\big(\vec{\beta},\vec{r},\vec{L}, M\big).
$$
Adapting the results obtained in the latter paper to the observation  model (\ref{eq:WGN-model}) we first state that the asymptotics of the minimax risk satisfies
$$
\phi_\e\big(\bN_{\vec{r},d}\big(\vec{\beta},\vec{L},M\big)\big)\;\gtrsim\; \mu_\e^{\nu}
$$
where
\begin{eqnarray*}
\nu=\left\{
\begin{array}{ccc}
\frac{\beta}{2\beta+1},\quad & \kappa(p)>0;
\\*[2mm]
\frac{\tau(p)}{2\tau(2)},\quad & \kappa(p)\leq 0,\;\tau(\infty)> 0;
\\*[2mm]
\frac{\omega}{p},\quad & \kappa(p)\leq 0,\;\tau(\infty)\leq 0;
\end{array}
\right.
\quad
\mu_\e=\left\{
\begin{array}{ll}
\e^{2}, \quad  \kappa(p)> 0\;\; \text{or}\; &\kappa(p)\leq 0,\;\tau(\infty)\leq 0;
\\*[2mm]
\e^{2}|\ln(\e)|,\quad & \kappa(p)\leq 0,\;\tau(\infty)> 0.
\end{array}
\right.
\end{eqnarray*}
It is important to note that the obtained lower bound remains true if $p=\infty$ that implies in particular  that under $\bL_\infty$-loss there is no a uniformly consistent on
$\bN_{\vec{r},d}\big(\vec{\beta},\vec{L},M\big)$ estimator if $\tau(\infty)\leq 0$ (note that $\kappa(\infty)=-\infty$).

The authors proposed {\it nearly} adaptive estimator, i.e. the estimator whose  maximal risk is proportional to
$
\big(\e^{2}|\ln(\e)|\big)^{\nu},
$
whatever the value of the nuisance parameter $\vartheta=\big(\vec{\beta},\vec{r},\vec{L}, M\big)$ and $p\in [1,\infty)$.

Thus, the existence of  optimally-adaptive estimators remains an open problem. Moreover, all discussed results are obtained under additional assumption that the underlying function is uniformly bounded. We will see that the situation change completely if this condition does not hold.
The optimally-adaptive estimator
over the scale of anisotropic Nikolskii classes under $\bL_\infty$-loss
was constructed in \cite{lepski13a} under assumption $\tau(\infty)>0$. Since $\tau(\infty)>0$ implies automatically that $\bN_{\vec{r},d}\big(\vec{\beta},\vec{L},M\big)=\bN_{\vec{r},d}\big(\vec{\beta},\vec{L}\big)$ for some $M$ completely determined by $\vec{L}$ the investigation under $\bL_\infty$-loss is finalized.

\smallskip

We would like to finish our short overview with mentioning  works where the adaptation is studied not only with respect to the smoothness properties of the underlying function but also with respect to some structural assumptions imposed on the statistical model.
\begin{itemize}
\item
Composite function structure, \cite{MamHor}, \cite{ioud}, \cite{B-B};
\item
Multi-index structure (single-index, projection pursuit etc), \cite{HJPS}, \cite{GL09}, \cite{LS14};
\item
Multiple index model in density estimation, \cite{samarov}.
\item
Independence structure in density estimation, \cite{lepski13a}.
\end{itemize}
The problems of adaptive estimation  over the scale of functional classes defined on some manifolds were studied \cite{kerk11}, \cite{kerk12}.

\subsection{{Objectives}} Considering the collection of functional classes
$$
\bF_\vartheta=\bN_{\vec{r},d}\big(\vec{\beta},\vec{L}\big),\;\vartheta=\big(\vec{\beta},\vec{r},\vec{L}\big),
$$
we want to answer on the following questions
\begin{enumerate}
\item
{\it What is the optimal decay of the minimax risk for any fixed value of the nuisance parameter $\vartheta$ and norm index $p\in [1,\infty]$?}
\item
{\it Do optimally-adaptive estimators always exist?}
\end{enumerate}
To realize this program we propose first a new data-driven selection rule from the family of kernel estimators with {\it varying bandwidths} and establish for it so-called $\bL_p$-norm
oracle inequality. Then, we use this inequality in order to prove the adaptivity properties of the proposed estimation procedure.

Let us discuss our approach more in detail.
Throughout of the paper we will use the following notations. For any $u,v\in\bR^d$ the operations and relations $u/v$, $uv$, $u\vee v$,$u\wedge v$,
$u<v$, $au, a\in\bR,$ are understood in coordinate-wise sense and $|u|$ stands for euclidian norm of $u$.
All integrals are taken over $\bR^d$
unless the domain of integration is specified explicitly. For any real $a$ its  positive part is denoted by  $(a)_+$ and $\lfloor a\rfloor$ is used for its integer part.

\paragraph{Kernel estimator with varying bandwidth}

Put $\mH=\{\mh_s=e^{-s-2}, s\in\bN\}$ and
denote by  $\mathfrak{S}_1$ the set of all measurable functions defined  on $(-b,b)^d$ and  taking values in $\mH$.
Introduce
$$
\mathfrak{S}_d=\Big\{\vec{h}:(-b,b)^{d}\to\mH^d:\quad \vec{h}(x)=\big(h_1(x),\ldots,h_d(x)\big),\;x\in(-b,b)^d,\;\;h_i\in\mS_1,\;i=\overline{1,d}\Big\}.
$$
Let $K:\bR^d\to\bR$ be a function satisfying $\int K=1$. With any
$\vec{h}\in\mathfrak{S}_d$ we associate the function
$$
K_{\vec{h}}(t,x)=V^{-1}_{\vec{h}}(x)K\bigg(\frac{t-x}{\vec{h}(x)}\bigg),\;\; t\in\bR^d,\; x\in(-b,b)^d,
$$
where $V_{\vec{h}}(x)=\prod_{i=1}^dh_i(x)$.
Let $\mS^*$ be a given subset of $\mathfrak{S}_d$. Consider the family of estimators
\begin{equation}
\label{eq:family-of-estimators}
\cF\left(\mathfrak{S}^*\right)=\left\{\widehat{f}_{\vec{h}}(x)=X_\e\left(K_{\vec{h}}(\cdot,x)\right),\;\;\vec{h}\in\mS^*,\;x\in(-b,b)^d\right\}.
\end{equation}
We will call these estimators \textsf{kernel estimators with varying bandwidth}. This type of estimators was introduced in \cite{M87} in the context of
cross-validation technique.

We will be particulary interested in the set $\mS^*=\mS_d^{\text{const}}\subset\mathfrak{S}_d$  which consists of constant functions.
 Note that if   $\vec{h}\in \mS_d^{\text{const}}$  we come to the standard definition of kernel estimator in white gaussian noise model.

In view of (\ref{eq:WGN-model}) we have the following decomposition which will be useful in the sequel
\begin{equation}
\label{eq:decomposition-of-kernel-estimator}
\widehat{f}_{\vec{h}}(x)-f(x)=\int K_{\vec{h}}(t,x)\big[f(t)-f(x)\big]\nu_d(\rd t) +\ve\xi_{\vec{h}}(x),\quad \xi_{\vec{h}}(x)=\int K_{\vec{h}}(t,x)W(\rd t).
\end{equation}
We note that $\xi_{\vec{h}}$ is centered  gaussian random field on $(-b,b)^d$ with the covariance function
$$
V^{-1}_{\vec{h}}(x)V^{-1}_{\vec{h}}(y)\int K\bigg(\frac{t-x}{\vec{h}(x)}\bigg)K\bigg(\frac{t-y}{\vec{h}(y)}\bigg)\nu_d(\rd t),\quad x,y\in(-b,b)^d.
$$

\paragraph{Oracle approach} Our goal is to propose data-driven (based on $X^\e$) selection procedure from the collection $\cF\left(\mS^*\right)$ and establish for it
 $\bL_p$-norm oracle inequality. More precisely we construct the random field $\big(\vec{\mathbf{h}}(x),\;x\in(-b,b)^d\big)$ completely determined by the observation  $X^\e$, such that  $x\mapsto \vec{\mathbf{h}}(x)$ belongs to $\mS^*$, and prove that for any  $p\in [1,\infty]$, $q\geq 1$ and $\e>0$ small enough
\begin{equation}
\label{eq:L_p-oracle-inequality}
\cR^{(p)}_\e\big[\hat{f}_{\vec{\mathbf{h}}}; f\big]\leq \Upsilon_1\inf_{\vec{h}\in\mS^*}A_{p,q}^{(\e)}\left(f,\vec{h}\right)+\Upsilon_2\e.
\end{equation}
Here $\Upsilon_1$ and $\Upsilon_2$ are numerical constants depending on $d,p,q,b$ and $K$ only and the inequality (\ref{eq:L_p-oracle-inequality})
is established for any function $f\in\bL_p\big(\bR^d,\nu_d\big)\cap\bL_2\big(\bR^d,\nu_d\big)$. We call (\ref{eq:L_p-oracle-inequality})
$\bL_p$-{\it norm oracle inequality}.

We provide with explicit expression of the functional
$A_{p,q}^{(\e)}(\cdot,\cdot)$ that allows us to derive different  minimax adaptive results from the unique $\bL_p$-norm oracle inequality. In this context
it is interesting to note that in the "extreme cases" $p=1$ and $p=\infty$ it suffices to select the estimator from the family
$\cF\left(\mathfrak{S}_d^{\text{const}}\right)$. When $p\in (1,\infty)$, the oracle inequality (\ref{eq:L_p-oracle-inequality}) as well as the selection from the family $\cF\left(\mS^*\right)$ will be done for some special
choice of bandwidth's set $\mS^*$. We will see that the restrictions imposed on $\mS^*$ are rather weak that will allow us
to prove very strong adaptive results presented in Section \ref{sec:adaptive-estimation}.

\subsection{{Organization of the paper}} In Section \ref{sec:selection-rule-and-L_p-norm-oracle-inequality} we present our selection rule
and formulate  for it $\bL_p$-norm oracle inequality, Theorem \ref{th:L_p-norm-oracle-inequality}. Its consequence related to the selection from the family
$\mS_d^{\text{const}}$ is established in Corollary \ref{cor:th:L_p-norm-oracle-inequality}.
Section \ref{sec:adaptive-estimation} is devoted to adaptive estimation over the collection of anisotropic Nikolskii classes. Lower bound result is formulated in Theorem \ref{th:no-consitency} and the adaptive upper bound is presented in Theorem \ref{th:adaptive-upper-bound}.
In Section \ref{sec:open-problems} we discuss open problems in adaptive minimax estimation in different statistical models.
Proofs of main results are given in Sections
\ref{sec:proofs}--\ref{sec:proof-th:adaptive-upper-bound} and all technical lemmas are proven in Appendix.

\medskip

\section{Selection rule and $\bL_p$-norm oracle inequality}
\label{sec:selection-rule-and-L_p-norm-oracle-inequality}

\subsection{Functional classes of bandwidths}


Put for any $\vec{h}\in\mathfrak{S}_d$ and any  $\mathbf{s}=(s_1,\ldots,s_d)\in\bN^d$
$$
\Lambda_{\mathbf{s}}\big[\vec{h}\big]=\cap_{j=1}^{d}\Lambda_{s_j}\big[h_j\big],\qquad \Lambda_{s_j}\big[h_j\big]=\big\{x\in(-b,b)^d:\;\; h_j(x)=\mh_{s_j}\big\}.
$$
Let $\varkappa\in (0,1)$ and $\mL>0$ be given constants. Define
\begin{eqnarray*}
\bH_d(\varkappa,\mL)=\bigg\{\vec{h}\in\mathfrak{S}_d:\;\; \sum_{\mathbf{s}\in\bN^d}\nu_d^{\varkappa}\Big(\Lambda_{\mathbf{s}}\big[\vec{h}\big]\Big)\leq \mL\bigg\}.
\end{eqnarray*}
We remark that obviously $\mathfrak{S}_d^{\text{const}}\subset\bH_d(\varkappa,\mL)$ for any $\varkappa\in (0,1)$ and $\mL=(2b)^\varkappa$.

\smallskip

Put $\bN^*_p=\big\{\lfloor p\rfloor+1,\lfloor p\rfloor+2,\ldots\big\}$ and define for any $\mathcal{A}\geq e^{d}$
$$
\mathbb{B}(\cA)=\bigcup_{r\in \bN^*_p}\mathbb{B}_r(\mathcal{A}),\qquad \mathbb{B}_r(\mathcal{A})=\Big\{\vec{h}\in\mathfrak{S}_d:\;\;
\Big\| V^{-\frac{1}{2}}_{\vec{h}}\Big\|_{\frac{rp}{r-p}}\leq \cA\Big\}.
$$
Later on in the case $p\in (1,\infty)$ we will be interested in selection from the family $\cF\left(\bH\right)$, where $\bH$ is an arbitrary subset
of $\bH_d(\kappa,\mL,\cA):=\bH_d(\kappa,\mL)\cap\mathbb{B}(\cA)$, $\kappa\in (0,1/d),$ with some special choice $\cA=\cA_\e\to\infty, \e\to 0$.

The following notations related to the functional class $\mathbb{B}(\cA)$ will be exploited in the sequel. For any $\vec{h}\in\mathbb{B}(\cA)$ define
\begin{equation}
\label{eq:def-N_p(h)}
\bN^*_p\big(\vec{h},\cA\big)= \bN^*_p\cap \big[r_{\cA}(\vec{h}),\infty\big),
\qquad
r_{\cA}(\vec{h})=\inf\big\{r\in \bN^*_p:\;\; \vec{h}\in\mathbb{B}_r(\mathcal{A})\big\}.
\end{equation}
Obviously   $r_{\cA}\big(\vec{h}\big)<\infty$ for any $\vec{h}\in\mathbb{B}(\mathcal{A})$.

\paragraph{Assumptions imposed on the kernel $K$}
Let $a\geq 1$ and $A>0$ be fixed.
\begin{assumption}
\label{ass:kernel}
 There exists
  $\cK:\bR\to\bR$  such that $\int \cK=1$, $\text{supp}(K)\subset[-a,a]$ and
\begin{eqnarray*}
&(\mathbf{i})&\qquad|\cK(s)-\cK(t)|\leq A|s-t|,\;\;\forall s,t\in\bR;
\\
&(\mathbf{ii})&\qquad K(x)=\prod_{i=1}^d\cK(x_i),\;\;\forall x=(x_1,\ldots,x_d)\in\bR^d
\end{eqnarray*}

\end{assumption}
Throughout the paper we will consider only kernel estimators with $K$ satisfying Assumption \ref{ass:kernel}.

\subsection{{Upper functions and the choice of parameters}}
Put
\begin{equation}
\label{eq1:choice-of-parameters}
\mh_\e:=e^{-\sqrt{|\ln(\e)|}}, \quad
\cA_\e:=e^{\ln^2(\e)},
\end{equation}
and let $\mS_d(\mh_\e)\subset\mS_d$ consists of the functions $\vec{h}$ taking  values in $\mH^d(\mh_\e):=\mH^d\cap (0,\mh_\e]^d$.

\smallskip

Set $ C_2(r)=C_2\big(r,d\kappa, (2\mL)^d\big)$ and define for any  $\vec{h}\in\bB(\cA_\e)$
\begin{eqnarray*}
\label{eq:def-of-Psi_lambda}
\widetilde{\Psi}_{\e,p}\big(\vec{h}\big)&=&C_1\Big\|\sqrt{\big|\ln{\big(\e V_{\vec{h}}\big)}\big|} V^{-\frac{1}{2}}_{\vec{h}}\Big\|_p,\qquad\; p\in [1,\infty]
\\
\overline{\Psi}_{\e,p}\big(\vec{h}\big)&=&\inf_{r\in\bN^*_p(\vec{h},\cA_\e)}C_2(r)\Big\| V^{-\frac{1}{2}}_{\vec{h}}\Big\|_{\frac{rp}{r-p}},\quad p\in[1,\infty).
\end{eqnarray*}
Introduce finally
\begin{eqnarray*}
\Psi_{\e,p}\big(\vec{h}\big)=\left\{
\begin{array}{ll}
\widetilde{\Psi}_{\e,p}\big(\vec{h}\big)\wedge\overline{\Psi}_{\e,p}\big(\vec{h}\big),\quad &\vec{h}\in\bB(\cA_\e)\cap\mS_d(\mh_\e),\;\;p\in [1,\infty);
\\*[2mm]
\widetilde{\Psi}_{\e,p}\big(\vec{h}\big),\quad &\vec{h}\in\bB(\cA_\e)\setminus\mS_d(\mh_\e),\;\;p\in [1,\infty].
\end{array}
\right.
\end{eqnarray*}
Some remarks are in order.

1) The constant $C_1$  depends  on  $\cK$, $d, p$ and $b$ and its explicit expression is given in Section \ref{sec:constants}. The explicit expression of the quantity
$C_2(r,\tau,\cL), r>p,\;\tau\in (0,1),\;\cL>0,$ can be found in \cite{lepski13b}, Section 3.2.2. Its definition is rather involved and since it will not be exploited in the sequel we omit the definition of the latter quantity  in the present  paper. Here we only mention that $C_2(\cdot,\tau,\cL):(p,\infty)\to\bR_+$ is bounded on each bounded interval.
However  $C_2(r,\tau,\cL)\to\infty,\;r\to\infty$.

\smallskip

2) The selection rule presented below exploits  heavily  the fact that $\big\{\Psi_{\e,p}\big(\vec{h}\big), \vec{h}\in\bH\big\}, p\in [1,\infty],$ is  the upper function   for the collection   $\big\{\|\xi_{\vec{h}}\|_p, \vec{h}\in\bH\big\}$. Here the random field $\xi_{\vec{h}}$ appeared in the decomposition $(\ref{eq:decomposition-of-kernel-estimator})$ of kernel estimator and $\bH$ is an arbitrary countable subset of $\bH_d(\kappa,\mL,\cA_\e)$. The latter result is recently proved in \cite{lepski13b}, and it is presented in Proposition \ref{prop:upper-function-for-l_p-oracle}, Section \ref{sec:Proof of_th:L_p-norm-oracle-inequality}.

\smallskip

3) The choice of $\mh_\e$ and $\cA_\e$ is mostly dictated by the following simple observation which will be used for proving  adaptive results presented in Section \ref{sec:adaptive-estimation}.
\begin{equation}
\label{eq2:choice-of-parameters}
\lim_{\e\to 0} \e^{-a}\mh_\e=\infty,\qquad \lim_{\e\to 0}\e^{a}\cA_\e=\infty,\;\;\forall a>0.
\end{equation}
The general relation between parameters $\mh_\e$ and $\cA_\e$ can be found in \cite{lepski13b}.

\subsection{{Selection rule}}
\label{sec:subsection-selection-rule}

Let  $\bH$ be  a countable subset of $\bH_d\big(\kappa,\mL,\cA_\e\big)$. Define
\begin{equation}
\label{eq:selection-rule1}
\widehat{\cR}_{\bH}\big(\vec{h}\big)=\sup_{\vec{\eta}\in\bH}\bigg[
\left\|\widehat{f}_{\vec{h}\vee\vec{\eta}}-\widehat{f}_{\vec{\eta}}\right\|_p-\e\Psi_{\e,p}\big(\vec{h}\vee\vec{\eta}\big)-
\e\Psi_{\e,p}\big(\vec{\eta}\big)
\bigg]_+,\quad \vec{h}\in\bH.
\end{equation}
Our selection rule is given now by
$
\vec{\mathbf{h}}_0=\arg\inf_{\vec{h}\in\bH}\left\{\widehat{\cR}_{\bH}\big(\vec{h}\big)+\e\Psi_{\e,p}\big(\vec{h}\big)\right\}.
$
Since $\vec{\mathbf{h}}_0$ does not necessarily belong to $\bH$ we define finally $\vec{\mathbf{h}}\in\bH$ from the relation
\begin{equation}
\label{eq:selection-rule2}
\widehat{\cR}_{\bH}\big(\vec{\mathbf{h}}\big)+\e\Psi_{\e,p}\big(\vec{\mathbf{h}}\big)\leq \widehat{\cR}_{\bH}\big(\vec{\mathbf{h}}_0\big)+\e\Psi_{\e,p}\big(\vec{\mathbf{h}}_0\big)+\e,
\end{equation}
that leads to the estimator $\widehat{f}_{\vec{\mathbf{h}}}$.

\begin{remark}
\label{rem:why-finite-set}
We restrict ourselves  by consideration of countable subsets of $\bH_d\big(\kappa,\mL,\cA_\e\big)$ in order not to discuss the measurability
of $\widehat{f}_{\vec{\mathbf{h}}}$. Formally, the proposed selection rule can be applied for any
$\bH\subseteq\bH_d\big(\kappa,\mL,\cA_\e\big)$ for which final estimator can be correctly defined.
\end{remark}

\subsection{$\bL_p$-{norm oracle inequality}}
\label{sec:subsection-L_p-norm-oracle-inequality}

For any $\vec{\mathrm{h}}\in\mS_d$ define
$$
S_{\vec{\mathrm{h}}}(x,f)=\int_{\bR^d}K_{\vec{\mathrm{h}}}(t-x)f(t)\nu_d(\rd t),\;\;x\in \bR^d,
$$
  which is understood as kernel approximation (smoother) of the function $f$ at a point $x$.

For any $\vec{h},\vec{\eta}\in\mS_d$ introduce also
\begin{equation}
\label{eq:def-B_{vec{h}}}
B_{\vec{h},\vec{\eta}}(x,f):=\big|S_{\vec{h}\vee \vec{\eta}}(x,f)-S_{ \vec{\eta}}(x,f)\big|,\qquad B_{\vec{h}}(x,f)=\big|S_{\vec{h}}(x,f)-f(x)\big|,
\end{equation}
and define finally for any $p\in[1,\infty]$
\begin{equation}
\label{eq:def-bias-norm}
\cB^{(p)}_{\vec{h}}(f)=\sup_{\vec{\eta}\in\bH}\big\|B_{\vec{h},\vec{\eta}}(\cdot,f)\big\|_p+\big\|B_{\vec{h}}(\cdot,f)\big\|_p.
\end{equation}

\begin{theorem}
\label{th:L_p-norm-oracle-inequality}
Let Assumption \ref{ass:kernel} be fulfilled and let $p\in[1,\infty]$, $q\geq 1$, $\kappa\in (0,1/d)$ and $\mL\geq 1$
be fixed.
Then, there exists $\e(q)>0$ such that for any
 any $\e\leq \e(q)$ and $\bH\subseteq\bH_d\big(\kappa,\mL,\cA_\e\big)$
\begin{equation*}
\cR^{(p)}_\e\big[\hat{f}_{\vec{\mathbf{h}}}; f\big]\leq 5\inf_{\vec{h}\in\bH}\bigg\{ \cB^{(p)}_{\vec{h}}(f)+\e\Psi_{\e,p}\big(\vec{h}\big)\bigg\}+9(C_3+C_4+2)\e,\quad
\forall f\in\bL_p\big(\bR^d,\nu_d\big)\cap\bL_2\big(\bR^d,\nu_d\big).
\end{equation*}
\end{theorem}
The quantities  $C_3$ and $C_4$ depend on $\cK, p, q, b$ and $d$ only and their explicit expressions are presented in the Section \ref{sec:constants}.

\paragraph{Some consequences} The selection rule  (\ref{eq:selection-rule2}) deals with the family of kernel estimators with varying bandwidths. This allows, in particular, to apply $\bL_p$-norm oracle inequality established in Theorem \ref{th:L_p-norm-oracle-inequality}  to adaptive estimation over the collection of inhomogeneous and anisotropic  functional classes. However in some cases it suffices to select from much less "massive" set of bandwidths namely from $\mathfrak{S}_d^{\text{const}}$. In this case one can speak about standard multi-bandwidth selection. In particular, in the next section we will show that the selection from $\mathfrak{S}_d^{\text{const}}$ leads to optimally adaptive estimator over anisotropic Nikolskii classes if $p=\{1,\infty\}$. Moreover, considering $\mathfrak{S}_d^{\text{const}}$ we simplify considerably the "approximation error" $\cB^{(p)}_{\vec{h}}(f)$ as well as the upper function $\Psi_{\e,p}(\cdot)$.
The following corollary of Theorem \ref{th:L_p-norm-oracle-inequality} will be proved in Section \ref{sec:Proof of_th:L_p-norm-oracle-inequality}.

 Set $C_{2,p}=(2b)^{\frac{d}{p}}\inf_{r\in\bN_p^*}C_2(r)$ and define for any $\vec{h}\in\mathfrak{S}_d^{\text{const}}(\mh_\e):=\mathfrak{S}_d^{\text{const}}\cap\mS_d(\mh_\e)$
\begin{equation}
\label{eq:def-Psi_constant}
\Psi^{(\text{const})}_{\e,p}\big(\vec{h}\big)=C_{2,p} V^{-\frac{1}{2}}_{\vec{h}},\;\; p\in[1,\infty),\qquad\Psi^{(\text{const})}_{\e,\infty}\big(\vec{h}\big)=C_1\sqrt{\big|\ln{\big(\e V_{\vec{h}}\big)}\big|} V^{-\frac{1}{2}}_{\vec{h}}.
\end{equation}
Let $\big\{\mathbf{e}_1,\ldots,\mathbf{e}_d\big\}$ be the canonical basis in $\bR^d$.
For any $\vec{h}\in\mathfrak{S}_d^{\text{const}}$    introduce
\begin{equation}
\label{eq:def-small-b_h}
b_{\vec{h},j}(x)=\sup_{s:\:\mh_s\leq h_j}
\left|\int_{\bR}\cK(u)f\big(x+u\mh_s\mathbf{e}_j\big)\nu_1(\rd u)-f(x) \right|,\quad j=1,\ldots,d.
\end{equation}
Define finally   $\bH^{\text{const}}_\e=\mathfrak{S}_d^{\text{const}}(\mh_\e)\cap\big\{\vec{h}: V_{\vec{h}}\geq (2b)^{\frac{d}{p}}\cA_\e^{-2}\big\}$ and
let $\hat{f}^{(\text{const})}_{\vec{\mathbf{h}}}$ be the estimator obtained by the selection rule (\ref{eq:selection-rule2}) with
$\bH=\bH^{\text{const}}_\e$ and $\Psi_{\e,p}\big(\vec{h}\big)$ replaced by $\Psi^{(\text{const})}_{\e,p}\big(\vec{h}\big)$ given in (\ref{eq:def-Psi_constant}).

\begin{corollary}
\label{cor:th:L_p-norm-oracle-inequality}
Let Assumption \ref{ass:kernel} be fulfilled and let $p\in[1,\infty]$ and  $q\geq 1$
be fixed.
Then, there exists $\e(q)>0$ such that for
 any $\e\leq \e(q)$,  $\bH\subseteq\bH^{\text{const}}_\e$ and $f\in\bL_p\big(\bR^d,\nu_d\big)\cap\bL_2\big(\bR^d,\nu_d\big)$
\begin{equation*}
\cR^{(p)}_\e\big[\hat{f}^{(\text{const})}_{\vec{\mathbf{h}}}; f\big]\leq 5\inf_{\vec{h}\in\bH}\bigg\{3\|\cK\|^{d}_{1,\bR}\sum_{j=1}^d\big\|b_{\vec{h},j}
\big\|_p+\e\Psi^{(\text{const})}_{\e,p}\big(\vec{h}\big)\bigg\}+9(C_3+C_4+2)\e.
\end{equation*}

\end{corollary}
We remark that since $\bH^{\text{const}}_\e$ is finite  a selected multi-bandwidth $\vec{\mathbf{h}}\in\bH$ is given by
$$
\vec{\mathbf{h}}=\arg\inf_{\vec{h}\in\bH}\left\{\widehat{\cR}_{\bH}\big(\vec{h}\big)+\e\Psi^{(\text{const})}_{\e,p}\big(\vec{h}\big)\right\}.
$$

\section{Adaptive estimation}
\label{sec:adaptive-estimation}

In this section
we study
properties of the estimator defined in
Section \ref{sec:subsection-selection-rule}.
The $\bL_p$-norm oracle inequalities
obtained Theorem~\ref{th:L_p-norm-oracle-inequality} and Corollary \ref{cor:th:L_p-norm-oracle-inequality} can be viewed as initial step
in 
bounding $\bL_p$-risk of this estimator on the anisotropic Nikol'skii classes.

\subsection{Anisotropic Nikolskii classes}
\label{sec:nikolski}
Recall that   $(\mathbf{e}_1,\ldots,\mathbf{e}_d)$ denotes the canonical basis of $\bR^d$.
 For function $g:\bR^d\to \bR^1$ and
real number $u\in \bR$ define
{\em the first order difference operator with step size $u$ in direction of the variable
$x_j$}~by
\[
 \Delta_{u,j}g (x)=g(x+u\mathbf{e}_j)-g(x),\;\;\;j=1,\ldots,d.
\]
By induction,
the $k$-th order difference operator with step size $u$ in direction of the variable $x_j$ is
defined~as
\begin{equation}
\label{eq:Delta}
 \Delta_{u,j}^kg(x)= \Delta_{u,j} \Delta_{u,j}^{k-1} g(x) = \sum_{l=1}^k (-1)^{l+k}\binom{k}{l}\Delta_{ul,j}g(x).
\end{equation}
\begin{definition}
\label{def:nikolskii}
For given  vectors $\vec{r}=(r_1,\ldots,r_d)$, $r_j\in [1,\infty]$, $\vec{\beta}=(\beta_1,\ldots,\beta_d)$,
$\beta_j>0$, and $\vec{L}=(L_1,\ldots, L_d)$, $L_j>0$, $j=1,\ldots, d$, we
say that function $g:\bR^d\to \bR^1$ belongs to the anisotropic
Nikolskii class $\bar{\bN}_{\vec{r},d}\big(\vec{\beta},\vec{L}\big)$ if
\begin{itemize}
\item[{\rm (i)}] $\|g\|_{r_j,\bR^d}\leq L_{j}$ for all $j=1,\ldots,d$;
\item[{\rm (ii)}]
for every $j=1,\ldots,d$ there exists natural number  $k_j>\beta_j$ such that
\begin{equation}\label{eq:Nikolski}
 \Big\|\Delta_{u,j}^{k_j} g\Big\|_{r_j,\bR^d} \leq L_j |u|^{\beta_j},\;\;\;\;
\forall u\in \bR,\;\;\;\forall j=1,\ldots, d.
\end{equation}
\end{itemize}
\end{definition}
Recall that the consideration of white gaussian noise model requires $f\in\bL_2\big(\bR^d\big)$ that is not always guaranteed
by $f\in \bar{\bN}_{\vec{r},d}\big(\vec{\beta},\vec{L}\big)$. So, later on we will study the functional classes $\bN_{\vec{r},d}\big(\vec{\beta},\vec{L}\big)=
\bar{\bN}_{\vec{r},d}\big(\vec{\beta},\vec{L}\big)\cap\bL_2\big(\bR^d\big)$ which we will also call anisotropic Nikolskii classes. Some conditions guaranteed
$\bN_{\vec{r},d}\big(\vec{\beta},\vec{L}\big)=\bar{\bN}_{\vec{r},d}\big(\vec{\beta},\vec{L}\big)$ can be found in Section \ref{sec:subsec-inclusions}.

\subsection{{Main results}}
\label{sec:adap-nikol}
Let $\bN_{\vec{r},d}\big(\vec{\beta},\vec{L}\big)$ be the
anisotropic Nikolskii functional class.
%
Put
\[
 \frac{1}{\beta} := \sum_{j=1}^d \frac{1}{\beta_j},\;\;\;\;\;
\frac{1}{\omega} := \sum_{j=1}^d \frac{1}{\beta_jr_j}, \;\;\;\;\;L_\beta:= \prod_{j=1}^d L_j^{1/\beta_j},
\]
and define for any $1\leq s\leq\infty$
\begin{eqnarray*}
 \label{eq:kappa&tau}
 \tau(s)=1-1/\omega+1/(s\beta),\qquad  \kappa(s)=\omega(2+1/\beta)-s.
\end{eqnarray*}
The following obvious relation will be useful in the sequel.
\begin{equation}
\label{eq:relation-kappa-tau}
 \frac{\kappa(s)}{\omega s}=\frac{2-s}{s}+\tau(s).
\end{equation}
Set finally $p^*=\big[\max_{j=1,\ldots,d}r_l\big]\vee p$ and introduce
\begin{eqnarray*}
\mathfrak{a}&=&\left\{
\begin{array}{cccc}
\frac{\beta}{2\beta+1},\quad & \kappa(p)>0;
\\*[2mm]
\frac{\tau(p)}{2\tau(2)},\quad & \kappa(p)\leq 0,\;\tau(p^*)> 0;
\\*[2mm]
\frac{\omega(p^*-p)}{p(p^*-\omega(2+1/\beta)},\quad & \kappa(p)\leq 0,\;\tau(p^*)\leq 0,\;p^*>p;
\\*[2mm]
0,\quad & \kappa(p)\leq 0,\;\tau(p^*)\leq 0;\; p^*=p.
\end{array}
\right.
\\*[2mm]
\delta_\e&=&\left\{
\begin{array}{cc}
L_\beta\e^{2}, \quad & \kappa(p)> 0;
\\*[2mm]
L_\beta\e^{2}|\ln(\e)|,\quad &  \kappa(p)\leq 0,\;\tau(p^*)\leq 0;
\\*[2mm]
L_\beta^{\frac{1-2/p}{\tau(p)}}\e^{2}|\ln(\e)|,\quad & \kappa(p)\leq 0,\;\tau(p^*)> 0.
\end{array}
\right.
\end{eqnarray*}

\subsubsection{\textsf{Lower  bound of minimax risk}}

\begin{theorem}
\label{th:no-consitency}
Let  $q\geq 1$, $L_0>0$ and $1\leq p\leq\infty$ be fixed.
Then for any  $\vec{\beta}\in (0,\infty)^d,\; \vec{r}\in [1,\infty]^d$ and  $\vec{L}\in [L_0,\infty)^d$
there exists $c>0$ independent of $\vec{L}$ such that
$$
\liminf_{\e\to 0}\;\inf_{\tilde{f}_\e}\sup_{f\in\bN_{\vec{r},d}\big(\vec{\beta},\vec{L}\big)}\delta_\e^{-\mathfrak{a}}\cR^{(p)}_\e\big[\tilde{f}_\e; f\big]\geq c,
$$
where infimum is taken over all possible estimators.

\end{theorem}

Let us make several remarks.

$1^0.\;$ {\it Case $p^*=p$}. Taking into account  (\ref{eq:relation-kappa-tau}) we note that
there is no a uniformly consistent estimator over $\bN_{\vec{r},d}\big(\vec{\beta},\vec{L}\big)$ if
\begin{equation}
\label{eq:nonexistence-uniform-consistency}
\tau(p)\mathrm{1}_{[2,\infty)}(p)+\kappa(p)\mathrm{1}_{[1,2)}(p)\leq 0,
\end{equation}
and this result  seems to be new.
As it will follow from the next theorem the latter condition is necessary and sufficient for nonexistence of uniformly consistent estimators over  $\bN_{\vec{r},d}\big(\vec{\beta},\vec{L}\big)$ under $\bL_p$-loss, $1\leq p\leq\infty$. In the case of $\bL_\infty$-loss, (\ref{eq:nonexistence-uniform-consistency}) is reduced to $\omega\leq 1$ and the similar result was recently proved
in \cite{GL13} for the density model.

$2^0.\;$ {\it Case $\kappa(p)\leq 0,\;\tau(p^*)\leq 0,\;p^*>p$.} The lower bound for minimax risk given in this case by
$$
\big(L_\beta\e^{2}|\ln(\e)|\big)^{\frac{\omega(p^*-p)}{p(p^*-\omega(2+1/\beta)}}
$$
is new. It is  interesting that the latter case does not appear in the dimension 1 or, more generally, when isotropic Nikolskii classes are considered. Indeed, if $r_l=r$ for all $l=1,\ldots d,$ then $p^*>p$ means $r>p$ that, in its turn, implies $\tau(p^*)=\tau(r)=1>0$.
It is worth mentioning that we improve in order the lower  bound recently found in \cite{GL13}, which corresponds  formally to our case $p^*=\infty$.

$3^0.\;$ {\it Case $\kappa(p)\leq 0,\;\tau(p^*)> 0$.} For the first time the same result was proved  in \cite{lepski-kerk-08} but under more restrictive assumption $\kappa(p)\leq 0,\;\tau(\infty)> 0$. Moreover, the dependence of the asymptotics of the minimax risk on $\vec{L}$ was not optimal.

$4^0.\;$ {\it Case $\kappa(p)>0$.} Presented lower bound   of minimax risk became the statistical folklore since it is the minimax rate of convergence over
anisotropic H\"older class ($r_l=\infty,\; l=1,\ldots d,$). If so, the required result can be easily deduced from the embedding of a H\"older class to
$\bN_{\vec{r},d}\big(\vec{\beta},\vec{L}\big)$ whatever the value of $\vec{r}$. However the author was enable to find  exact references and derived
the announced result from the general construction used in the proof of Theorem \ref{th:no-consitency}. Moreover we are interested in finding not only the optimal decay of the minimax risk with respect to $\e\to 0$ but also its correct dependence of the radii $\vec{L}$.

\subsubsection{\textsf{Upper bound for minimax risk. Optimally-adaptive estimator}}

The results of this section will be derived from $\bL_p$-norm oracle inequalities proved in Theorem \ref{th:L_p-norm-oracle-inequality} and Corollary
\ref{cor:th:L_p-norm-oracle-inequality}.

\paragraph{Construction of kernel $K$}
We will use the following specific kernel $K$ [see, e.g., \cite{lepski-kerk} or \cite{GL11}] in the definition of the estimator's family (\ref{eq:family-of-estimators}).

\smallskip

 Let  $\ell$ be an integer number,
and let $w:[-1/(2\ell), 1/(2\ell)]\to \bR^1$ be a function satisfying  $\int w(y)\rd y=1$,
and $w\in\bC^{1}(\bR^1)$. Put
\begin{equation}
\label{eq:w-function}
 w_\ell(y)=\sum_{i=1}^\ell \binom{\ell}{i} (-1)^{i+1}\frac{1}{i}w\Big(\frac{y}{i}\Big),\qquad
 K(t)=\prod_{j=1}^d w_\ell(t_j),\;\;\;\;t=(t_1,\ldots,t_d).
\end{equation}

\paragraph{Set of bandwidths}

Set $t_{k,n}=-(b+1)+(b+1)k2^{1-n}$, $k=0,\ldots 2^n$, $n\in\bN^*$ and let $\Delta_{k,n}=\big[t_{k,n},t_{k+1,n}\big), k=0,\ldots 2^{n-1}$,
$\Delta_{k,n}=\big(t_{k,n},t_{k+1,n}\big], k=2^{n-1}+1,\ldots, 2^{n}-1$. Thus, $\{\Delta_{k,n}, k=0,\ldots 2^n\}$ forms the partition of $(-b-1,b+1)$
whatever  $n\in\bN^*$.

\smallskip

For any  $n\in\bN^*$ set also  $\mathfrak{K}_{n}=\{0,\ldots 2^{n}\}^d$ and define
\begin{equation}
\label{eq:def-partition-collection}
\Gamma_{d}(n)=\Big\{\Delta^{(d)}_{\mathbf{k},n}=\Delta_{k_1,n}\times\cdots\times\Delta_{k_d,n},\;\mathbf{k}=(k_1,\ldots, k_d)\in\mK_{n}\Big\}.
\end{equation}
For any  $n\in\bN^{*}$ the collection of cubs $\Gamma_{d}(n)$ determines  the partition of $(-b-1,b+1)^d$.

Denote by $\mS_d^{(n)}, n\in\bN^{*},$ the set of all step functions defined on $(-b,b)^d$ with the  steps belonging to $\Gamma_{d}(n)\cap (-b,b)^d$ and taking values in $\mH^d$.

Introduce finally for any $R>0$
$$
\bH_\e(R)=\bH_d\big(1/(2d),R,\cA_\e\big)\cap\big\{\cup_{n\in\bN^{*}}\mS^{(n)}_{d}\big\},
$$
where $\cA_\e$  is given in  (\ref{eq1:choice-of-parameters}).

\smallskip

Let $\widehat{f}^{(R)}_{\vec{\mathbf{h}}},\;R>0,$ denote the estimator obtained by the selection rule (\ref{eq:selection-rule1})--(\ref{eq:selection-rule2}) from the family
of kernel estimators $\cF\left(\bH_\e(R)\right)$ and $\hat{f}^{(\text{const})}_{\vec{\mathbf{h}}}$ denote the estimator constructed in Corollary \ref{cor:th:L_p-norm-oracle-inequality}. Both constructions are made with
 the kernel $K$ satisfying (\ref{eq:w-function}).

\paragraph{Adaptive upper bound}

For any  $\ell\in\bN^*$ and  $L_0>0$ set $\Theta=(0,\ell]^d\times[1,\infty]^d\times[L_0,\infty)^d$
and later on we will use the notation $\vartheta\in\Theta$ for the triplet  $\big(\vec{\beta},\vec{r},\vec{L}\big)$.
Denote $\cP=\Theta\times[1,\infty]$ and introduce
$$
\cP^{\text{consist}}=\big\{(\vartheta,p)\in\cP:
\tau(p)\mathrm{1}_{[2,\infty)}(p)+\kappa(p)\mathrm{1}_{[1,2)}(p)> 0\big\}\cup
\big\{(\vartheta,p)\in\cP:
 p^*>p\big\}.
$$
The latter set consists of the class parameters and norm indexes for which a uniform consistent estimation is possible.

Let $V_p\big(\vec{L}\big)$ be the  quantity whose presentation is postponed  to Section \ref{sec:subsec-proof of the theorem-finite-p} since its expression is rather cumbersome. Put $L^*=\min_{j: r_j=p^*} L_j$ and introduce
\begin{eqnarray*}
\bde_\e&=&\left\{
\begin{array}{lll}
L_\beta\e^{2}, \quad & \kappa(p)\geq 0;
\\*[2mm]
L_\beta(L^*)^{\frac{1}{\ma}}\e^{2}|\ln(\e)|, \quad &  \kappa(p)\leq 0,\;\tau(p^*)\leq 0;
\\*[2mm]
V_p\big(\vec{L}\big)\e^{2}|\ln(\e)|,\quad & \kappa(p)\leq 0,\;\tau(p^*)> 0.
\end{array}
\right.
\end{eqnarray*}

\begin{theorem}
\label{th:adaptive-upper-bound}
Let  $q\geq 1$, $L_0>0$ and  $\ell>0$  be fixed and let $R=3+\sqrt{2b}$.

1) For any $(\vartheta,p)\in\cP^{\text{consist}}$ such that  $p\in(1,\infty)$, $\vec{r}\in(1,\infty]^d$ and $\kappa(p)\neq0$
there exists  $C>0$ independent of $\vec{L}$ for which
$$
\limsup_{\e\to 0}\;
\sup_{f\in\bN_{\vec{r},d}\big(\vec{\beta},\vec{L}\big)}\bde_\e^{-\mathfrak{a}}\cR^{(p)}_\e\big[\widehat{f}^{(R)}_{\vec{\mathbf{h}}}; f\big]\leq C.
$$

2) For any $(\vartheta,p)\in\cP^{\text{consist}}$,  $p\in\{1,\infty\}$ there exists  $C>0$ independent of $\vec{L}$ for which
$$
\limsup_{\e\to 0}\;
\sup_{f\in\bN_{\vec{r},d}\big(\vec{\beta},\vec{L}\big)}\bde_\e^{-\mathfrak{a}}\cR^{(p)}_\e\big[\hat{f}^{(\text{const})}_{\vec{\mathbf{h}}}; f\big]\leq C.
$$

3) For any $(\vartheta,p)\in\cP^{\text{consist}}$ such that  $p\in(1,\infty)$, $\vec{r}\in(1,\infty]^d$ and $\kappa(p)=0$
there exists  $C>0$ independent of $\vec{L}$ for which
$$
\limsup_{\e\to 0}\;
\sup_{f\in\bN_{\vec{r},d}\big(\vec{\beta},\vec{L}\big)}\bde_\e^{-\mathfrak{a}}(|\ln(\e)|)^{\frac{1}{p}}\cR^{(p)}_\e\big[\widehat{f}^{(R)}_{\vec{\mathbf{h}}}; f\big]\leq C.
$$
\end{theorem}

Some remarks are in order.

$1^0.\;$ Combining the results of Theorems \ref{th:no-consitency} and \ref{th:adaptive-upper-bound} we  conclude that  optimally-adaptive
estimators under $\bL_p$-loss exist over all parameter set $\cP^{\text{consist}}$ if $p\in\{1,\infty\}$. If $p\in(1,\infty)$ such estimators exist as well
except the boundary cases  $\kappa(p)=0$ and $\min_{j=1,\ldots,d}r_j=1$.

$2^0.\;$ We remark  that the upper and lower bound for minimax risk differ each other on the boundary $\kappa(p)=0$ only by $(|\ln(\e)|)^{\frac{1}{p}}$-factor. Using $(1,1)$-weak type inequality for strong maximal operator,  \cite{Guzman}, one can prove
adaptive upper bound on the boundary $\min_{j=1,\ldots,d}r_j=1$ containing additional $(|\ln(\e)|)^{\frac{d-1}{p}}$-factor. Note, nevertheless,
that exact asymptotics of minimax risk on both boundaries  remains an open problem.

$3^0.\;$ We obtain  full classification of minimax rates over anisotropic Nikolskii classes if $p\in\{1,\infty)\}$ and "almost" full one
(except the boundaries mentioned above) if $p\in(1,\infty)$.
We can assert that $\delta^{\ma}_\e$ is minimax rate of convergence on $\bN_{\vec{r},d}\big(\vec{\beta},\vec{L}\big)$ for any
$\vec{\beta}\in (0,\infty)^d,\vec{r}\in(1,\infty]^d$ and $\vec{L}\in(0,\infty)^d$.  Indeed,  for given $\vec{\beta}$ and $\vec{L}$ one can choose
$L_0=\min_{j=1,\ldots d}L_j$ and the number $\ell$, used in the  kernel construction (\ref{eq:w-function}), as an any integer strictly larger than $\max_{j=1,\ldots d}\beta_j$.

$4^0.\;$ We remark that the dependence of minimax rate on $\vec{L}$ is correct ( $\delta_\e=\bde_\e$) if $\kappa(p)\geq 0$. In spite of the cumbersome expression of the quantity $V_p\big(\vec{L}\big)$ one can easily check that
$$
V_p\big(\vec{L}\big)=L_\beta^{\frac{1-2/p}{\tau(p)}}
$$
if $L_j=L$ for any $j=1,\ldots,d$. Hence, under this restriction $\delta_\e=\bde_\e$ if
$\kappa(p)\leq 0, \tau(p^*)>0$ as well.

\section{Open problems in adaptive estimation}
\label{sec:open-problems}

The goal of this section is to discuss the directions in adaptive multivariate function estimation to be developed. We do not pretend here to cover whole  specter of existing problems and mostly restrict ourselves by consideration of the  adaptation over the scale of anisotropic classes. Moreover we will be concentrated on principal difficulties and the mathematical aspect of the problem and we
will not pay much attention to the technical details and practical applications.
Although we will speak  about adaptive estimation it is important to realize   that for the majority of problems discussed below very little is known about the  minimax approach.

\subsection{{Abstract statistical model}}

 Let $\big(\cY^{(n)},\mA^{(n)},\bP_f^{(n)}, f\in\bF\big)$  be the sequence of statistical experiments generated by observation $Y^{(n)}, n\in\bN^*$.  Let $\Lambda$ be a set and
$\rho:\Lambda\times\Lambda\to\bR_+$ be a loss functional.
The goal is to estimate the mapping $G:\bF\to\Lambda$ and as an estimator we understand an $Y^{(n)}$-measurable $\Lambda$-valued map.

The quality of an estimation procedure $\tilde{G}_n$ on $\bF$ is measured by the maximal risk
$$
\cR_n\big[\tilde{G}_n; \bF\big]=\Big\{\sup_{f\in\bF}\bE^{(n)}_f\rho^q\big(\tilde{G}_n,G(f)\big)\Big\}^{\frac{1}{q}},\;\; q\geq 1.
$$
and as previously   $\phi_n(\bF)=\inf_{\tilde{G}_n}\cR_n\big[\tilde{G}_n; \bF\big]$ denotes the minimax risk.

Assume that $\bF\supset\cup_{\vartheta\in\Theta}\bF_\vartheta$, where $\{\bF_\vartheta, \vartheta\in\Theta\}$ is a given collection of sets.

\smallskip

\noindent \verb"PROBLEM I (fundamental)":\;\; \verb"Find" \verb"necessary" \verb"and" \verb"sufficient" \verb"conditions" \verb"of" \verb"existence" \verb"of"\; \verb"optimally"-\verb"adaptive" \verb"estimators", \verb"i.e." \verb"the existence of an estimator" $\hat{G}_n$ \verb"satisfying"
$$
\cR_n\big[\hat{G}_n; \bF_\vartheta\big]\sim \phi_n\big(\bF_\vartheta\big),\quad \forall\vartheta\in\Theta.
$$

It is well-known that optimally-adaptive estimators do not always exist, see \cite{lepski90}, \cite{lepski92a}, \cite{EfLow}, \cite{CLow05}. Hence, the goal is to understand how the answer on aforementioned question depends on the statistical model, underlying estimation problem (mapping $G$) and the collection of considered classes. The attempt to provide with such classification was undertaken in \cite{lepski92a}, but found there sufficient conditions of existence as well as of nonexistence of optimally-adaptive estimators are too restrictive.

It is important to realize that the answers on formulated problem may be different even if the statistical model and the collection of functional classes are the same and  estimation problems have  "similar nature". Indeed, consider univariate model (\ref{eq:WGN-model}) and let $\bF_\vartheta=\bN_{\infty,1}(\beta,L),\;\vartheta=(\beta,L),$ be the collection of H\"older classes. Set
$$
G_\infty(f)=\|f\|_\infty,\qquad G_2(f)=\|f\|_2.
$$
As we know the optimally-adaptive estimator of $f$, say $\hat{f}_n$,  under $\bL_\infty$-loss was constructed in \cite{lepski91}. Moreover  the asymptotics of minimax risk under $\bL_\infty$-loss on $\bN_{\infty,1}(\beta,L)$ coincides with asymptotics corresponding to the estimation of $G_\infty(\cdot)$.  Therefore, $\hat{G}_n:=G_\infty(\hat{f}_n)$ is an optimally-adaptive estimator for $G_\infty(\cdot)$. On the other hand, there is no optimally-adaptive estimator for $G_2(\cdot)$, \cite{CLow06}.

\subsection{{White gaussian noise model}} Let us return to the problems studied in the present paper. Looking at the optimally-adaptive estimator proposed in Theorem \ref{th:adaptive-upper-bound} we conclude that its construction is  not feasible. Indeed, it is based on the selection from very huge set of parameters, sometimes even infinite.

\smallskip

\noindent \verb"PROBLEM II (feasible estimator)":\; \verb"Find" \verb"optimally"-\verb"adaptive" \verb"estimator" \verb"whose" \verb"construction" \verb"would be"
\verb"computationally" \verb"reasonable".

At our glance the interest to this problem is not related to the "practical applications" since the pointwise bandwidths selection rule from \cite{GL13} will do this job although it is not
theoretically optimal. We think that "feasible solution" could bring new ideas and approaches to the construction of estimation procedures.

Another source of problems is structural adaptation. Let us consider one of the possible directions. Denote by $\cE$ the set of all $d\times s$ real matrices,
$1\leq s< d$.  Introduce the following collection of functional classes
$$
\bF_\vartheta=\bS_{\vec{r},d}\big(\vec{\beta},\vec{L},E\big):=\Big\{f:\bR^d\to\bR: \;\; f(x)=g(Ex),\;g\in\bN_{\vec{r},p}\big(\vec{\beta},\vec{L}\big),\;E\in\cE\Big\},\;\;\vartheta=\big(\vec{\beta},\vec{r},\vec{L},E\big).
$$


\noindent \verb"PROBLEM III (structural adaptation)":\; \verb"Prove" \verb"or" \verb"disprove" \verb"the" \verb"existence" \verb"of"  \verb"optimally"-\verb"adaptive" \verb"estimators" \verb"over" \verb"the" \verb"collection"
$\bS_{\vec{r},d}\big(\vec{\beta},\vec{L},E\big)$ \verb"under" $\bL_p$-\verb"loss".

\smallskip

Note that if $\vec{r}=(\infty,\ldots,\infty)$ (H\"older case) the optimally adaptive estimator was constructed in \cite{GL09}. {\it Nearly} adaptive estimator
in the case $s=1$ (single index constraint) and  $d=2$   was proposed in \cite{LS14}. Many other structural models like additive, projection pursuit or their generalization, see \cite{GL09}, can be studied as well.

\subsection{{Density model}} Let $X_i, i=1,\ldots,n,$ be {\it i.i.d.} $d$-dimensional random vectors with common probability density $f$.  The goal is to estimate $f$  under $\bL_p$-loss on $\bR^d$.

\smallskip

\noindent \verb"PROBLEM IV ":\; \verb"Prove" \verb"or" \verb"disprove" \verb"the" \verb"existence" \verb"of"  \verb"optimally"-\verb"adaptive" \verb"estimators" \verb"over" \verb"the" \verb"collection" \verb"of" \verb"anisotropic" \verb"Nikolskii" \verb"classes"
$\bN_{\vec{r},d}\big(\vec{\beta},\vec{L}\big)$ \verb"under" $\bL_p$-\verb"loss".

\smallskip

The last advances in this task  were made in \cite{GL13}. However, as it was conjectured in this paper, the developed there local approach cannot lead to the construction of optimally-adaptive estimators. On the other hand it is not clear how to adapt the approach developed in the present paper to the density estimation on $\bR^d$. Indeed, the key element of our procedure is  upper functions for $\bL_p$-norm of random fields found in \cite{lepski13b}. These results are heavily based on the fact that corresponding norm is defined on a bounded interval of $\bR^d$.

The same problem can be formulated for more complicated {\it deconvolution model}. Recent advances in the estimation under $\bL_2$-loss in this model can be found in  \cite{comte}.

\subsection{{Regression model}} Let $\xi_i, i\in\bN^*,$ be {\it i.i.d. symmetric} random variables with common probability density $\varrho$ and let $X_i, i\in\bN^*,$ be  {\it i.i.d.}  $d$-dimensional random vectors with common probability density $g$.
Suppose that we observe the pairs $(X_1,Y_1),\ldots, (X_n,Y_n)$ satisfying
$$
Y_i=f(X_i)+\xi_i,\quad i=1,\ldots,n.
$$
The goal is to estimate function $f$ under $\bL_p$-loss on $(-b,b)^d$, where $b>0$ is a given number.

We will suppose that the sequences $\xi_i, i\in\bN^*,$ and $X_i, i\in\bN^*,$ are mutually independent and that the design density $g$ (known or unknown) is separated away from zero on  $(-b,b)^d$.

\paragraph{Regular noise}
Suppose  that there exists $a>0$ and $A>0$ such that for any $u,v\in [-a,a]$
\begin{equation}
\label{eq:reg-noise}
\int_{\mathbb{R}} \frac{\varrho(y+u)\varrho(y+v)}{\varrho(y)} \rd y\leq 1+A|uv|.
\end{equation}
Assume also that $\bE|\xi_1|^\alpha<\infty$ for some $\alpha\geq 2$.

\smallskip

\noindent \verb"PROBLEM V:"\; \verb"Prove" \verb"or" \verb"disprove" \verb"the" \verb"existence" \verb"of"  \verb"optimally"-\verb"adaptive" \verb"estimators" \verb"over" \verb"the" \verb"collection" \verb"of" \verb"anisotropic" \verb"Nikolskii" \verb"classes"
$\bN_{\vec{r},d}\big(\vec{\beta},\vec{L}\big)$ \verb"under" $\bL_p$-\verb"loss".

\smallskip

The interesting question arising in this context is what is the minimal value of $\alpha$ under which the formulated problem can be solved. In particular, is it related to the norm index $p$ or not?

\paragraph{Cauchy noise} Let $\varrho(x)=\big\{\pi(1+x^2)\big\}^{-1}$. In this case the noise is of course regular, i.e. (\ref{eq:reg-noise}) holds, but the moment assumption fails. At our knowledge there is no minimax and minimax adaptive results in the multivariate regression model with  noise "without moments" when anisotropic functional classes are considered.

\smallskip

\noindent \verb"PROBLEM VI:"\; \verb"Propose" \verb"the" \verb"construction"  \verb"of"  \verb"optimally"-\verb"adaptive" \verb"estimators" \verb"over" \verb"the" \verb"scale" \verb"of" \verb"anisotropic" \verb"H"\"o\verb"lder" \verb"classes"
$\bN_{\vec{\infty},d}\big(\vec{\beta},\vec{L}\big)$ \verb"under" $\bL_p$-\verb"loss".

\smallskip

The same problem can be of course formulated over the scale of anisotropic Nikolskii classes but it seems that nowadays neither probabilistic nor the tools from functional analysis  are
sufficiently developed in order to proceed to this task.

\paragraph{Irregular noise} Consider two particular examples:
$$
\varrho(x)=2^{-1}\mathrm{1}_{[-1,1]}(x),\qquad \varrho_\gamma(x)=C_\gamma e^{-|x|^{\gamma}},\;\; \gamma\in (0,1/2).
$$
In both cases the condition (\ref{eq:reg-noise}) is not fulfilled. In parametric case $f(\cdot)\equiv \mathbf{f}\in\bR$ the minimax rate of convergence is faster than  $n^{-\frac{1}{2}}$, \cite{IH81}.

\smallskip

\noindent \verb"PROBLEM VII:"\; \verb"Find" \verb"minimax" \verb"rate" \verb"of" \verb"convergence" \verb"on"   \verb"anisotropic" \verb"H"\"o\verb"lder" \verb"class"
$\bN_{\vec{\infty},d}\big(\vec{\beta},\vec{L}\big)$ \verb"under" $\bL_p$-\verb"loss". \verb"Propose" \verb"an" \verb"aggregation" \verb"scheme" \verb"for" \verb"these" \verb"estimators" \verb"led" \verb"to" \verb"the" \verb"construction" \verb"of" \verb"optimally"-\verb"adaptive" \verb"estimators".

\smallskip

One of the possible approaches to solving \verb"Problems VI" and \verb"VII" could be an $\bL_p$-aggregation of locally-bayesian or $M$-estimators. Some recent  advances in this direction can be found in \cite{Chi}, \cite{ChiLed}.

\paragraph{Unknown distribution of the noise} Suppose now that the density $\varrho$ is unknown or even does not exist. The goal is to consider simultaneously
the noises with and without moments, regular or irregular etc. Even if the regression function belong to {\it known} functional class the different noises may lead to different minimax rates of convergence.

\smallskip

\noindent \verb"PROBLEM VIII:"\; \verb"Build" \verb"an" \verb"estimator" \verb"which" \verb"would" \verb"simultaneously"   \verb"adapt" \verb"to" \verb"a"
 \verb"a"  \verb"given" \verb"scale" \verb"of" \verb"functional" \verb"classes" \verb"and" \verb"to" \verb"the" \verb"noise" \verb"distribution".

\smallskip

We do not precise here the collection of classes since the formulated problem seems extremely complicated.  Any solution even in the dimension 1 can be considered as the great progress. In this context let us mention   very promising results recently obtained in \cite{B-B-S}.

We finish this section with following remarks. The regression model is very reach and many other problems can be formulated in the framework of it. For instance, the discussed problems can be mixed with imposing structural assumptions in the model. On the other hand aforementioned problems are not directly related to the concrete statistical model. In particular, almost all of them can be postulated in the {\it inverse problem} estimation context or in nonparametric auto-regression.

\section{Proof of Theorem \ref{th:L_p-norm-oracle-inequality} and Corollary \ref{cor:th:L_p-norm-oracle-inequality}}
\label{sec:proofs}

We start this section with presenting the constants appearing in the assertion of Theorem \ref{th:L_p-norm-oracle-inequality}.

\subsection{{Important quantities}}
\label{sec:constants}

Put
\begin{eqnarray*}
C_1&=&2\big(q\vee\big[p\mathrm{1}\{p<\infty\}+\mathrm{1}\{p=\infty\}\big]\big)+2\sqrt{2d}\Big[\sqrt{\pi}+ \|K\|_2\Big(\sqrt{\big|\ln{\big(4bA\|K\|_2\big)}\big|} +1\Big)\Big];
\\*[2mm]
C_3&=&C_3(\tilde{q},p)\mathrm{1}\{p<\infty\}+C_3(q,1)\mathrm{1}\{p=\infty\},
\quad\;\; \tilde{q}=(q/p)\vee 1;
\\*[2mm]
C_4&=&\Big(\boldsymbol{\gamma}_{q+1}\sqrt{(\pi/2)}\big[1\vee(2b)^{qd}\big]\sum_{r\in\bN_p^*}e^{-e^{r}}\big[\big(r\sqrt{e}\big)^d
\|\cK\|_{\frac{2r}{r+2}}^d\big]^{\frac{q}{2}}\Big)^{\frac{1}{q}},
\end{eqnarray*}
where $\boldsymbol{\gamma}_{q+1}$ is the $(q+1)$-th absolute moment of the standard normal distribution and
$$
C_3(u,v)=2^{\frac{d}{v}}\bigg[2u\int_0^\infty z^{u-1}\exp\bigg(-\frac{z^{\frac{2}{v}}}{8\|K\|^{2}_2}\bigg)\rd z\bigg]^{\frac{1}{uv}},\; u,v\geq 1.
$$

\subsection{{Auxiliary results}}
\label{sec:Proof of_th:L_p-norm-oracle-inequality}

As it was already mentioned the main ingredient of the proof of Theorem \ref{th:L_p-norm-oracle-inequality} is the fact that $\big\{\Psi_\e\big(\vec{h}\big), \vec{h}\in\bH\big\}$ is  the upper function   for the collection   $\big\{\|\xi_{\vec{h}}\|_p, \vec{h}\in\bH\big\}$. The corresponding result  formulated below for citation convenience as Proposition \ref{prop:upper-function-for-l_p-oracle}  is proved  in Theorem 1 and in   Corollary 1 of Theorem 2, \cite{lepski13b}.

Set for any $p\in [1,\infty)$, $\tau\in(0,1)$ and $\cL>0$
\begin{eqnarray*}
&&\psi_\e\big(\vec{h}\big)=
\widetilde{\Psi}_{\e,p}\big(\vec{h}\big)\wedge
\Big(\inf_{r\in\bN^*_p(\vec{h},\cA_\e)}C_2(r,\tau,\cL)\Big\| V^{-\frac{1}{2}}_{\vec{h}}\Big\|_{\frac{rp}{r-p}}\Big),\quad \vec{h}\in\bB(\cA_\e)
\end{eqnarray*}

\begin{proposition}
\label{prop:upper-function-for-l_p-oracle}
Let    $\cL>0$  be fixed and let $\mh_\e$ and $\cA_\e$ are defined in (\ref{eq1:choice-of-parameters}). Suppose also that
 $K$ satisfies Assumption \ref{ass:kernel}.

 Then for any
$q\geq 1$ and $\tau\in (0,1)$ one can find $\e(\tau,q)$ such that

1) for any $p\in [1,\infty)$,  $\e\leq\e(\tau,q)$
and   any countable  $\mathrm{H}\subset\mS_d(\mh_\e)\cap\bH_d\big(\tau,\cL,\cA_\e\big)$  one has
\begin{eqnarray*}
&&\bE\bigg\{\sup_{\vec{h}\in\mathrm{H}}
\Big[\big\|\xi_{\vec{h}}\big\|_p-\psi_\e\big(\vec{h}\big)\Big]_+\bigg\}^{q}
\leq \big\{(C_3+C_4)\e\big\}^q;
\end{eqnarray*}
\par
2) for any $p\in [1,\infty]$,  $\e\leq\e(\tau,q)$ and   any countable  $\mathrm{H}\subset\mS_d$
\begin{eqnarray*}
&&\bE\bigg\{\sup_{\vec{h}\in\mathrm{H}}
\Big[\big\|\xi_{\vec{h}}\big\|_p-\widetilde{\Psi}_{\e,p}\big(\vec{h}\big)\Big]_+\bigg\}^{q}
\leq \big\{C_3\e\big\}^q.
\end{eqnarray*}

\end{proposition}

We will need also the following technical result.
\begin{lemma}
\label{lem:about-maximum}
For any  $d\geq 1$, $\kappa\in (0,1/d)$, $\mL> 0$ and $\cA\geq e^{d}$
\begin{eqnarray*}
&(\mathbf{i})&\qquad \bH_d(\kappa,\mL,\cA)\subseteq\bH_d\Big(d\kappa,\mL^{d},\cA\Big);
\\
&(\mathbf{ii})&\qquad
\vec{h}\vee\vec{\eta}\in\bH_d\Big(d\kappa,(2\mL)^d,\cA\Big),\quad\forall\; \vec{h},\vec{\eta}\in\bH_d(\kappa,\mL,\cA).
\end{eqnarray*}

\end{lemma}

The  first statement of the lemma is obvious and the second one will be proved in Appendix.

\subsection{{Proof of Theorem \ref{th:L_p-norm-oracle-inequality}}}
Let $\vec{h}\in\bH$ be fixed.
We have  in view of the triangle inequality
\begin{equation}
\label{eq1:proof-th:L_p-norm-oracle-inequality}
\left\|\widehat{f}_{\vec{\mathbf{h}}}-f\right\|_p\leq\left\|\widehat{f}_{\vec{\mathbf{h}}\vee\vec{h}}-\widehat{f}_{\vec{\mathbf{h}}}\right\|_p+
\left\|\widehat{f}_{\vec{\mathbf{h}}\vee\vec{h}}-\widehat{f}_{\vec{h}}\right\|_p+\left\|\widehat{f}_{\vec{h}}-f\right\|_p.
\end{equation}
First, note that $\widehat{f}_{\vec{\mathbf{h}}\vee\vec{h}}\equiv \widehat{f}_{\vec{h}\vee\vec{\mathbf{h}}}$ and, therefore,
\begin{eqnarray}
\label{eq2:proof-th:L_p-norm-oracle-inequality}
\left\|\widehat{f}_{\vec{\mathbf{h}}\vee\vec{h}}-\widehat{f}_{\vec{\mathbf{h}}}\right\|_p&=&
\left\|\widehat{f}_{\vec{h}\vee\vec{\mathbf{h}}}-\widehat{f}_{\vec{\mathbf{h}}}\right\|_p\leq \widehat{\cR}_{\bH}\big(\vec{h}\big)+ \e\Psi_{\e,p}\big(\vec{h}\vee\vec{\mathbf{h}}\big)+\e\Psi_{\e,p}\big(\vec{\mathbf{h}}\big)
\nonumber\\
&\leq&\widehat{\cR}_{\bH}\big(\vec{h}\big)+ 2\e\Psi_{\e,p}\big(\vec{\mathbf{h}}\big).
\end{eqnarray}
Here we used first,  $\vec{\mathbf{h}}\in\bH$, and then   that  $V_{\vec{h}\vee\vec{\eta}}\geq V_{\vec{h}}\vee V_{\vec{\eta}}$    implies
$\Psi_{\e,p}\big(\vec{h}\vee\vec{\eta}\big)\leq \Psi_{\e,p}\big(\vec{h}\big)\wedge\Psi_{\e,p}\big(\vec{\eta}\big)$ for any $\vec{h}$ and $\vec{\eta}$.
Similarly we have
\begin{eqnarray}
\label{eq3:proof-th:L_p-norm-oracle-inequality}
&&\left\|\widehat{f}_{\vec{\mathbf{h}}\vee\vec{h}}-\widehat{f}_{\vec{h}}\right\|_p\leq \widehat{\cR}_{\bH}\big(\vec{\mathbf{h}}\big)+ \e\Psi_{\e,p}\big(\vec{h}\vee\vec{\mathbf{h}}\big)+\e\Psi_{\e,p}\big(\vec{h}\big)
\leq\widehat{\cR}_{\bH}\big(\vec{\mathbf{h}}\big)+2\e\Psi_{\e,p}\big(\vec{h}\big).
\end{eqnarray}
The definition of $\vec{\mathbf{h}}$ implies
$$
\widehat{\cR}_{\bH}\big(\vec{h}\big)+2\e\Psi_{\e,p}\big(\vec{h}\big)+\widehat{\cR}_{\bH}\big(\vec{\mathbf{h}}\big)+
2\e\Psi_{\e,p}\big(\vec{\mathbf{h}}\big)
\leq 4\widehat{\cR}_{\bH}\big(\vec{h}\big)+4\e\Psi_{\e,p}\big(\vec{h}\big)+2\e,
$$
and we get from (\ref{eq1:proof-th:L_p-norm-oracle-inequality}), (\ref{eq2:proof-th:L_p-norm-oracle-inequality}) and (\ref{eq3:proof-th:L_p-norm-oracle-inequality})
\begin{equation}
\label{eq4:proof-th:L_p-norm-oracle-inequality}
\big\|\widehat{f}_{\vec{\mathbf{h}}}-f\big\|_p\leq  4\widehat{\cR}_{\bH}\big(\vec{h}\big)+4\e\Psi_{\e,p}\big(\vec{h}\big)+\big\|\widehat{f}_{\vec{h}}-f\big\|_p+2\e.
\end{equation}
We obviously have for any $\vec{h},\vec{\eta}\in\bH$
$$
\left\|\widehat{f}_{\vec{h}\vee\vec{\eta}}-\widehat{f}_{\vec{\eta}}\right\|_p\leq \big\|B_{\vec{h},\vec{\eta}}\big\|_p+\e\big\|\xi_{\vec{h},\vec{\eta}}\big\|_p+\e\big\|\xi_{\vec{\eta}}\big\|_p.
$$
  Denote $\bH^*=\left\{\vec{v}:\;\;\vec{v}=\vec{h}\vee\vec{\eta},\;\vec{h},\vec{\eta}\in\bH\right\}$ and remark that
  $\bH\subseteq\bH^*\subseteq \bH\big(d\kappa,(2\mL)^d,\cA_\e\big)$. The latter inclusion follows from assertions of Lemma \ref{lem:about-maximum}.
  Moreover since $\bH$ is countable $\bH^*$ is countable as well.
Putting
$
\zeta=\sup_{\vec{v}\in\bH^*}\Big[\big\|\xi_{\vec{v}}\big\|_p- \Psi_{\e,p}\big(\vec{v}\big)\Big]_+,
$
we obtain
$$
\widehat{\cR}_{\bH}\big(\vec{h}\big)\leq \sup_{\vec{\eta}\in\bH}\big\|B_{\vec{h},\vec{\eta}}(\cdot,f)\big\|_p+2\e\zeta
$$
and, therefore, in view of  (\ref{eq4:proof-th:L_p-norm-oracle-inequality})
\begin{equation*}
\big\|\widehat{f}_{\vec{\mathbf{h}}}-f\big\|_p\leq 4\sup_{\vec{\eta}\in\bH}\big\|B_{\vec{h},\vec{\eta}}(\cdot,f)\big\|_p+4\e\Psi_{\e,p}(\vec{h})+8\e\zeta+\big\|\widehat{f}_{\vec{h}}-f\big\|_p
+2\e.
\end{equation*}
Taking into account that
$
\big\|\widehat{f}_{\vec{h}}-f\big\|_p\leq\big\|B_{\vec{h}}\big\|_p+\e\big\|\xi_{\vec{h}}\big\|_p,
$
we obtain
\begin{equation*}
\label{eq12:proof-th:L_p-norm-oracle-inequality}
\big\|\widehat{f}_{\vec{\mathbf{h}}}-f\big\|_p\leq 5\cB^{(p)}_{\vec{h}}(f)+5\e\Psi_{\e,p}\big(\vec{h}\big)+9\e\zeta+2\e.
\end{equation*}
It remains to note that if $p\in [1,\infty)$ in view of the definition of $\Psi_{\e,p}(\cdot)$
$$
\zeta=\bigg(\sup_{\vec{v}\in\bH^*\cap\mS_d(\mh_\e)}\Big[\big\|\xi_{\vec{v}}\big\|_p- \overline{\Psi}_{\e,p}\big(\vec{v}\big)\Big]_+\bigg)\vee
\bigg(\sup_{\vec{v}\in\bH^*\setminus\mS_d(\mh_\e)}\Big[\big\|\xi_{\vec{v}}\big\|_p- \widetilde{\Psi}_{\e,p}\big(\vec{v}\big)\Big]_+\bigg)
$$

Applying the first and the second assertions of Proposition \ref{prop:upper-function-for-l_p-oracle} with $\tau=d\kappa$, $\cL=(2\mL)^d$,   $\mathrm{H}=\bH^{*}\cap\mS_d(\mh_\e)$ and
$\mathrm{H}=\bH^{*}\setminus\mS_d(\mh_\e)$ respectively,  we obtain
$$
\cR^{(p)}_\e\big[\tilde{f}_{\vec{\mathbf{h}}}; f\big]\leq 5\cB^{(p)}_{\vec{h}}(f)+
5\e\Psi_{\e,p}\big(\vec{h}\big)+18(C_3+C_4+2)\e.
$$
It remains to note that the left hand side of the obtained inequality is independent of $\vec{h}$ and we come to the assertion of the theorem
with  $\Upsilon=18(C_3+C_4+2)$ where, recall,  $C_3$  and $C_4$ are given in Section \ref{sec:constants}.

If $p=\infty$ the second assertion of  Proposition \ref{prop:upper-function-for-l_p-oracle} with $\mathrm{H}=\bH^{*}$ is directly applied to the random variable $\zeta$ and the statement of the theorem follows.
\epr

\subsection{{Proof of Corollary \ref{cor:th:L_p-norm-oracle-inequality}}}
\label{sec:proof-cor:th:L_p-norm-oracle-inequality}

The proof of the corollary consists  mostly in  bounding from above the quantity
$\cB^{(p)}_{\vec{h}}(f)$. This, in its turn, is based on the technical result presented in Lemma \ref{lem:techniclal-profcor:th:L_p-norm-oracle-inequality} below which will be used in the proof of Proposition \ref{prop:key-prop}, Section \ref{sec:subsec-proof-prop:key-prop}, as well.

\subsubsection{\textsf{Auxiliary lemma}}
The following notations will be exploited in the sequel.

For any $J\subseteq\{1,\ldots d\}$ and  $y\in\bR^d$ set $y_J=\{y_j,\;j\in J\}\in\bR^{|J|}$ and we will write $y=\big(y_J,y_{\bar{J}}\big)$, where as usual $\bar{J}=\{1,\ldots d\}\setminus J$.

For any
$j=1,\ldots,d$ introduce $\mathbf{E}_j=(\mathbf{0},\ldots,\mathbf{e}_j,\ldots,\mathbf{0})$ and set $\mathbf{E}[J]=\sum_{j\in J}\mathbf{E}_j$.
Later on $\mathbf{E}_0=\mathbf{E}[\emptyset]$ denotes the matrix with zero entries.

To any  $J\subseteq\{1,\ldots d\}$  and  any $\lambda:\bR^d\to\bR_+$ such that $\lambda\in\bL_p(\bR^d)$,  associate the function
$$
\lambda\big(y_J,z_{\bar{J}}\big)=\lambda\big(z+\mathbf{E}[J](y-z)\big),\quad y,z\in\bR^d,
$$
with the obvious agreement $\lambda_{J}\equiv\lambda$ if $J=\{1,\ldots d\}$ that is always the case if $d=1$.

At last for any  $\vec{\mathrm{h}}=\big(\mathrm{h}_1,\ldots,\mathrm{h}_d\big)\in\mS^{\text{const}}_d$ and $J\subseteq\{1,\ldots d\}$ set
$
K_{\vec{\mathrm{h}},J}(u_J)=\prod_{j\in J}\mathrm{h}^{-1}_j\cK\big(u_j/\mathrm{h}_j\big)
$
and define for any $y\in\bR^d$
$$
\big[K_{\vec{\mathrm{h}}}\star\lambda\big]_{J}(y)=\int_{\bR^{|\bar{J}|}}K_{\vec{\mathrm{h}},\bar{J}}(u_{\bar{J}}-y_{\bar{J}})\lambda\big(y_J,u_{\bar{J}}\big)
\nu_{|\bar{J}|}\big(\rd u_{\bar{J}}\big)
$$
The following result is a trivial consequence of the Young inequality and Fubini theorem. For any $J\subseteq\{1,\ldots d\}$ and $p\in[1,\infty]$.
\begin{equation}
\label{eq5:proof-th:L_p-norm-oracle-inequality}
\Big\|\big[K_{\vec{\mathrm{h}}}\star\lambda\big]_{J}\Big\|_{p}\leq \|\cK\|^{d-|J|}_{1,\bR}\|\lambda\|_{p,\cA_{J}},\quad\forall\mathrm{h}\in\mS^{\text{const}}_d,
\end{equation}
where we have denoted $\cA_J=(-b,b)^{|J|}\times\bR^{|\bar{J}|}$.

\begin{lemma}
\label{lem:techniclal-profcor:th:L_p-norm-oracle-inequality}
 For any $\vec{h},\vec{\eta}\in\mS^{\text{const}}_d$ one can find $k=1,\ldots d,$ and the collection of indexes
 
 \noindent 
 $\left\{j_1<j_2<\cdots<j_k\right\}\in \{1,\ldots,d\}$ such that for any  $x\in\bR^d$ and any  $f:\bR^d\to\bR$
\begin{eqnarray*}
B_{\vec{h},\vec{\eta}}(x,f)&\leq& \sum_{l=1}^k\Big(\left[\big|K_{\vec{h}\vee\vec{\eta}}\big|\star b_{\vec{h},j_l} \right]_{\cJ_l}(x)+
\left[\big|K_{\vec{\eta}}\big|\star b_{\vec{h},j_l} \right]_{\cJ_l}(x)\Big);
\\
B_{\vec{h}}(x,f)&\leq& \sum_{l=1}^k
\left[\big|K_{\vec{h}}\big|\star b_{\vec{h},j_l} \right]_{\cJ_l}(x),\qquad \cJ_l=\{j_1,\ldots, j_l\}.
\end{eqnarray*}

\end{lemma}

The proof of the lemma is postponed to Appendix.

\subsubsection{\textsf{Proof of the corollary}}

We obtain in view of the first assertion of  Lemma \ref{lem:techniclal-profcor:th:L_p-norm-oracle-inequality} and  (\ref{eq5:proof-th:L_p-norm-oracle-inequality})
\begin{eqnarray*}
 \Big\|B_{\vec{h},\vec{\eta}}(\cdot,f)\Big\|_p&\leq& 2\sum_{l=1}^k\|\cK\|^{d-l}_{1,\bR}\big\|b_{\vec{h},j_l}\big\|_{p,\cA_{J_l}}=2\sum_{l=1}^k\|\cK\|^{d-l}_{1,\bR}\big\|b_{\vec{h},j_l}\big\|_{p}.
\end{eqnarray*}
The latter equality follows from the fact that $f$ is compactly supported on $(-b,b)^d$ that implies that $b_{\vec{h},j_l}(x_{\cJ_l},\cdot)$ is compactly supported on $(-b,b)^{d-l}$.
Taking into account that $\|\cK\|_{1,\bR}\geq 1$  we get
$$
\big\|B_{\vec{h},\vec{\eta}}(\cdot,f)\big\|_p\leq 2\|\cK\|^{d}_{1,\bR}\sum_{l=1}^k\big\|b_{\vec{h},j_l}\big\|_{p}\leq  2\|\cK\|^{d}_{1,\bR}\sum_{j=1}^d\big\|b_{\vec{h},j}\big\|_{p},\quad \forall \vec{h},\vec{\eta}\in\mS^{\text{const}}_d.
$$
Since the right hand side of the latter inequality is independent of $\vec{\eta}$
we obtain
\begin{gather*}
\label{eq10000:proof-th:L_p-norm-oracle-inequality}
\sup_{\vec{\eta}\in\mS^{\text{const}}_d} \big\|B_{\vec{h},\vec{\eta}}(\cdot,f)\big\|_p\leq 2\|\cK\|^{d}_{1,\bR}\sum_{j=1}^d\big\|b_{\vec{h},j}\big\|_p.
\end{gather*}
Repeating previous computations and using the second assertion of  Lemma \ref{lem:techniclal-profcor:th:L_p-norm-oracle-inequality}  we have
\begin{gather}
\label{eq1:proof-cor}
\big\|B_{\vec{h}}(\cdot,f)\big\|_p\leq \|\cK\|^{d}_{1,\bR}\sum_{j=1}^d\big\|b_{\vec{h},j}\big\|_p
\end{gather}
for any $\vec{h}\in\mS^{\text{const}}_d$ and any $p\geq 1$. We obtain finally
\begin{gather}
\label{eq10:proof-th:L_p-norm-oracle-inequality}
\cB^{(p)}_{\vec{h}}(f)\leq 3\|\cK\|^{d}_{1,\bR}\sum_{j=1}^d\big\|b_{\vec{h},j}\big\|_p.
\end{gather}

$3^0.\;$ We obviously have
$
r_{\cA}(\vec{h})=\lfloor p\rfloor+1
$
for any  $\vec{h}\in\bH_\e$ and, therefore, for any $p\in [1,\infty)$
$$
\Psi_{\e,p}\big(\vec{h}\big)\leq (2b)^{\frac{d}{p}} V^{-\frac{1}{2}}_{\vec{h}}\inf_{r\in\bN^*_p}C_2(r)=\Psi^{(\text{const})}_{\e,p}\big(\vec{h}\big).
$$
It is also obvious that
\begin{equation*}
\Psi_{\e,\infty}\big(\vec{h}\big)=\Psi^{(\text{const})}_{\e,\infty}\big(\vec{h}\big),\quad\forall \vec{h}\in\bH_\e.
\end{equation*}
As it was already mentioned  $\mathfrak{S}_d^{\text{const}}\subset\bH_d(\varkappa,\mL)$ for any $\varkappa\in (0,1)$ and $\mL=(2b)^\varkappa$. Thus, choosing for example $\kappa=(2d)^{-1}$ we constat that $\bH^{\text{const}}_\e\subset\bH_d\big((2d)^{-1},(2b)^{\frac{1}{2d}},\cA_\e\big)$ and, moreover, $\bH^{\text{const}}_\e$ is obviously finite set.

The assertion of the corollary follows now from (\ref{eq10:proof-th:L_p-norm-oracle-inequality}) and Theorem \ref{th:L_p-norm-oracle-inequality}.
\epr


\section{Proof of Theorem \ref{th:no-consitency}}

The proof is organized as follows. First, we formulate
two auxiliary statements,
Lemmas~\ref{lem:tsyb_book-result} and~\ref{lem:lepski}.
Second, we present a general construction
of a finite set of functions employed in the proof of lower
bounds. Then we specialize the constructed set of functions
in different regimes and derive the announced lower bounds.
\subsection{{Proof of Theorems \ref{th:no-consitency}. Auxiliary lemmas}}
The first statement given
in Lemma~\ref{lem:tsyb_book-result} is a
simple consequence of Theorem~2.4 from \cite{Tsybakov}. Let $\bF$ be a given
set of real functions defined on $(-b,b)^d$.
\begin{lemma}
\label{lem:tsyb_book-result}
Assume that for any sufficiently small
 $\e>0$ one can find a positive real number $\rho_\e$
and  a finite subset of functions
$\big\{f^{(0)}, f^{(j)},\;j\in\cJ_\e\big\}\subset \bF
$
such that
\begin{eqnarray}
\label{eq:ass1-klp-lemma}
  &&\big\|f^{(i)}- f^{(j)}\big\|_p \geq 2\rho_\e,\qquad\;
\forall  i, j\in \cJ_\e\cup\{0\}:\;i\neq j;
\\*[2mm]
\label{eq:ass2-klp-lemma}
&& \limsup_{\e\to 0}\frac{1}{|\cJ_\e|^{2}}
\sum_{j\in\cJ_\e}\bE_{f^{(0)}}\Bigg\{
\frac{\rd \bP_{f^{(j)}}}{\rd \bP_{f^{(0)}}}(X^{(\e)})\Bigg\}^{2}=:C <\infty.
\end{eqnarray}
Then for any $q\geq 1$
$$
\liminf_{\e\to 0}
\inf_{\tilde f}\;
\sup_{f \in \bF}
\rho^{-1}_\e\left(\bE_f \big\|\tilde{f} - f\big\|^{q}_p\right)^{1/q}
\geq \left(\sqrt{C} +\sqrt{C+1} \right)^{-2/q},
$$
where infimum on the left hand side is taken over all possible estimators.
\end{lemma}
We will apply Lemma \ref{lem:tsyb_book-result} with $\bF= \bN_{\vec{r},d}(\vec{\beta},\vec{L}, M)$.

Next, we will need the result  being a generalization  of
the Varshamov--Gilbert lemma. It can be found in  \cite{rigollet-tsybakov11}, Lemma A3. In the  version established in Lemma \ref{lem:lepski} below we only provide with particular choice of the constants appeared in the latter result.

Let $\varrho_n$ be the  Hamming distance on $\{0,1\}^n$, $n\in\bN^*$, i.e.
$$
\varrho_n(a,b)=\sum_{j=1}^n  {\bf 1}\left\{a_j\neq b_j\right\}=\sum_{j=1}^n |a_j-b_j|,\quad a,b\in\{0,1\}^n.
$$
\begin{lemma}
\label{lem:lepski}

For any $m\geq 4$ there exist a subset $\cP_{m,n}$ of $\{0,1\}^n$ such that
$$
\big|\cP_{m,n}\big|\geq  2^{-m}(n/m-1)^{\frac{m}{2}},\qquad\sum^{n}_{k=1}a_k=m,\qquad\varrho_m\big(a,a^\prime\big)\geq m/2,\;\;\;\;\;\forall a,a^\prime\in\cP_{m,n}.
$$

\end{lemma}

\subsection{{Proof of Theorem \ref{th:no-consitency}. General construction of a finite set of functions}}
 This part of the proof is mostly based on the constructions and computations made in \cite{GL13}, proof of Theorem 3.
 For any $t\in\bR$ set
$$
g(t)=
e^{-1/(1-t^2)}\; \mathbf{1}_{[-1,1]}(t).
$$
 For any
$l=1,\ldots,d$ let
$b/2>\sigma_l=\sigma_l(\e)\to 0$, $\e\to 0$,  be the sequences
to be specified later. Let
$M_l=\sigma_l^{-1}$,
and without loss of generality assume that $M_l$, $l=1,\ldots, d$ are integer numbers.

Define  also
$$
x_{j,l}=-b+2j\sigma_l,\;\;\; j=1,\ldots, M_l,\;\;l=1,\ldots,d,
$$
and let $\cM=\{1,\ldots, M_1\}\times\cdots\times\{1,\ldots, M_d\}$.
For any $\mm=(\mm_1,\ldots,\mm_d)\in \cM$ define
\begin{eqnarray*}
&&\pi(\mm)=\sum_{j=1}^{d-1}(\mm_j-1)\bigg(\prod_{l=j+1}^d M_l\bigg)+\mm_d,\qquad G_\mm(x)=
\prod_{l=1}^dg\left(\frac{x_l-x_{\mm_l,l}}{\sigma_l}\right),\quad x\in\bR^d.
\end{eqnarray*}
Let $W$ be a subset of $\{0,1\}^{|\cM|}$.
Define a family of functions $\{f_w, w\in W\}$ by
$$
f_w(x)=A\sum_{\mm\in\cM}w_{\pi(\mm)}G_\mm(x),\;\;\;\;x\in \bR^d,
$$
where $w_j$, $j=1,\ldots, |\cM|$ are
the coordinates of $w$,  and $A$ is a parameter to be specified.

 Suppose  that the set $W$ is chosen so that
\begin{eqnarray}
\label{eq18:proof-th:lower-bound-in-L_p}
&& \varrho_{|\cM|}\big(w,w^\prime\big)\geq B,\quad\forall w,w^\prime\in W,
\end{eqnarray}
where, we remind, $\varrho_{|\cM|}$ is the Hamming distance on $\{0,1\}^{|\cM|}$.
Here $B=B(\e)\geq 1$ is a parameter to be specified. Let also $S_W:=\sup_{w\in W}|\{j:\;w_j\neq 0\}|$.
Note finally that $f_w, w\in W$, are compactly supported on $(-b,b)^d$.

\smallskip

Repeating the computations made in  \cite{GL13}, proof of Theorem 3, we  assert first that
if
\begin{eqnarray}
\label{eq14:proof-th:lower-bound-in-L_p}
&&A\sigma_l^{-\beta_l}
\bigg(S_W\prod_{j=1}^d\sigma_j\bigg)^{1/r_l}\leq C_1^{-1} L_l,\quad\forall l=1,\ldots, d
\end{eqnarray}
then $f_{w}\in\bN_{\vec{r},d}(\vec{\beta},\vec{L})$ for any $w\in W$. Here $C_1$ as well as $C_2$ and $C_3$ defined below  are the numerical constants completely determined by the function $g$.

Next, the condition (\ref{eq:ass1-klp-lemma})
of Lemma~\ref{lem:tsyb_book-result} is fulfilled with
\begin{eqnarray}
\label{eq20:proof-th:lower-bound-in-L_p}
&& \rho_\e=C_2A\bigg(B\prod_{j=1}^d\sigma_j\bigg)^{1/p},
\end{eqnarray}
which remains true if $p=\infty$ as well.
At last, we have
$
\big\|f_w\big\|_2^2\leq C_3A^{2}S_W\prod_{j=1}^d\sigma_j.
$

Set $f^{(0)}\equiv 0$ and let us verify condition (\ref{eq:ass2-klp-lemma}) of Lemma~\ref{lem:tsyb_book-result}.
First observe that in view of Girsanov formulae
$$
\frac{\rd \bP_{f_w}}{\rd \bP_{f^{(0)}}}\big(X^{(\e)}\big)=
\exp\bigg\{\e^{-1}\int f_wb(\rd t)-(2\e^2)^{-1}\big\|f_w\big\|^2_2\bigg\}.
$$
It yields for any $w\in W$
\begin{eqnarray*}
\bE_{f^{(0)}}\bigg\{
\frac{\rd \bP_{f_w}}{\rd \bP_{f^{(0)}}}\big(X^{(\e)}\big)
\bigg\}^{2}=\exp\Big\{\e^{-2}\big\|f_w\big\|_2^2
\Big\}\leq \exp\bigg\{\e^{-2}C_3A^{2}S_W\prod_{j=1}^d\sigma_j
\bigg\}.
\end{eqnarray*}
The right hand side of the latter inequality does not depend on $w$; hence we have
$$
\frac{1}{|W|^{2}}\sum_{w\in W}\bE_{f^{(0)}}
\Bigg\{
\frac{\rd \bP_{f_w}}{\rd \bP_{f^{(0)}}}\big(X^{(\e)}\big)
\Bigg\}^{2}\leq
\exp\bigg\{C_3\e^{-2} A^2S_W \bigg(\prod_{j=1}^d\sigma_j\bigg)-\ln{\big(|W|\big)}
\bigg\}.
$$
Therefore, the condition
(\ref{eq:ass2-klp-lemma}) of Lemma \ref{lem:tsyb_book-result} is fulfilled with $C=1$ if
\begin{eqnarray}
\label{eq22:proof-th:lower-bound-in-L_p}
C_3\e^{-2}A^2S_W \prod_{j=1}^d \sigma_j\leq\ln{\big(|W|\big)}.
\end{eqnarray}

\par
In order to apply Lemma~\ref{lem:tsyb_book-result}
it remains to specify the parameters $A$,
 $\sigma_l$, $l=1,\ldots,d$, and  the set $W$  so that
the relationships
(\ref{eq18:proof-th:lower-bound-in-L_p}), (\ref{eq14:proof-th:lower-bound-in-L_p})
and (\ref{eq22:proof-th:lower-bound-in-L_p}) are simultaneously fulfilled.
According to Lemma \ref{lem:tsyb_book-result},
 under these conditions the lower bound is given by $\rho_\e$ in
(\ref{eq20:proof-th:lower-bound-in-L_p}).

\subsection{{Proof of Theorems \ref{th:no-consitency}. Choice of the parameters}}

We begin with the construction of the set $W$.
Let $m\geq 4$ be an integer number whose choice will be made later,
and, without loss of generality, assume that $|\cM|/m\geq 9$ is integer.
Let $\cP_{m,|\cM|}$ be a subset of $\{0,1\}^{|\cM|}$ defined in Lemma \ref{lem:lepski}, where we put $n=|\cM|$.

Set
$W=\cP_{m,|\cM|}\cup{\mathbf{0}}$, where $\mathbf{0}$ is the zero sequence of the size $|\cM|$.
With such a set  $W$
$$
S_W\leq  m, \qquad \ln(|W|)\geq (m/2)\big[\ln_2\big(|\cM|/m-1\big)-2\big]
$$
and, therefore, condition (\ref{eq22:proof-th:lower-bound-in-L_p}) holds true if
\begin{eqnarray}
\label{eq24:proof-th:lower-bound-in-L_p}
&&A^{2}\e^{-2}\prod_{j=1}^d \sigma_j \leq (2C_3)^{-1}\big[\ln_2\big(|\cM|/m-1\big)-2\big].
\end{eqnarray}
  We also note that condition (\ref{eq14:proof-th:lower-bound-in-L_p}) is fulfilled if we require
\begin{eqnarray}
\label{eq26-new:proof-th:lower-bound-in-L_p}
A\sigma_l^{-\beta_l} \bigg(m\prod_{j=1}^d \sigma_j\bigg)^{1/r_l}\leq C_1^{-1} L_l,\quad\forall
l=1,\ldots, d.
\end{eqnarray}
In addition,
(\ref{eq18:proof-th:lower-bound-in-L_p})
holds with $B=m/2$ and, therefore
\begin{eqnarray}
\label{eq20-new2:proof-th:lower-bound-in-L_p}
&& \rho_\e=2^{-1/p}C_2A\bigg(m\prod_{j=1}^d\sigma_j\bigg)^{1/p}.
\end{eqnarray}

\subsection{{Proof of Theorem \ref{th:no-consitency}. Derivation of lower bounds in different zones}}

Let  $\mathbf{c_i}, \mathbf{i}=1,\ldots, 6,$
be constants those choice will be made later.

\paragraph{Case: $\kappa(p)\leq 0,\; \tau(p^*)\leq 0$}

Set
$$
\varpi_\e=\left\{
\begin{array}{ll}
\big(L_\beta\e^2|\ln(\e)|\big)^{\frac{\omega}{\kappa(p^*)}},\quad &\kappa(p^*)<0;
\\*[2mm]
L_\beta e^{-\e^{-2}},\quad &\kappa(p^*)=0;
\end{array}
\right.
$$
 and note that $\varpi_\e\to\infty,\;\e\to 0$.
 In view of the latter remark we will assume that $\e$ is small enough provided $\varpi_\e>1$.
 We start our considerations  with the following remark. The case $\kappa(p^*)=0$ is possible only if $p^*=p$ since $\kappa(\cdot)$ is strictly decreasing. Moreover,
in view of the relation (\ref{eq:relation-kappa-tau})
$\kappa(p^*)=0$ is possible only if $p\leq 2$ since $\tau(p^*)\leq 0.$
Choose
\begin{equation*}
A=\mathbf{c_1}\varpi_\e,\qquad
m=\mathbf{c_2} L_\beta\varpi_\e^{- p^*\tau(p^*)},\qquad
\sigma_l=
\mathbf{c_3}L_l^{-\frac{1}{\beta_l}}\varpi_\e^{\frac{r_l-p^*}{\beta_lr_l}}.
\end{equation*}
With this choice,
we have
$
\frac{|\cM|}{m}=m^{-1}\prod_{j=1}^d\sigma^{-1}_j=\mathbf{c_2}^{-1}\mathbf{c_3}^{-d}
\varpi_\e^{p^*}\to\infty,\;\e\to 0.
$
  Hence , for any $\e$ small enough one has
$$
\big[\ln_2\big(|\cM|/m-1\big)-2\big]\geq Q_1\left\{
\begin{array}{ll}
|\ln(\e)|,\quad & \kappa(p^*)<0;
\\
\e^{-2},\quad & \kappa(p^*)=0,
\end{array}
\right.
$$
where $Q_1$ is  independent of  of $\e$ and $\vec{L}$.
This yields that (\ref{eq24:proof-th:lower-bound-in-L_p}) and  (\ref{eq26-new:proof-th:lower-bound-in-L_p}) will be fulfilled if
\begin{eqnarray}
\label{eq27:proof-th:lower-bound-in-L_p}
&&\mathbf{c_1}^2\mathbf{c_3}^d\leq (2C_3)^{-1}Q_1,\qquad
\mathbf{c_1}\mathbf{c_2}^{\frac{1}{r_l}}\mathbf{c_3}^{\frac{d}{r_l}-\beta_l}
\leq C_1^{-1}.
\end{eqnarray}

Some remarks are in order. First, since $r_l\leq p^*$ for any $l=1,\ldots, d$ and  $\kappa(p^*)<0$  we have
$$
\sigma_l\leq \mathbf{c_3}L_l^{-\frac{1}{\beta_l}}\leq \mathbf{c_3}\Big[\min_{l=1,\ldots, d}L_0^{-\frac{1}{\beta_l}}\Big].
$$
Here we also used  $\varpi_\e>1$. Thus,  choosing $\mathbf{c_3}$ small enough we can guarantee that $\sigma_l\leq b/2$ for any   $l=1,\ldots, d$ that was the unique restriction imposed
on the choice of the latter sequence.

Next   $\tau(p^*)\leq 0$,  $\varpi_\e>1$  and $p^*\tau(p^*)=2-p^*$, when $\kappa(p^*)=0$, imply that
$
m\geq \mathbf{c_2}L_0^{\frac{1}{\beta}}
$
and, therefore, choosing $\mathbf{c_2}$ large enough we guarantee that  $m\geq 4$.
At last, choosing
$\mathbf{c_1}$ small enough we can assert that (\ref{eq27:proof-th:lower-bound-in-L_p}) is satisfied.

Thus, it remains to compute $\rho_\e$. We get from (\ref{eq20-new2:proof-th:lower-bound-in-L_p})
\begin{eqnarray}
\label{eq28:proof-th:lower-bound-in-L_p}
&& \rho_\e=C_22^{-1/p}\mathbf{c_1}\big(\mathbf{c_2}
\mathbf{c_3}^d\big)^{1/p}\varpi_\e^{1-\frac{p^*}{p}}=:\mathrm{1}_{(p,\infty]}(p^*)
Q_2\big(L_\beta\e^2|\ln(\e)|\big)^{\frac{\omega(p^*-p)}{p(p^*-\omega(2+1/\beta)}}.
\end{eqnarray}
We remark that there is no uniformly consistent estimators if $p^*=p$.

\paragraph{Case: $\kappa(p)\leq 0,\; \tau(p^*)> 0$} First note, that
the case  $\kappa(p)\leq 0,\; \tau(p^*)> 0$, is possible only if $p>2$. It follows from  (\ref{eq:relation-kappa-tau}) and $\tau(p^*)\leq \tau(p)$ since
$\tau(\cdot)$ is decreasing. It implies
  $\tau(2)>0$ and
 $\tau(r_l)>0$ for any $l=1,\ldots d,$ since $r_l\leq p^*$.

\smallskip

Set $\varpi_\e=\e^2|\ln(\e)|$ and choose
\begin{equation*}
A=\mathbf{c_4}L_\beta^{\frac{1}{2\tau(2)}}\varpi_\e^{\frac{1-1/\omega}{2\tau(2)}},\qquad
m=4,\qquad
\sigma_l=
L_l^{-\frac{1}{\beta_l}}L_\beta^{\frac{r_l-2}{2\beta_lr_l\tau(2)}}\varpi_\e^{\frac{\tau(r_l)}{2\beta_l\tau(2)}}.
\end{equation*}
We remark, first that
$$
\sigma_l\to 0,\;\; \e\to 0,\;\; \forall l=1,\ldots d,
$$
and, therefore, $\sigma_l\leq b/2$ for all $\e>0$ small enough.

Next,
$
|\cM|/m=4^{-1}L_\beta^{\frac{1}{\tau(2)}}
\varpi_\e^{-\frac{1}{2\beta\tau(2)}}\geq 4^{-1}L_0^{\frac{1}{\beta\tau(2)}}
\varpi_\e^{-\frac{1}{2\beta\tau(2)}}
$
 and, hence , for any $\e$ small enough
$$
\big[\ln_2\big(|\cM|/m-1\big)-2\big]\geq Q_3|\ln(\e)|,
$$
where $Q_3$ is  independent of  of $\e$ and $\vec{L}$.
This yields that (\ref{eq24:proof-th:lower-bound-in-L_p}) and  (\ref{eq26-new:proof-th:lower-bound-in-L_p}) will be fulfilled if
\begin{eqnarray*}
&&\mathbf{c_4}^2\leq (2C_3)^{-1}Q_3,\qquad
\mathbf{c_4}4^{\frac{1}{r_l}}
\leq C_1^{-1}.
\end{eqnarray*}
Choosing $\mathbf{c_4}$ small enough we satisfy the latter restrictions. Thus, it remains to compute $\rho_\e$. We get from (\ref{eq20-new2:proof-th:lower-bound-in-L_p})
\begin{eqnarray}
\label{eq30:proof-th:lower-bound-in-L_p}
&& \rho_\e=C_22^{1/p}\mathbf{c_4}L_\beta^{\frac{1-2/p}{2\tau(2)}}\varpi_\e^{\frac{\tau(p)}{2\tau(2)}}=:
Q_4\bigg(L_\beta^{\frac{1-2/p}{\tau(p)}}\e^2|\ln(\e)|\bigg)^{\frac{1-1/\omega+1/(\beta p)}{2-2/\omega+1/\beta}}.
\end{eqnarray}

\paragraph{Case: $\kappa(p)>0$}

Choose
\begin{equation*}
A=\mathbf{c_6}\big(L_\beta\e^2\big)^{\frac{\beta}{2\beta+1}},\qquad
m=9^{-1}L_\beta\big(L_\beta\e^2\big)^{-\frac{\beta}{2\beta+1}},\qquad
\sigma_l=
L_l^{-\frac{1}{\beta_l}}\big(L_\beta\e^2\big)^{\frac{\beta}{\beta_l(2\beta+1)}}.
\end{equation*}
We remark that $|\cM|/m=9$ and  $m\to\infty,\;\e\to 0$;  hence $m>4$. Moreover
$$
\sigma_l\to 0,\;\; \e\to 0,\;\; \forall l=1,\ldots d,
$$
and, therefore, $\sigma_l\leq b/2$ for all $\e>0$ small enough.

We obviously get that (\ref{eq24:proof-th:lower-bound-in-L_p}) and  (\ref{eq26-new:proof-th:lower-bound-in-L_p}) will be fulfilled if
\begin{eqnarray*}
\label{eq29:proof-th:lower-bound-in-L_p}
&&\mathbf{c_6}^2\leq (2C_3)^{-1},\qquad
\mathbf{c_6}
\leq (9C_1)^{-1}.
\end{eqnarray*}
Choosing $\mathbf{c_6}$ small enough we satisfy the latter restrictions. Finally,  we get  from (\ref{eq20-new2:proof-th:lower-bound-in-L_p})
\begin{eqnarray*}
\label{eq300:proof-th:lower-bound-in-L_p}
&& \rho_\e=C_2\mathbf{c_4}18^{-1/p}\big(L_\beta\e^2\big)^{\frac{\beta}{2\beta+1}}.
\end{eqnarray*}
\epr

\section{Proof of Theorem \ref{th:adaptive-upper-bound}}
\label{sec:proof-th:adaptive-upper-bound}

Later on  $\mathbf{c}_i,\;i=1,2,\ldots,$ denote numerical constants   independent of
$\vec{L}$.  Moreover without further mentioning we will assume that all quantities those definitions involve the kernel $\cK$
 are defined with $\cK=w_\ell$.

\subsection{{Preliminary facts. Embedding of Nikolskii classes}}
\label{sec:subsec-inclusions}

For any  $\vec{\beta}\in (0,\infty)^d$,  $\vec{r}\in [1,\infty]^d$ and $s\geq 1$ define
\begin{eqnarray}
\label{eq:gamma-and-q-new}
\bga_j(s)&=&\frac{\beta_j\tau(s)}{\tau(r_j)},\;\;j=1,\ldots, d,\quad\; \vec{\bga}(s)=\big(\bga_1(s)\wedge\beta_1,\ldots,\bga_d(s)\wedge\beta_d\big);
\\
\label{eq:gamma-and-q-new-new}
r^*(s)&=&\big[\max_{j=1,\ldots d}r_j\big]\vee s,\qquad\qquad\;\vec{r}(s)=\big(r_1\vee s,\ldots,r_d\vee s\big).
\end{eqnarray}

\begin{lemma}
\label{lem:nik-lep}
For any $s\geq 1$ provided $\tau\big(r^*(s)\big)>0$
\begin{equation*}
\label{eq:embedd-nik-new}
 \bN_{\vec{r},d}\big(\vec{\beta},\vec{L}\big) \subseteq
\bN_{\vec{r}(s),d}\big(\vec{\bga}(s),\mathbf{c}\vec{L}\big),
\end{equation*}
where constant $\mathbf{c}>0$ is independent of $\vec{L}$,  $\vec{r}$ and $\vec{\beta}$.
\end{lemma}
The statement of the lemma  is a generalization of
 the embedding theorem for anisotropic Nikol'skii classes $\bar{\bN}_{\vec{r},d}(\vec{\beta},\vec{L})$. Indeed, if $r^*(s)=s$
 the assertion of the lemma can be found in
\cite{Nikolski}, Section~6.9.
The proof of this lemma as well as whose of Lemma \ref{lem:bias-norm-bound} below is postponed to Appendix.

Define $\cJ_\pm=\{j=1,\ldots,d: r_j\neq\infty\}$,  $p_\pm=[\sup_{j\in\cJ_\pm}r_j]\vee p$ and  introduce
\begin{equation}
\label{eq:gamma-and-q}
q_j=\left\{\begin{array}{ll}
p_\pm,\quad & j\in\cJ_\pm,
\\
\infty,\quad & j\neq \cJ_\pm,
\end{array}
\right.,\;\;\;
\qquad\;\;
\gamma_j=\left\{
\begin{array}{ll}
\bga_j(p_\pm),\quad & j\in\cJ_\pm,
\\
\beta_j,\quad & j\neq \cJ_\pm.
\end{array}
\right.
\end{equation}
Note that $p^*\geq p_\pm$ and, therefore, if $\tau(p^*)>0$ we have in view of Lemma \ref{lem:nik-lep} with $s=p_\pm$
\begin{equation}
\label{eq:embedd-nik}
 \bN_{\vec{r},d}\big(\vec{\beta},\vec{L}\big) \subseteq
\bN_{\vec{q},d}\big(\vec{\gamma},\mathbf{c}\vec{L}\big),
\end{equation}

\begin{lemma}
\label{lem:bias-norm-bound}
 Let $f\in \bN_{\vec{r},d}\big(\vec{\beta},\vec{M}\big)$
and let $\ell> \max_{j=1,\ldots,d}\beta_j$. Then for any $\vec{\mathrm{h}}\in\mS^{\text{const}}_d$
\begin{equation}
\label{eq:bias-norms-1}
 \big\|b_{\vec{\mathrm{h}},j}(\cdot,f)\big\|_{\mathbf{r},\bR^d} \leq (2b+1)^d\|w_\ell\|_{1,\bR^d}\big(1-e^{-\beta_j}\big)^{-1}M_j \mathrm{h}_j^{\beta_j},\;\;\;\forall \mathbf{r}\in[1,r_j],\; j=1,\ldots,d.
\end{equation}
Moreover,  if   $\tau(p^*)>0$  then for any $p\geq 1$
\begin{equation}\label{eq:bias-norms-2}
\big\|b_{\vec{\mathrm{h}},j}(\cdot,f)\big\|_{q_j,\bR^d} \leq (2b+1)^d\|w_\ell\|_{1,\bR^d}\big(1-e^{-\gamma_j}\big)^{-1} M_j \mathrm{h}_j^{\gamma_j} ,\;\;\;\forall j=1,\ldots, d,
\end{equation}
where $\vec{\gamma}$ and $\vec{q}$ are defined in (\ref{eq:gamma-and-q}).

\end{lemma}

\subsection{{Preliminary facts. Maximal operator}}
Let $\lambda:\bR^m\to\bR, m\geq 1,$  be a locally integrable function. We define  the strong maximal function $M[\lambda]$ of $\lambda$ by formula
\begin{equation}\label{eq2:maximal-function}
 M[\lambda](x):= \sup_{\bK_m} \frac{1}{\nu_m(\bK_m)} \int_{\bK_m} \lambda(t) \nu_m(\rd t),\;\;\;x\in \bR^m,
\end{equation}
where the supremum is taken over all
possible hyper-rectangles $\bK_m$  in $\bR^m$ with sides parallel to the coordinate
axes, containing point $x$.
It is worth noting that the {\em Hardy--Littlewood maximal function} is defined by
(\ref{eq2:maximal-function}) with the supremum taken over all cubes with sides parallel to the
coordinate axes,
centered at $x$.
\par
It is well known that the strong maximal operator $\lambda\mapsto M[\lambda]$
is of the strong $(r,r)$--type for all $1<r\leq \infty$, i.e.,
if $\lambda\in \bL_r(\bR^m)$ then $M[\lambda]\in \bL_r(\bR^m)$ and for any $r>1$ there exists a constant $\bar{C}(r)$
depending on $r$ only such that
\begin{equation}
\label{eq:strong-max}
 \big\|M[\lambda]\big\|_{r,\bR^d} \leq \bar{C}(r) \|\lambda\|_{r,\bR^d}.
\end{equation}
Using the notations from Section \ref{sec:proof-cor:th:L_p-norm-oracle-inequality},  to   any $J\subseteq\{1,\ldots d\}\cup\emptyset$  and  locally integrable function $\lambda:\bR^d\to\bR_+$ we associate the operator
$$
M_{J}[\lambda](x)=\sup_{\bK_{|\bar{J}|}}\frac{1}{\nu_{|\bar{J}|}(\bK_{|\bar{J}|})}\int_{\bK_{|\bar{J}|}}\lambda\big(t+\mathbf{E}[J][x-t]\big)
\nu_{|\bar{J}|}(\rd t_{\bar{J}})
$$
where
 the supremum is taken over all hyper-rectangles in $\bR^{|\bar{J}|}$ with center  $x_{\bar{J}}=(x_j, j\in \bar{J})$ and with sides parallel to the axis.

As we see $M[\lambda]$ is the strong maximal operator applied to the function obtained from $\lambda$ by fixing of coordinates
those indices   belong to $J$. It is obvious  that $M_{\emptyset}[\lambda]\equiv M[\lambda]$ and $M_{\{1,\ldots,d\}}[\lambda]\equiv\lambda$.

The following result is the direct consequence of (\ref{eq:strong-max}) and Fubuni theorem.
For any  $r>1$  there exists $\mathbf{C}_{r}$  such that for any  $d\geq 1$,  $\lambda$, $J\subseteq\{1,\ldots d\}\cup\emptyset$ and $\mathrm{y}\in (0,\infty]$
\begin{equation}
\label{eq:strong-max-partial}
\big\|M_J[\lambda]\big\|_{r,(-\mathrm{y},\mathrm{y})^{d}}\leq \mathbf{C}_{r}\|\lambda\|_{r,\cT_J(\mathrm{y})},
\end{equation}
where we have denoted $\cT_J(\mathrm{y})=(-\mathrm{y},\mathrm{y})^{|J|}\times\bR^{|\bar{J}|}$.
Note also that $\mathbf{C}_\infty=1$.

\subsection{{Preliminary facts. Key proposition}}
The result presented in Proposition \ref{prop:key-prop} below is the milestone for the proof of Theorem \ref{th:adaptive-upper-bound}.
For any $(\vartheta,p)\in\cP^{\text{consist}}$ define
$$
\varphi:=\varphi_\e(\vartheta,p)=\left\{
\begin{array}{ll}
(L_\beta \e^2)^{\beta/(2\beta+1)},\quad&\kappa(p)>0;
\\*[2mm]
(L_\beta \e^2|\ln(\e)|)^{\beta/(2\beta+1)},\quad& \kappa(p)\leq 0.
\end{array}
\right.
$$
\paragraph{Special set of bandwidths}
For any  $(\vartheta,p)\in\cP^{\text{consist}}$, $m\in\bN$ and any  $j=1,\ldots,d$   set
\begin{eqnarray}
\label{eq1:def-tilde(mu)_j}
\widetilde{\bleta}_j(m)&=&e^{-2}\big(L_j^{-1}\varphi\big)^{1/\beta_j}
e^{2dm\big(\frac{1}{\beta_j}-\frac{\omega(2+1/\beta)}{\beta_jr_j}\big)}.
\\
\label{eq2::def-bar(mu)_j}
\widehat{\bleta}_j(m)&=&e^{-2}
\big(L_j^{-1}\varphi\big)^{1/\gamma_j}
e^{2dm\left(\frac{1}{\gamma_j}-\frac{\upsilon(2+1/\gamma)}{\gamma_jq_j}\right)}
\left[\frac{L_\gamma\varphi^{1/\beta}}{L_\beta\varphi^{1/\gamma}}\right]^{\frac{\upsilon}{\gamma_jq_j}},
\end{eqnarray}
where $\gamma_j$, $q_j$ are defined in (\ref{eq:gamma-and-q}) and
$\gamma$, $\upsilon$ and $L_\gamma$ are given by
\begin{equation}\label{eq:gamma-upsilon}
\frac{1}{\gamma}:=\sum_{j=1}^d \frac{1}{\gamma_j},\qquad
\frac{1}{\upsilon}:=\sum_{j=1}^d \frac{1}{\gamma_j q_j},\qquad L_\gamma:=\prod_{j=1}^dL_j^{1/\gamma_j}.
\end{equation}
Introduce  the integer $\widehat{\mathbf{m}}=\widehat{\mathbf{m}}(\vartheta,p),\; (\vartheta,p)\in\cP^{\text{consist}},$ satisfying
\begin{eqnarray*}
e^{-2d}\Big[(L_\gamma/L_\beta)^{\frac{1}{(1/\gamma-1/\beta)}}\varphi^{-1}\Big]^{\frac{1}{2\beta\omega\tau(2)}} &\leq&
e^{2d\widehat{\mathbf{m}}}\;\leq\;\Big[(L_\gamma/L_\beta)^{\frac{1}{(1/\gamma-1/\beta)}}\varphi^{-1}\Big]^{\frac{1}{2\beta\omega\tau(2)}}.
\end{eqnarray*}
 Later on $\widehat{\mathbf{m}}$ will be used only if  $\kappa(p)< 0$ and $\tau(p^*)>0$. Note that in this case $\widehat{\mathbf{m}}\geq 1$ for all $\e>0$ small enough since $\tau(2)>0$ (see, e.g., the proof of Theorem \ref{th:no-consitency}).

\smallskip

Introduce also the integer $\widetilde{\mathbf{m}}=\widetilde{\mathbf{m}}(\vartheta,p),\; (\vartheta,p)\in\cP^{\text{consist}}$ as follows.

\smallskip

\noindent{\bf Case $\kappa(p)> 0,\; \kappa(p^*)\geq 0$}:\;  $\widetilde{\mathbf{m}}=+\infty$.

\smallskip

\noindent{\bf Case $\kappa(p)> 0,\; \kappa(p^*)<0$}:\;\;\;
$
e^{-2d} \big(\mh^{-\ell}_\e L_0^{-1}\varphi\big)^{\frac{p^*}{\kappa(p^*)}}\leq
e^{2d\widetilde{\mathbf{m}}}\;\leq\;
\big(\mh^{-\ell}_\e L_0^{-1}\varphi\big)^{\frac{p^*}{\kappa(p^*)}}.
$

\smallskip

\noindent{\bf Case $\kappa(p)\leq 0,\;\tau(p^*)\leq 0$}:\;\;\;\;
$
e^{-2d} \big(L_0^{-1}\varphi\big)^{\frac{p^*}{\kappa(p^*)}}\leq
e^{2d\widetilde{\mathbf{m}}}\;\leq\;
\big(L_0^{-1}\varphi\big)^{\frac{p^*}{\kappa(p^*)}}.
$

\medskip

\noindent{\bf Case $\kappa(p)\leq 0,\;\tau(p^*)>0$}: $\widetilde{\mathbf{m}}=\widehat{\mathbf{m}}+1$ if $p^*=p$; $\widetilde{\mathbf{m}}=\widehat{\mathbf{m}}+\overline{\mathbf{m}}$ if $p^*>p$, where
\begin{eqnarray*}
e^{-2d}\varphi^{-\frac{1+(1/\gamma-1/\beta)\upsilon(1/p-1/p^*)}{(2+1/\gamma)\upsilon(1/p-1/p^*)}} &\leq&
e^{2d\overline{\mathbf{m}}}\;\leq\;\varphi^{-\frac{1+(1/\gamma-1/\beta)\upsilon(1/p-1/p^*)}{(2+1/\gamma)\upsilon(1/p-1/p^*)}},\;p^*>p.
\end{eqnarray*}



Some remarks are in order. First we note that
$\widetilde{\mathbf{m}}\geq 1$ for all $\e>0$ small enough. Indeed, $\varphi\to 0, \e\to 0,$ and
 $\kappa(p)\leq0$ implies $\kappa(p^*)<0$ if $p^*>p$. Moreover, since $(\vartheta,p)\in\cP^{\text{consist}}$ the case
 $\kappa(p)\leq 0,\;\tau(p^*)\leq 0$ is possible only if $p^*>p$ that, in its turn, implies $\kappa(p^*)<0$.

\smallskip

For any $(\vartheta,p)\in\cP^{\text{consist}}$ and any   $0\leq m\leq \widetilde{\mathbf{m}}$ introduce
\begin{eqnarray}
\label{eq2::def-hat(mu)_j}
&&\bar{\bleta}_j(m)=\left\{
\begin{array}{ll}
\widetilde{\bleta}_j(m)\mathrm{1}_{\{m\leq\widehat{\mathbf{m}}\}}+
\widehat{\bleta}_j(m)_{\mathrm{1}\{m>\widehat{\mathbf{m}}\}},\quad &\kappa(p)\leq 0,\;\tau(p^*)>0;
\\*[2mm]
\widetilde{\bleta}_j(m),&\text{otherwise},
\end{array}
\right.
\end{eqnarray}
and define $\vec{\bleta}(m)=\big(\bleta_1(m),\ldots,\bleta_d(m)\big)$
as follows.

For any $m\in\bN$ set $\bleta_j(m)=\mh_{s_j(m)}\in\mH$, where $\mathbf{s}(m)=\left(s_1(m),\ldots,s_d(m)\right)\in\bN^d$ is given by
\begin{eqnarray}
\label{eq2::def-projection}
s_j(m)=\min\{s\in\bN:\;\; \mh_{s}\leq \bar{\bleta}_j(m) \}.
\end{eqnarray}
Introduce finally the  set of bandwidths
$
\mH_\e(\vartheta,p)=\big\{\vec{\bleta}(m),\quad m=0,\ldots, \widetilde{\mathbf{m}}\big\}.
$

\begin{lemma}
\label{lem:set-of-band}
For any $(\vartheta,p)\in\cP^{\text{consist}}$ and any $\e>0$ small enough one has
$$
\mH_\e(\vartheta,p)\subset
\left\{
\begin{array}{ll}
\mH^d(\mh_\e),\quad& \kappa(p)>0;
\\*[2mm] \mH^d,\quad&\text{otherwise}.
\end{array}
\right.
$$
Moreover,
$
\mathbf{s}(m)\neq \mathbf{s}(n),\quad \forall m\neq n,\; m,n=0,\ldots, \widetilde{\mathbf{m}}.
$
\end{lemma}

\paragraph{Result formulation} For any $\vec{h}\in\mS_d$ and any $x\in\bR^d$ put
 $$
\Phi_\e\big(V_{\vec{h}}(x)\big)
 =\left\{
\begin{array}{ll}
 V^{-\frac{1}{2}}_{\vec{h}}(x),\quad & \kappa(p)>0;
\\*[2mm]
\Big[V^{-1}_{\vec{h}}(x)\big|\ln\big(\e V_{\vec{h}}(x)\big)\big|\Big]^{\frac{1}{2}},\quad & \kappa(p)\leq 0.
\end{array}
\right.
$$
For any $g:\bR^d\to\bR$  and any $\vec{h}\in\mS^{\text{const}}_d$ introduce
\begin{eqnarray*}
B^*_{\vec{h}}(x,g)&=&\sup_{\eta\in\mS^{\text{const}}_d}B_{\vec{h},\vec{\eta}}(x,g)+B_{\vec{h}}(x,g);
\\
b^*_{\vec{h}}(x,g)&=&\sup_{J\in\mJ}\sup_{j=1,\ldots d}M_{J}\big[b_{\vec{h},j}\big](x).
\end{eqnarray*}
 Here $B_{\vec{h},\vec{\eta}}$ and $ B_{\vec{h}}$ are defined in (\ref{eq:def-B_{vec{h}}}) and   $b_{\vec{h},j}$ is defined in (\ref{eq:def-small-b_h}), where $f$  is replaced  by $g$.

Let $\bC_\bK\big(\bR^d\big)$ denote the set of   continuous functions on $\bR^d$ compactly supported on $\bK=(-b-1,b+1)^d$ and let $\bN^*_{\vec{r},d}\big(\vec{\beta},\vec{L}\big)=\bN_{\vec{r},d}\big(\vec{\beta},\vec{L}\big)\cap\bC_\bK(\bR^d)$.
Remark that $\bN_{\vec{r},d}\big(\vec{\beta},\vec{L}\big)\subset \bC\big(\bR^d\big)$ if $\omega>1$  in view of (\ref{eq:embedd-nik}).

\smallskip

For any $\vec{\beta}\in (0,\infty)^d$ and $\vec{r}\in (1,\infty]^d$ set $\beta_*=\min_{j=1,\ldots,d}\beta_j$, $\mathbf{C}\big(\vec{r}\big)=\max_{j=1,\ldots,d}\mathbf{C}_{r_j}$ and define
$$
\Upsilon_1=3d\big(1\vee\|w_\ell\|_{\infty,\bR^d}\big)^{d},\qquad \Upsilon_2=4\Upsilon_1 \mathbf{C}\big(\vec{r}\big) (2b+1)^d\|w_\ell\|_{1,\bR^d}\big(1-e^{-\beta_*}\big)^{-1}.
$$

\begin{proposition}
\label{prop:key-prop}
For any  $\mathbf{a}\geq 1$,  $\ell\in\bN^*$,  $L_0>0$,
 any  $(\vartheta,p)\in\cP^{\text{consist}}$ and  $\e>0$ one can find $\mS_\e^*(\vartheta,p)=\big\{\vec{\mathbf{h}}:(-b,b)^d\to \mH_\e(\vartheta,p)\big\}$ such that
for any $\e>0$ small enough

\smallskip

1) $\mS_\e^*(\vartheta,p)\subset\bH_\e\big(3+\sqrt{2b}\big)$;

\smallskip

2) for
any  $g\in\bN^*_{\vec{r},d}\big(\vec{\beta},\mathbf{a}\vec{L}\big)$ there exists $\vec{\blh}_g\in\mS_\e^*(\vartheta,p)$ such that
$$
(\mathbf{i})\; B^*_{\vec{\blh}_g}(x,g)+\mathbf{a}\Upsilon_2\e \Phi_\e\big(V_{\vec{\blh}_g}(x)\big)\leq \inf_{\vec{h}\in\mH_{\e}(\vartheta,p)}\Big[\Upsilon_1 b^*_{\vec{h}}(x,g)+ \mathbf{a}\Upsilon_2\e \Phi_\e\big( V_{\vec{h}}\big)\Big]+\e,\;\forall x\in(-b,b)^d;
 $$

$\;\;(\mathbf{ii})$ if, $\kappa(p)>0$ there exists $\mathfrak{r}\in\bN^*_p$ such that $\kappa\big(\frac{\mathfrak{r}p}{\mathfrak{r}-p}\big)>0$ and
$
\mathfrak{r}\in\bN^*_p\big(\vec{\blh}_g,\cA_\e\big).
$

\end{proposition}

\subsubsection{\textsf{Proof of Proposition \ref{prop:key-prop}}} We break up the proof on several steps.
\label{sec:subsec-proof-prop:key-prop}

\smallskip

$1^0.\;$ The condition  $g\in\bC_\bK\big(\bR^d\big)$  implies that $g$ is uniformly continuous on $\bR^d$ and, therefore, for any $\e>0$ there exists
$\delta(\e)$ such that
\begin{equation}
\label{eq01:proof-key-lemma}
|g(y)-g(y^\prime)|\leq \e^2, \quad\forall y,y^\prime\in\bR^d:\: |y-y^\prime|\leq \delta(\e).
\end{equation}
Let $\mathrm{x}_{\mathbf{k},n}, \mathbf{k}\in\mK_{n}, n\in\bN^{*}$ denote the center of the cube
$\Delta^{(d)}_{\mathbf{k},n}$ defined in (\ref{eq:def-partition-collection}). Introduce
for any $\vec{h}\in\mS^{\text{const}}_d$
$$
\widetilde{B}^*_{\vec{h}}(x,g)=\sum_{\mathbf{k}\in\mK_{\widetilde{n}}}B^*_{\vec{h}}\big(\mathrm{x}_{\mathbf{k},\widetilde{n}},g\big)
\mathrm{1}_{\Delta^{(d)}_{\mathbf{k},\widetilde{n}}}(x),
$$
where $\widetilde{n}$ is chosen from the relation $2^{-\widetilde{n}}<\delta(\e)\leq 2^{-\widetilde{n}+1}$.

Our first goal is to prove that
\begin{equation}
\label{eq02:proof-key-lemma}
\sup_{\vec{h}\in\mS_d^{\text{const}}}\big\|B^*_{\vec{h}}(\cdot,g)-\widetilde{B}^*_{\vec{h}}(\cdot,g)\big\|_{\infty,\bR^d}\leq \mathbf{c}_1\e^2.
\end{equation}
Indeed, for any  $\vec{h}\in\mS^{\text{const}}_d$ since $g$ is compactly supported on $\bK$ one has
\begin{equation}
\label{eq03:proof-key-lemma}
\big\|B^*_{\vec{h}}(\cdot,g)-\widetilde{B}^*_{\vec{h}}(\cdot,g)\big\|_{\infty,\bR^d}=\sup_{\mathbf{k}\in\mK_{\widetilde{n}}}
\sup_{x\in\Delta^{(d)}_{\mathbf{k},\widetilde{n}}}\big|B^*_{\vec{h}}(x,g)-B^*_{\vec{h}}\big(\mathrm{x}_{\mathbf{k},\widetilde{n}},g\big)\big|.
\end{equation}
In view of the definition of $B^*_{\vec{h}}(\cdot,g)$ we have for any $x\in\Delta^{(d)}_{\mathbf{k},\widetilde{n}}$
\begin{eqnarray*}
\big|B^*_{\vec{h}}(x,g)-B^*_{\vec{h}}\big(\mathrm{x}_{\mathbf{k},\widetilde{n}},g\big)\big|&\leq& 3\sup_{\vec{\mathrm{h}}\in\mS_d^{\text{const}}}
\big|S_{\vec{\mathrm{h}}}(x,g)-S_{\vec{\mathrm{h}}}\big(\mathrm{x}_{\mathbf{k},\widetilde{n}},g\big)\big|+
\big|g(x)-g\big(\mathrm{x}_{\mathbf{k},\widetilde{n}}\big)\big|
\\
&\leq&3\sup_{\vec{\mathrm{h}}\in\mS_d^{\text{const}}}\big|S_{\vec{\mathrm{h}}}(x,g)-
S_{\vec{\mathrm{h}}}\big(\mathrm{x}_{\mathbf{k},\widetilde{n}},g\big)\big|+\e^2.
\end{eqnarray*}
The last inequality follows from   (\ref{eq01:proof-key-lemma}) and the definition of $\widetilde{n}$.

Recall that $K$ is given in (\ref{eq:w-function}) and
$$
S_{\vec{\mathrm{h}}}(x,g)=\int_{\bR^d}K_{\vec{\mathrm{h}}}(t-x)g(t)\nu_d(\rd t)=
\int_{\bR^d}K(u)g\big(x+u\vec{\mathrm{h}}\big)\nu_d(\rd u).
$$
Hence, for any $\vec{\mathrm{h}}\in\mS_d^{\text{const}}$ and any $x\in\Delta^{(d)}_{\mathbf{k},\widetilde{n}}$
\begin{eqnarray*}
\big|S_{\vec{\mathrm{h}}}(x,g)-S_{\vec{\mathrm{h}}}\big(\mathrm{x}_{\mathbf{k},\widetilde{n}},g\big)\big|
\leq \int_{\bR^d}\big|K(u)\big|\big|g\big(x+u\vec{\mathrm{h}}\big)-g\big(\mathrm{x}_{\mathbf{k},\widetilde{n}}+u\vec{\mathrm{h}}\big)\big|\nu_d(\rd u)
 \leq \|K\|_{1,\bR^d}\e^2
\end{eqnarray*}
in view of (\ref{eq01:proof-key-lemma}) and the definition of $\widetilde{n}$.

 Since the latter bound is independent of $\vec{\mathrm{h}}$ we obtain for any $\vec{h}\in\mS_d^{\text{const}}$ and any $x\in\Delta^{(d)}_{\mathbf{k},\widetilde{n}}$
$$
\big|B^*_{\vec{h}}(x,g)-B^*_{\vec{h}}\big(\mathrm{x}_{\mathbf{k},\widetilde{n}},g\big)\big|\leq \big(1+ 3\|K\|_{1,\bR^d}\big)\e^2.
$$
Taking into account that  the right hand side of the latter inequality is independent of $\vec{h}$, $\mathbf{k}$ and $x$ we deduce (\ref{eq02:proof-key-lemma})  from
(\ref{eq03:proof-key-lemma}).

One of the immediate  consequences of (\ref{eq02:proof-key-lemma})   is that for any $x\in\bR^d$
\begin{equation}
\label{eq04:proof-key-lemma}
\bigg|\inf_{\vec{h}\in\mH_\e(\vartheta,p)}\Big[B^*_{\vec{h}}(x,g)+ \mathbf{a}\Upsilon_2\e\Phi_\e\big( V_{\vec{h}}(x)\big)\Big]-
\inf_{\vec{h}\in\mH_\e(\vartheta,p)}\Big[\widetilde{B}^*_{\vec{h}}(x,g)+ \mathbf{a}\Upsilon_2\e\Phi_\e\big( V_{\vec{h}}(x)\big)\Big]\bigg|\leq\mathbf{c}_1\e^2.
\end{equation}

$2^0.\;$ For any $x\in(-b,b)^d$ introduce
\begin{equation}
\label{eq999:proof-key-lemma}
\vec{\blh}_g(x)=\arg\inf_{\vec{\mathrm{h}}\in\mH_\e(\vartheta,p)}\Big[ \widetilde{B}^*_{\vec{\mathrm{h}}}(x,g)+ \mathbf{a}\Upsilon_2\e\Phi_\e\big( V_{\vec{h}}(x)\big)\Big],
\end{equation}
and define $\mS^*_\e(\vartheta,p)=\big\{\vec{\blh}_g, \; g\in\bN^*_{\vec{r},d}\big(\vec{\beta},\mathbf{a}\vec{L}\big)\big\}$.

First, we deduce from (\ref{eq02:proof-key-lemma})  and (\ref{eq04:proof-key-lemma}) that for any $x\in(-b,b)^d$
\begin{equation}
\label{eq004:proof-key-lemma}
B^*_{\vec{\blh}_g}(x,g)+ \mathbf{a}\Upsilon_2\e\Phi_\e\big( V_{\vec{\blh}_g}(x)\big)\leq
\inf_{\vec{h}\in\mH_\e(\vartheta,p)}\Big[B^*_{\vec{h}}(x,g)+ \mathbf{a}\Upsilon_2\e\Phi_\e\big( V_{\vec{h}}(x)\big)\Big]+2\mathbf{c}_1\e^2.
\end{equation}
Next, since $\widetilde{B}^*_{\vec{\mathrm{h}}}(\cdot,g)$ is a piecewise constant on $\big\{\Delta^{(d)}_{\mathbf{k},\widetilde{n}}\cap(-b,b)^d,  \mathbf{k}\in\mK_{\widetilde{n}}\big\}$  one has  $\vec{\blh}_g\in\mS^{\e}_{\widetilde{n}}$ and, hence, we can assert that
\begin{equation}
\label{eq05:proof-key-lemma}
\vec{\blh}_g\in\bigcup_{n\in\bN^{*}}\mS^{\e}_{n},\quad\forall g\in\bC_\bK\big(\bR^d\big).
\end{equation}
Our  goal now is to prove that for any $(\vartheta,p)\in\cP^{\text{consist}}$ one can find $\e(\vartheta,p)>0$ such that for any $\e<\e(\vartheta,p)$
\begin{equation}
\label{eq1:proof-key-lemma}
\vec{\blh}_g\in \bH_d\big(1/(2d),3+\sqrt{2b},\;\cA_\e\big),\quad\forall g\in\bN^*_{\vec{r},d}\big(\vec{\beta},\vec{L}\big).
\end{equation}

$3^0a.\;$ Note that the definition of the function $w_\ell$ together with the assumption $g\in\bC_\bK\big(\bR^d\big)$ implies that
$\sup_{x\in\bR^d}\big|B^*_{\vec{h}}(x,g)\big|<\infty$ that implies in view of (\ref{eq02:proof-key-lemma}) $\sup_{x\in\bR^d}\big|\widetilde{B}^*_{\vec{h}}(x,g)\big|<\infty$.

Hence, $\blh_{j,g}(x)<\infty$ for any $x\in(-b,b)^d$ and any $j=1,\ldots d,$, where  $\blh_{j,g}(\cdot)$ is $j$-th coordinate of the vector-function $\vec{\blh}_{g}$. It implies, in particular, that the  infimum in (\ref{eq999:proof-key-lemma}) is achievable and, therefore, for any $x\in(-b,b)^d$
\begin{equation}
\label{eq2:proof-key-lemma}
\vec{\blh}_g(x)\in\mH_\e(\vartheta,p),\quad\forall g\in\bN^*_{\vec{r},d}\big(\vec{\beta},\mathbf{a}\vec{L}\big).
\end{equation}
 By the same reason $B^*_{\vec{h}}(\cdot,g)$ as well as  $\widetilde{B}^*_{\vec{h}}(\cdot,g)$ are Borel functions and since $\mH_\e(\vartheta,p)$ is countable we assert in view of (\ref{eq2:proof-key-lemma}) that for any  $\mathbf{s}\in\bN^d$ such that $\vec{\mh}_{\mathbf{s}}:=\big(\mh_{s_1},\ldots,\mh_{s_d}\big)\in\mH_\e(\vartheta,p)$
\begin{equation}
\label{eq3:proof-key-lemma}
\Lambda_{\mathbf{s}}\big[\vec{\blh}_g\big]\in\mB\big(\bR^d\big),\quad\forall g\in\bC_\bK\big(\bR^d\big).
\end{equation}
It implies, in particular, that $\vec{\blh}_g$ is Borel function.

\smallskip

$3^0b.\;$
Taking into account that $w_\ell$ is compactly supported on $[-1/2,1/2]^d$ we easily deduce from the assertions of Lemma 2 that
 for any $\vec{h},\vec{\eta}\in\mS^{\text{const}}_d$
\begin{eqnarray*}
B_{\vec{h},\vec{\eta}}(x,g)&\leq& 2d\big(1\vee\|w_\ell\|_{\infty,\bR^d}\big)^{d}\sup_{J\in\mJ}\sup_{j=1,\ldots d}M_{J}\big[b_{\vec{h},j}\big](x);
\\
B_{\vec{h}}(x,g)&\leq& d\big(1\vee\|w_\ell\|_{\infty,\bR^d}\big)^{d}\sup_{J\in\mJ}\sup_{j=1,\ldots d}M_{J}\big[b_{\vec{h},j}\big](x)
\end{eqnarray*}
Since the right hand side of the first inequality is independent of $\vec{\eta}$ we obtain
for any $\vec{h}\in\mS^{\text{const}}_d$
\begin{equation}
\label{eq5:proof-key-lemma}
B^*_{\vec{h}}(x,g)\leq \Upsilon_1\sup_{J\in\mJ}\sup_{j=1,\ldots d}M_{J}\big[b_{\vec{h},j}\big](x),\quad \forall x\in\bR^d.
\end{equation}
In particular, it yields together with (\ref{eq004:proof-key-lemma}) for any $ x\in\bR^d$
\begin{equation}
\label{eq004-new:proof-key-lemma}
B^*_{\vec{\blh}_g}(x,g)+ \mathbf{a}\Upsilon_2\e\Phi_\e\big( V_{\vec{\blh}_g}(x)\big)\leq
\inf_{\vec{h}\in\mH_\e(\vartheta,p)}\Big[\Upsilon_1b^*_{\vec{h}}(x,g)+ \mathbf{a}\Upsilon_2\e\Phi_\e\big( V_{\vec{h}}(x)\big)\Big]+2\mathbf{c}_1\e^2.
\end{equation}

\smallskip

$4^{0}.\;$ To get (\ref{eq1:proof-key-lemma}) let us first prove that for any $(\vartheta,p)\in\cP^{\text{consist}}$ one can find $\e(\vartheta,p)>0$ such that for any $\e<\e(\vartheta,p)$
\begin{equation}
\label{eq501:proof-key-lemma}
\vec{\blh}_g\in \bH_d\big(1/(2d), 3+\sqrt{2b}\big),\quad\forall g\in\bN^*_{\vec{r},d}\big(\vec{\beta},\mathbf{a}\vec{L}\big).
\end{equation}
For any $\mathbf{s}\in\bN^*$   recall that $\vec{\mh}_\mathbf{s}=(\mh_{s_1},\ldots,\mh_{s_d})$ and
$V_{\mathbf{s}}=\prod_{j=1}^d\mh_{s_j}$. Denote  $\cS^d=\{\mathbf{s}(m),\; m=0,\ldots,\widetilde{\mathbf{m}}\}$ and remark that $\Lambda_{\mathbf{s}}\big[\vec{\blh}_g\big]:=\left\{x\in (-b,b)^d:\;\vec{\blh}_g(x)=\vec{\mh}_\mathbf{s}\right\}=\emptyset$ for any $\mathbf{s}\in\bN^d,\; \mathbf{s}\neq \cS^d,$
in view of the definition of $\vec{\blh}_g$.


Taking into account   (\ref{eq3:proof-key-lemma}) we have for any
$1\leq m\leq \widetilde{\mathbf{m}}$ in view of the definition $\vec{\blh}_g$ and the second assertion of Lemma \ref{lem:set-of-band}
\begin{eqnarray*}
\Lambda_{\mathbf{s}(m)}\big[\vec{\blh}_g\big]
&\subseteq&\Big\{x\in (-b,b)^d: \widetilde{B}^*_{\vec{\mh}_{\mathbf{s}(m-1)}}(x,g)+ \mathbf{a}\Upsilon_2\e \Phi_\e\big( V_{\mathbf{s}(m-1)}\big)\geq
\widetilde{B}^*_{\vec{\mh}_{\mathbf{s}(m)}}(x,g)+ \mathbf{a}\Upsilon_2\e \Phi_\e\big( V_{\mathbf{s}(m)}\big) \Big\}
\nonumber\\
&\subseteq&\Big\{x\in (-b,b)^d:\;
\mathbf{c}_3\e^2+B^*_{\vec{\mh}_{\mathbf{s}(m-1)}}(x,g)\geq  \mathbf{a}\Upsilon_2\big[\e\Phi_\e\big( V_{\mathbf{s}(m)}\big) -\e \Phi_\e\big( V_{\mathbf{s}(m-1)}\big)\big] \Big\}.
\end{eqnarray*}
To get the last inclusion we have taken into account (\ref{eq02:proof-key-lemma}). The definition of $\mathbf{s}(m)$ implies that
\begin{eqnarray}
\label{eq5020:proof-key-lemma}
e^{-d}\prod_{j=1}^d\mh_{s_j(m)}\leq \prod_{j=1}^d\bar{\bleta}_j(m)=e^{-2d}L_\beta^{-1}\varphi^{\frac{1}{\beta}}e^{-4dm}\leq\prod_{j=1}^d\mh_{s_j(m)},\quad\forall m=0,\ldots,\widetilde{\mathbf{m}}
\end{eqnarray}
and, therefore,
$
 V^{-1}_{\mathbf{s}(m-1)} V_{\mathbf{s}(m)}\leq e^{-3d}.
$
It yields
$$
\Phi_\e\big( V_{\mathbf{s}(m)}\big) - \Phi_\e\big( V_{\mathbf{s}(m-1)}\big)\geq
2^{-1} \Phi_\e\big( V_{\mathbf{s}}\big)
$$
for any $\e>0$ small enough.

Putting  $\mathbf{c}_2=\mathbf{a}(2b+1)\mathbf{C}\big(\vec{r}\big)\|w\|_{1,\bR^d}\big(1-e^{-\beta_*}\big)^{-1}$ and using  (\ref{eq5:proof-key-lemma}) we have for any $\e>0$ provided $\e<\mathbf{c}^{-1}_1\mathbf{c}_2$
\begin{eqnarray}
\label{eq51:proof-key-lemma}
\Lambda_{\mathbf{s}(m)}\big[\vec{\blh}_g\big]
&\subseteq&\Big\{x\in (-b,b)^d:\;
\mathbf{c}_1\e^2+B^*_{\vec{\mh}_{\mathbf{s}(m-1)}}(x,g)\geq 2^{-1} \mathbf{a}\Upsilon_2\e\Phi_\e\big( V_{\mathbf{s}}\big) \Big\}
\nonumber\\
&\subseteq&\bigcup_{J\in\mJ}\bigcup_{j=1}^d\Big\{x\in (-b,b)^d:\; \e^{-1} \Phi^{-1}_\e\big( V_{\mathbf{s}}\big)
M_{J}\big[b_{\vec{\mh}_{\mathbf{s}(m-1)},j}\big](x)> \mathbf{c}_2\Big\}.
\end{eqnarray}
Here we have also used  that $\Upsilon_1\geq 1$ as well as  $\Phi_\e\big( V_{\mathbf{s}}\big)>1$ for any $\mathbf{s}\in\cS^d$.

Introduce  $\cJ_\infty=\{j=1,\ldots, d:\; r_j=\infty\}$ and recall that
 $\cJ_{\pm}=\{1,\ldots,d\}\setminus\cJ_\infty$.
In view of
(\ref{eq:strong-max-partial}) and the bound (\ref{eq:bias-norms-1}) of Lemma \ref{lem:bias-norm-bound}  with $\mathbf{r}=\infty$ and $\vec{M}=\mathbf{a}\vec{L}$ we obtain for any $j\in \cJ_\infty$ and any $J\in\mJ$
\begin{eqnarray}
\label{eq52:proof-key-lemma}
\Big\|M_{J}\big[b_{\vec{\mh}_{\mathbf{s}(m-1)},j}\big]\Big\|_{\infty,\bR^d}\leq \mathbf{c}_2L_j\mh^{\beta_j}_{s_j(m-1)}\leq \mathbf{c}_2L_j\bar{\bleta}^{\beta_j}_j(m-1)\leq \mathbf{c}_2\varphi e^{2d(m-1)}.
\end{eqnarray}
Here we have used that $\bar{\bleta}^{\beta_j}_j(m-1)=e^{-2}\big(L_j^{-1}\varphi\big)^{1/\beta_j}
e^{2d(m-1)}$ if $r_j=\infty$.

Set $\mu_\e=1$ if $\kappa(p)>0$ and $\mu_\e=\sqrt{|\ln(\e)|}$ if $\kappa(p)\leq 0$. We obtain for any $m=0,\ldots,\widetilde{\mathbf{m}}$ in view of (\ref{eq5020:proof-key-lemma}) and the definition of $\varphi$
$$
(\e \mu_\e)^{-1} \sqrt{V_{\mathbf{s}(m)}}\varphi e^{2d(m-1)}\leq e^{-\frac{5d}{2}}.
$$
Moreover, we obviously have that $\Phi\big(V_{\mathbf{s}}\big) \geq V^{-\frac{1}{2}}_{\mathbf{s}}\mu_\e$ for any $\mathbf{s}\in\bN^d$. Thus, we have
\begin{eqnarray}
\label{eq53:proof-key-lemma}
&&\e^{-1}\Phi^{-1}_\e\big( V_{\mathbf{s}(m)}\big) \varphi e^{2d(m-1)}\leq e^{-\frac{5d}{2}}
\end{eqnarray}
and, therefore, for any $j\in J_\infty$ and any $J\in\mJ$
$$
\Big\|M_{J}\big[b_{\vec{\mh}_{\mathbf{s}(m-1)},j}\big]\Big\|_{\infty,\bR^d}\leq \mathbf{c}_2 e^{-\frac{5d}{2}}<\mathbf{c}_2.
$$
It yields together with (\ref{eq51:proof-key-lemma})
\begin{eqnarray}
\label{eq54:proof-key-lemma}
\Lambda_{\mathbf{s}(m)}\big[\vec{\blh}_g\big]
&\subseteq&\bigcup_{J\in\mJ}\bigcup_{j\in\cJ_{\pm}}^d\Big\{x\in (-b,b)^d:\; \e^{-1} \Phi^{-1}_\e\big( V_{\mathbf{s}}\big)
M_{J}\big[b_{\vec{\mh}_{\mathbf{s}(m-1)},j}\big](x)> \mathbf{c}_2\Big\}.
\end{eqnarray}
We remark also that if $J_\pm=\emptyset$ then only $\Lambda_{\mathbf{s}(0)}\big[\vec{\blh}_g\big]\neq\emptyset$.
Let us consider now separately two cases.

\smallskip

$4^{0}a.\;$ Suppose  that either $\kappa(p)>0$ or $\kappa(p)\leq 0, \tau(p^*)\leq 0$ and remind that $\bar{\bleta}_j(m)=\widetilde{\bleta}(m), j=1,\ldots,d,$ for all values of $m$.
Applying the Markov inequality we get for any $m=1,\ldots \widetilde{\mathbf{m}}$  in view of
(\ref{eq:strong-max-partial}) and the bound (\ref{eq:bias-norms-1}) of Lemma \ref{lem:bias-norm-bound}  with $\mathbf{r}=r_j$ and $\vec{M}=\mathbf{a}\vec{L}$
\begin{eqnarray}
\label{eq61:proof-key-lemma}
2^{-d}\nu_d\Big(\Lambda_{\mathbf{s}(m)}\big[\vec{\blh}_g\big]\Big)&\leq&
\sum_{j\in\cJ_\pm}^d \big[\mathbf{c}_2\e\Phi_\e\big( V_{\mathbf{s}}\big)\big]^{-r_j}\big\|b_{\vec{\mh}_{\mathbf{s}(m-1)},j}\big\|^{r_j}_{r_j,\bR^d}
\nonumber\\
&\leq& \sum_{j\in\cJ_\pm}^d \big[\mathbf{c}_2\e\Phi_\e\big( V_{\mathbf{s}}\big)\big]^{-r_j}\Big(\mathbf{c}_2L_j \mh_{s_j(m-1)}^{\beta_j}\Big)^{r_j}
\nonumber\\
&\leq& \sum_{j\in\cJ_\pm}^d \Big[\e^{-1}\Phi^{-1}_\e\big( V_{\mathbf{s}}\big)L_j\widetilde{\bleta}^{\beta_j}_j(m-1)\Big]^{r_j}
\nonumber\\
&\leq& \sum_{j\in\cJ_\pm}^d \Big[\e^{-1}\Phi^{-1}_\e\big( V_{\mathbf{s}}\big)\varphi e^{2d(m-1)}\Big]^{r_j}e^{-2d\omega(2+1/\beta)(m-1)}.
\end{eqnarray}

Taking into account that $\omega\geq\beta$, we obtain in view of (\ref{eq53:proof-key-lemma})  for any $\e<\mathbf{c}^{-1}_1\mathbf{c}_2$
\begin{eqnarray}
\label{eq6100:proof-key-lemma}
&&\nu_d\Big(\Lambda_{\mathbf{s}(m)}\big[\vec{\blh}_g\big]\Big)\leq de^{-\frac{d}{2}}e^{-2d\omega(2+1/\beta)(m-1)}\leq de^{-\frac{d}{2}}e^{-2d(m-1)},\quad \forall m=1,\ldots,\widetilde{\mathbf{m}}.
\end{eqnarray}
Remembering that $\Lambda_{\mathbf{s}}\big[\vec{\blh}_g\big]=\emptyset$ for any $\mathbf{s}\notin \cS_d$ and that
$\nu_d\big(\Lambda_{\mathbf{s}(0)}\big[\vec{\blh}_g\big]\big)\leq (2b)^{d}$, taking into account the second assertion of Lemma \ref{lem:set-of-band},
we obtain putting $\tau=(2d)^{-1}$ for any  $\e<\mathbf{c}^{-1}_1\mathbf{c}_2$
\begin{eqnarray*}
\sum_{\mathbf{s}\in\bN^d}\nu^{\tau}_d\Big(\Lambda_{\mathbf{s}}\big[\vec{\blh}_g\big]\Big)=\sum_{m=1}^{\widetilde{\mathbf{m}}} \nu^{\tau}_d\Big(\Lambda_{\mathbf{s}(m)}\big[\vec{\blh}_g\big]\Big)+(2b)^{d\tau}\leq d^{\frac{1}{2d}} \big(1-e^{-1}\big)^{-1}+\sqrt{2b}\leq 2+\sqrt{2b}.
\end{eqnarray*}
Here we have used that $\sup_{d\geq 1} d^{\frac{1}{2d}} \big(1-e^{-1}\big)^{-1}<2$.

Thus, we assert that (\ref{eq501:proof-key-lemma}) is established
if either $\kappa(p)>0$ or $\kappa(p)\leq 0, \tau(p^*)\leq 0$.

\smallskip

 $4^{0}b.\;$ Let now $\kappa(p)\leq 0, \tau(p^*)> 0$.
  Since $\bar{\bleta}_j(m)=\widetilde{\bleta}(m), j=1,\ldots,d,$ if $m=0,\ldots,\widehat{\mathbf{m}}$, (\ref{eq6100:proof-key-lemma}) remains true for any $m=0,\ldots,\widehat{\mathbf{m}}$.
Similarly to (\ref{eq61:proof-key-lemma}) we obtain for any $m>\widehat{\mathbf{m}} $ in view of
(\ref{eq:strong-max-partial}) and the bound (\ref{eq:bias-norms-2}) of Lemma \ref{lem:bias-norm-bound}  with $\vec{M}=\mathbf{a}\vec{L}$
\begin{eqnarray}
\label{eq62:proof-key-lemma}
2^{-d}\nu_d\Big(\Lambda_{\mathbf{s}(m)}\big[\vec{\blh}_g\big]\Big)&\leq&
\sum_{j\in\cJ_\pm}^d \big[\mathbf{c}_2\e\Phi_\e\big( V_{\mathbf{s}}\big)\big]^{-q_j}\big\|b_{\vec{\mh}_{\mathbf{s}(m-1)},j}\big\|^{q_j}_{q_j,\bR^d}
\nonumber\\
&\leq& \sum_{j\in\cJ_\pm}^d \big[\mathbf{c}_2\e\Phi_\e\big( V_{\mathbf{s}}\big)\big]^{-q_j}\Big(\mathbf{c}_3(1-e^{-\gamma_j})^{-1}L_j \mh_{s_j(m-1)}^{\gamma_j}\Big)^{q_j}
\nonumber\\
&\leq& \sum_{j\in\cJ_\pm}^d \Big[\mathbf{c}_4(1-e^{-\gamma_j})^{-1}\e^{-1}\Phi^{-1}_\e\big( V_{\mathbf{s}}\big)L_j\widehat{\bleta}^{\gamma_j}_j(m-1)\Big]^{q_j},
\end{eqnarray}
where we have put $\mathbf{c}_3=(1-e^{\beta_*})\mathbf{c}_2$ and $\mathbf{c}_4=(1-e^{\beta_*})$.

Using (\ref{eq53:proof-key-lemma}) we get
$$
\Big[\e^{-1}\Phi^{-1}_\e\big( V_{\mathbf{s}}\big)L_j\widehat{\bleta}^{\gamma_j}_j(m-1)\Big]^{q_j}\leq e^{-\frac{3d p_\pm}{2}}e^{-2d\upsilon(2+1/\gamma)(m-1)}
\left[\frac{L_\gamma\varphi^{1/\beta}}{L_\beta\varphi^{1/\gamma}}\right]^{\upsilon}
$$
Moreover, the definition of $\widehat{\mathbf{m}}$ implies that
$$
e^{-2d\upsilon(2+1/\gamma)\widehat{\mathbf{m}}}\varphi^{(1/\beta-1/\gamma)\upsilon}
\leq e^d\big(L_\beta/L_\gamma\big)^{-\frac{\upsilon(2+1/\gamma)}{2\beta\omega\tau(2)(1/\gamma-1/\beta)}}
\varphi^{\frac{\upsilon(2+1/\gamma)}{2\beta\omega\tau(2)}-\upsilon(1/\gamma-1/\beta)}.
$$
Below we prove (see, formulae (\ref{eq00449:proof-key-lemma})), that
$
\upsilon(2+1/\gamma)-\omega(2+1/\beta)
=2\beta\tau(2)\omega\upsilon\big(1/\gamma-1/\beta\big),
$
and we obtain (recall that $\tau(2)>0$ in the considered case, see, e.g. proof of Theorem \ref{th:no-consitency})
$$
e^{-2d\upsilon(2+1/\gamma)\widehat{\mathbf{m}}}\varphi^{(1/\beta-1/\gamma)\upsilon}\leq e^d\big(L_\beta/L_\gamma\big)^{-\frac{\upsilon(2+1/\gamma)}{2\beta\omega\tau(2)(1/\gamma-1/\beta)}}
\varphi^{\frac{2+1/\beta}{2\beta\tau(2)}}\to 0,\;\e\to 0.
$$
The latter bound together with (\ref{eq62:proof-key-lemma}) yields for any $m\geq \widehat{\mathbf{m}}+1$
\begin{eqnarray}
\label{eq63:proof-key-lemma}
\nu_d\Big(\Lambda_{\mathbf{s}(m)}\big[\vec{\blh}_g\big]\Big)\leq e^{-\frac{3d}{2}}\big(L_\beta/L_\gamma\big)^{\frac{2+1/\beta}{2(\beta/\gamma-1)\tau(2)}}
\varphi^{\frac{2+1/\beta}{2\beta\tau(2)}}e^{-2d\upsilon(2+1/\gamma)(m-\widehat{\mathbf{m}}-1)}
\end{eqnarray}
and, therefore, putting $\tau=(2d)^{-1}$, we can assert that one can find $\e(\vartheta,p)>0$ such that for any $\e<\e(\vartheta,p)$
$$
\sum_{m=\widehat{\mathbf{m}}+1}^{\widetilde{\mathbf{m}}} \nu^{\tau}_d\Big(\Lambda_{\mathbf{s}(m)}\big[\vec{\blh}_g\big]\Big)\leq 1.
$$
It yields together with (\ref{eq6100:proof-key-lemma}) for all $\e<\min\{\mathbf{c}^{-1}_1\mathbf{c}_2,\e(\vartheta,p)\}$
\begin{eqnarray*}
\sum_{\mathbf{s}\in\bN^d}\nu^{\tau}_d\Big(\Lambda_{\mathbf{s}}\big[\vec{\blh}_g\big]\Big)\leq\sum_{m=1}^{\widetilde{\mathbf{m}}} \nu^{\tau}_d\Big(\Lambda_{\mathbf{s}(m)}\big[\vec{\blh}_g\big]\Big)+(2b)^{d\tau}+1\leq 3+\sqrt{2b}.
\end{eqnarray*}
Thus, we assert that (\ref{eq501:proof-key-lemma}) is established
in the case  $\kappa(p)\leq 0, \tau(p^*)>0$ as well.

\smallskip

 $5^{0}.\;$ To get (\ref{eq1:proof-key-lemma}) it remains to prove  that for any $(\vartheta,p)\in\cP^{\text{consist}}$ one can find $\e(\vartheta,p)>0$ such that for any $\e<\e(\vartheta,p)$
 \begin{equation}
\label{eq11:proof-key-lemma}
\vec{\blh}_g\in\bB(\cA_\e),\quad\forall g\in\bN^*_{\vec{r},d}\big(\vec{\beta},\mathbf{a}\vec{L}\big).
\end{equation}
The proof of (\ref{eq11:proof-key-lemma}) is mostly based on the choice of $\cA_\e$ given in (\ref{eq1:choice-of-parameters}) which, in its turn,  guarantees
(\ref{eq2:choice-of-parameters}).
We will consider separately 2 cases.

\smallskip

$5^{0}a.\;$ Let   $\kappa(p)>0$.
Obviously we can find $\mathfrak{r}\in \bN^*_p$ such that $\mathfrak{p}:=\frac{p\mathfrak{r}}{\mathfrak{r}-p}$ satisfies
$
\kappa(\mathfrak{p})>0
$
and we  have in view of (\ref{eq5020:proof-key-lemma})
\begin{eqnarray}
\label{eq110:proof-key-lemma}
&&\Big\|V^{-\frac{1}{2}}_{\vec{\blh}_g}\Big\|_{\mathfrak{p}}^{\mathfrak{p}}\leq e^{d\mathfrak{p}}L_\beta^{\mathfrak{p}/2}\varphi^{-\frac{\mathfrak{p}}{2\beta}}\bigg[(2b)^d + \sum_{m=1}^{\mathbf{m}}e^{2\mathfrak{p}dm}\nu_d\Big(\Lambda_{\mathbf{s}(m)}\big[\vec{\blh}_g\big]\Big)\bigg].
\end{eqnarray}
Using  the first bound established in (\ref{eq6100:proof-key-lemma}) we obtain
\begin{eqnarray}
\label{eq111:proof-key-lemma}
&&\Big\|V^{-\frac{1}{2}}_{\vec{\blh}_g}\Big\|_{\mathfrak{p}}^{\mathfrak{p}}\leq e^{d\mathfrak{p}}L_\beta^{\mathfrak{p}/2}\varphi^{-\frac{\mathfrak{p}}{2\beta}}\bigg[(2b)^d +d2^{-\frac{d}{2}}e^{-2d\omega(2+1/\beta)} \sum_{m=1}^{\mathbf{m}}e^{2(\mathfrak{p}-\omega(2+1/\beta))dm}\bigg].
\end{eqnarray}
Remembering that $\mathfrak{p}-\omega(2+1/\beta)=:-\kappa(\mathfrak{p})<0$ and that $\varphi^{\frac{\mathfrak{p}}{2\beta}}\cA_\e\to\infty$, we assert that
there exists $\e(\theta,p)>0$ such that
$
\big\|V^{-\frac{1}{2}}_{\vec{\blh}_g}\big\|_{\mathfrak{p}}^{\mathfrak{p}}\leq \cA_\e
$
 for any $\e<\e(\theta,p)$.

 Thus, (\ref{eq11:proof-key-lemma}) is proved if $\kappa(p)>0$. Moreover, since the right hand side of the inequality (\ref{eq111:proof-key-lemma}) as well as the choice of $\mathfrak{r}$ is independent of $g$ we can assert that for all  $\e<\e(\theta,p)$
 $$
 \mathfrak{r}\in\bN_p^*\big(\vec{\blh}_g,\cA_\e\big),\quad \forall g\in\bN^*_{\vec{r},d}\big(\vec{\beta},\mathbf{a}\vec{L}\big),
 $$
and the  assertion 2$(\mathbf{ii})$ of the proposition follows.

\smallskip

$5^{0}b.\;$ In all other cases the set $\mH_\e(\vartheta,p)$ is finite and  we obviously have
in view of (\ref{eq5020:proof-key-lemma})
\begin{eqnarray}
&&\Big\|V^{-\frac{1}{2}}_{\vec{\blh}_g}\Big\|_{t}\leq (2b)^d e^{d}L_\beta^{1/2}\varphi^{-\frac{1}{2\beta}}e^{2d\widetilde{\mathbf{m}}},\quad\forall t\geq 1.
\end{eqnarray}
This, together with the definition of $\widetilde{\mathbf{m}}$ and $\varphi$ implies that the right hand side of the latter inequality increases to infinity polynomially in $\e^{-1}$. Thus, there exists $\e(\theta,p)>0$ such that
$
\big\|V^{-\frac{1}{2}}_{\vec{\blh}_g}\big\|_{\mathfrak{p}}^{\mathfrak{p}}\leq \cA_\e
$
 for any $\e<\e(\theta,p)$ and  (\ref{eq11:proof-key-lemma}) follows.

\smallskip

$6^0.\;$
 Thus, (\ref{eq1:proof-key-lemma}) follows from  (\ref{eq501:proof-key-lemma}) and (\ref{eq11:proof-key-lemma}) and it yields together with (\ref{eq05:proof-key-lemma}) that $\mS_\e^*(\vartheta,p)\subset\bH_\e(R)$ for all $\e>0$ small enough and, therefore, the first assertion of the proposition is proved.
 We note that the  assertion 2$(\mathbf{i})$ of the proposition follows from (\ref{eq004-new:proof-key-lemma}) for any $\e>0$ such that $2\mathbf{c}_1\e\leq 1$.
 Recall, at last, that in view of (\ref{eq2:proof-key-lemma}) any $\vec{\blh}\in\mS_\e^*(\vartheta,p)$ takes values in $\mH_\e(\vartheta,p)$.
\epr

\subsection{{Proof of the theorem. Case $p\in(1,\infty)$}}
\label{sec:subsec-proof of the theorem-finite-p}

We will need  some technical results presented in Lemmas \ref{lem:technical} and \ref{lem:technical-new}. Whose proofs are postponed to Appendix.

Recall that the quantity $\cB^{(p)}_{\vec{h}}(\cdot)$ is defined in (\ref{eq:def-bias-norm})
with $K$ given in (\ref{eq:w-function}). Also, furthermore  $\mathbf{a}=\|K\|_{1,\bR^d}=\|w_\ell\|^d_{1,\bR}$.

\begin{lemma}
\label{lem:technical}
For any   $(\vartheta,p)\in\cP^{\text{consist}}$
 and any $\bH\subseteq\mS_d$
$$
\sup_{f\in\bN_{\vec{r},d}\big(\vec{\beta},\vec{L}\big)}\inf_{\vec{h}\in\bH}\Big[\cB^{(p)}_{\vec{h}}(f)+\e\Psi_{\e,p}\big(\vec{h}\big)\Big]\leq
\sup_{g\in\bN^*_{\vec{r},d}\big(\vec{\beta},\mathbf{a}\vec{L}\big)}\inf_{\vec{h}\in\bH}\Big[\cB^{(p)}_{\vec{h}}(g)+\e\Psi_{\e,p}\big(\vec{h}\big)\Big].
$$
\end{lemma}

For any $x\in (-b,b)^d$ and any $g\in\bN^*_{\vec{r},d}\big(\vec{\beta},\mathbf{a}\vec{L}\big)$ define
 $$
 U_{\vartheta,p}(x,g)=\inf_{\vec{h}\in \mH_{\e}(\vartheta,p)}\Big[ b^*_{\vec{h}}(x,g)+\varpi_\e V^{-\frac{1}{2}}_{\vec{h}}\Big],
 $$
where  $\varpi_\e=\e$ if $\kappa(p)>0$ and $\varpi_\e=\e\sqrt{|\ln(\e)|}$ if $\kappa(p)\leq 0$.
\begin{lemma}
\label{lem:technical-new}
For any   $(\vartheta,p)\in\cP$ provided $p^*>p$ and any $\e>0$ small enough
$$
\sup_{g\in\bN^*_{\vec{r},d}\big(\vec{\beta},\mathbf{a}\vec{L}\big)}\big\|U_{\vartheta,p}(\cdot,g)\big\|_{p^*}\leq \Upsilon_3 L^*,
$$
where, recall, $L^*=\min_{j: r_j=p^*} L_j$ and $\Upsilon_3=\mathbf{a}d2^d\mathbf{C}_{p^*}\big(\mathbf{C}_{p^*}\|w_\ell\|_{\infty,\bR^d}+1\big)+1$.
\end{lemma}

Let $(\vartheta,p)\in\cP^{\text{consist}}$ be fixed.
Later on  $R=3+\sqrt{2b}$ and  without further mentioning we will assume that $\e>0$ is sufficiently small in order to provide the results of Proposition
\ref{prop:key-prop}, Lemmas \ref{lem:set-of-band} and \ref{lem:technical-new}. Set also
$$
V_p\big(\vec{L}\big)=(L_\gamma/L_\beta)^{\frac{p-\omega(2+1/\beta)}{2p\beta\omega\tau(2)(1/\gamma-1/\beta)}}L^{\frac{\tau(p)}{2\tau(2)}}_\beta,\;\;p<\infty;
\quad V_\infty\big(\vec{L}\big)=L_\gamma.
$$


\subsubsection{\textsf{Proof of the theorem. Preliminaries}}
We deduce from Theorem \ref{th:L_p-norm-oracle-inequality} and Lemma \ref{lem:technical}
\begin{equation*}
\cR:=\mathbf{c}_5^{-1}\sup_{f\in\bN_{\vec{r},d}\big(\vec{\beta},\vec{L}\big)}\cR^{(p)}_\e\big[\hat{f}^{(R)}_{\vec{\mathbf{h}}}; f\big]\leq \sup_{g\in\bN^*_{\vec{r},d}\big(\vec{\beta},\mathbf{a}\vec{L}\big)}\inf_{\vec{h}\in\bH_\e(R)}\bigg\{ \cB^{(p)}_{\vec{h}}(g)+\e\Psi_{\e,p}\big(\vec{h}\big)\bigg\}+\e.
\end{equation*}
In view of the first assertion of  Proposition \ref{prop:key-prop}   $\mS_\e^*(\theta,p)\subset \bH_\e(R)$.
\begin{eqnarray}
\label{eq1:proof-th:adaptive-upper-bound}
\cR&\leq& \sup_{g\in\bN^*_{\vec{r},d}\big(\vec{\beta},\mathbf{a}\vec{L}\big)}\inf_{\vec{h}\in\mS_\e^*(\theta,p)}\Big\{ \cB^{(p)}_{\vec{h}}(g)+\e\Psi_{\e,p}\big(\vec{h}\big)\Big\}+\e,
\end{eqnarray}
 Note also the following obvious inequality: for any $p\geq 1$ and any $g:\bR^d\to\bR$
\begin{equation*}
\sup_{\vec{\eta}\in\mS_d}\big\|B_{\vec{h},\vec{\eta}}(\cdot,g)\big\|_p\leq \bigg\|\sup_{\eta\in\mS^{\text{const}}_d}B_{\vec{h},\vec{\eta}}(\cdot,g)\bigg\|_p,\quad\forall \vec{h}\in\mS_d.
\end{equation*}
This yields, in particular, for any $\vec{h}\in\mS_d$ and any $g:\bR^d\to\bR$
\begin{equation}
\label{eq2:proof-th:adaptive-upper-bound}
\cB^{(p)}_{\vec{h}}(g)\leq 2\big\|B^*_{\vec{h}}(\cdot,g)\big\|_{p}.
\end{equation}
Combining (\ref{eq1:proof-th:adaptive-upper-bound}) and (\ref{eq2:proof-th:adaptive-upper-bound}) we get
\begin{equation}
\label{eq3:proof-th:adaptive-upper-bound}
\cR\leq \sup_{g\in\bN^*_{\vec{r},d}\big(\vec{\beta},\mathbf{a}\vec{L}\big)}\Big\{2\big\|B^*_{\vec{\blh}_g}(\cdot,g)\big\|_{p} +\e\Psi_{\e,p}\big(\vec{\blh}_g\big)\Big\}+\e,
\end{equation}
where $\vec{\blh}_g$ satisfies the second  assertion of Proposition \ref{prop:key-prop}.
Consider separately two cases.

\paragraph{Case $\kappa(p)>0$}

 Recall that $\vec{\blh}_g(x)$  takes values in $\mH_\e(\vartheta,p)$ for any $g\in\bN^*_{\vec{r},d}\big(\vec{\beta},\mathbf{a}\vec{L}\big)$ and $x\in(-b,b)^d$. Additionally, $\mH_\e(\vartheta,p)\subset \mH^d(\mh_\e)$
 in view of the first assertion of Lemma \ref{lem:set-of-band} since $\kappa(p)>0$.

 It implies  $\vec{\blh}_g\in\mS_d(\mh_\e)$  and we can assert that
 $$
 \Psi_{\e,p}\big(\vec{\blh}_g\big)\leq \inf_{r\in\bN^*_p(\vec{\blh}_g,\cA_\e)}C_2(r)\Big\| V^{-\frac{1}{2}}_{\vec{\blh}_g}\Big\|_{\frac{rp}{r-p}},\quad \forall
 g\in\bN^*_{\vec{r},d}\big(\vec{\beta},\mathbf{a}\vec{L}\big).
 $$
Applying the assertion 2($\mathbf{i}\mathbf{i}$) of Proposition \ref{prop:key-prop} we can state that
 for some $\mathfrak{r}$ provided $\kappa\big(\frac{\mathfrak{r}p}{\mathfrak{r}-p}\big)>0$
$$
\Psi_{\e,p}\big(\vec{\blh}_g\big)\leq  C_2(\mathfrak{r})\big\| V^{-\frac{1}{2}}_{\vec{\blh}}\big\|_{\frac{\mathfrak{r}p}{\mathfrak{r}-p}},\quad \forall
 g\in\bN^*_{\vec{r},d}\big(\vec{\beta},\mathbf{a}\vec{L}\big).
$$
Denoted  $\mmp=\frac{\mathfrak{r}p}{\mathfrak{r}-p}$. Since $\mmp>p$ in view of H\"older inequality
 $\big\|B^*_{\vec{\blh}_g}(\cdot,g)\big\|_{p}\leq (2b)^d\big\|B^*_{\vec{\blh}_g}(\cdot,g)\big\|_{\mmp}$ and we deduce from (\ref{eq3:proof-th:adaptive-upper-bound}) (remembering that we consider here the norms of positive functions) that
\begin{eqnarray}
\label{eq4:proof-th:adaptive-upper-bound}
\cR&\leq& \mathbf{c}_6\sup_{g\in\bN^*_{\vec{r},d}\big(\vec{\beta},\mathbf{a}\vec{L}\big)}\Big\{\big\|B^*_{\vec{\blh}_g}(\cdot,g)\big\|_{\mmp} +\e\big\| V^{-\frac{1}{2}}_{\vec{\blh}_g}\big\|_{\mathfrak{p}}\Big\}+\e
\nonumber\\
&\leq&\mathbf{c}_7\sup_{g\in\bN^*_{\vec{r},d}\big(\vec{\beta},\mathbf{a}\vec{L}\big)}\Big\|B^*_{\vec{\blh}_g}(\cdot,g)+ \mathbf{a}\Upsilon_2\e V^{-\frac{1}{2}}_{\vec{\blh}_g}(\cdot)\Big\|_{\mathfrak{p}}+\e.
\end{eqnarray}
Applying the  assertion 2($\mathbf{i}$) of  Proposition \ref{prop:key-prop} we obtain
\begin{equation}
\label{eq5:proof-th:adaptive-upper-bound}
\cR\leq \mathbf{c}_8\sup_{g\in\bN^*_{\vec{r},d}\big(\vec{\beta},\mathbf{a}\vec{L}\big)}\bigg\|\inf_{\vec{h}\in\mH_\e(\vartheta,p)}\Big[ b^*_{\vec{h}}(\cdot,g)+\e \Phi_\e\big( V_{\vec{h}}\big)\Big]\bigg\|_{\mathfrak{p}}+\e.
\end{equation}

\paragraph{Case $\kappa(p)\leq 0$}
Since
$
\Psi_{\e,p}\big(\vec{h}\big)\leq\bigg(C_1\Big\|\sqrt{\big|\ln{\big(\e V_{\vec{h}}\big)}\big|} V^{-\frac{1}{2}}_{\vec{h}}\Big\|_p\bigg)
$
for any  $\vec{h}\in\bB(\cA_\e)$, we deduce from (\ref{eq3:proof-th:adaptive-upper-bound}) similarly to  (\ref{eq4:proof-th:adaptive-upper-bound})
$$
\cR\leq \mathbf{c}_{9}\sup_{g\in\bN^*_{\vec{r},d}\big(\vec{\beta},\mathbf{a}\vec{L}\big)}\Big\|B^*_{\vec{\blh}_g}(\cdot,g)+ \mathbf{a}\Upsilon_2\e \sqrt{\big|\ln{\big(\e V_{\vec{\blh}_g}(\cdot)\big)}\big|} V^{-\frac{1}{2}}_{\vec{\blh}_g}(\cdot)\Big\|_{p}+\e
$$
Applying the first assertion 2($\mathbf{i}$) of  Proposition \ref{prop:key-prop} we have
\begin{equation*}
\cR\leq \mathbf{c}_{10}\sup_{g\in\bN^*_{\vec{r},d}\big(\vec{\beta},\mathbf{a}\vec{L}\big)}\bigg\|\inf_{\vec{h}\in\mH_\e(\vartheta,p)}\Big[ b^*_{\vec{h}}(\cdot,g)+\e \Phi_\e\big( V_{\vec{h}}\big)\Big]\bigg\|_{p}+\e.
\end{equation*}
This together with (\ref{eq5:proof-th:adaptive-upper-bound}) allows us to assert that
\begin{eqnarray}
\label{eq6:proof-th:adaptive-upper-bound}
\cR&\leq& \mathbf{c}_{11}\sup_{g\in\bN^*_{\vec{r},d}\big(\vec{\beta},\mathbf{a}\vec{L}\big)}\bigg\|\inf_{\vec{h}\in\mH_{\e}(\vartheta,p)}\Big[ b^*_{\vec{h}}(\cdot,g)+\e \Phi_\e\big( V_{\vec{h}}\big)\Big]\bigg\|_{\mathbf{p}}+\e,
\end{eqnarray}
where we have denoted $\mathbf{p}=\mathfrak{p}$ if $\kappa(p)> 0$ and $\mathbf{p}=p$ if $\kappa(p)\leq 0$.

The definition of  $\widetilde{\mathbf{m}}$ allows us to assert that if $\kappa(p)\leq 0$
$$
\Phi_\e\big( V_{\vec{h}}\big)\Big]\leq \mathbf{c}_{12}\sqrt{\big|\ln{(\e)}\big|} V^{-\frac{1}{2}}_{\vec{h}},\quad\forall \vec{h}\in \mH_{\e}(\vartheta,p).
$$
Hence we get from (\ref{eq6:proof-th:adaptive-upper-bound})
\begin{eqnarray}
\label{eq7:proof-th:adaptive-upper-bound}
\cR&\leq& \mathbf{c}_{13}\sup_{g\in\bN^*_{\vec{r},d}\big(\vec{\beta},\mathbf{a}\vec{L}\big)}\bigg\|\inf_{\vec{h}\in \mH_{\e}(\vartheta,p)}\Big[ b^*_{\vec{h}}(\cdot,g)+\varpi_\e V^{-\frac{1}{2}}_{\vec{h}}\Big]\bigg\|_{\mathbf{p}}+\e
\nonumber\\
&=&\mathbf{c}_{13}\sup_{g\in\bN^*_{\vec{r},d}\big(\vec{\beta},\mathbf{a}\vec{L}\big)}\big\|U_{\vartheta,p}(\cdot,g)\big\|_{\mathbf{p}}+\e
=:\mathbf{c}_{13}\sup_{g\in\bN^*_{\vec{r},d}\big(\vec{\beta},\mathbf{a}\vec{L}\big)}\cR_{\mathbf{p}}(g)+\e,
\end{eqnarray}
where, recall, $\varpi_\e=\e$ if $\kappa(p)>0$ and $\varpi_\e=\e\sqrt{|\ln(\e)|}$ if $\kappa(p)\leq 0$.

\subsubsection{\textsf{Proof of the theorem. Slicing}} For any $g\in\bN^*_{\vec{r},d}\big(\vec{\beta},\mathbf{a}\vec{L}\big)$ we have
\begin{eqnarray}
\label{eq8:proof-th:adaptive-upper-bound}
\cR^{\mathbf{p}}_{\mathbf{p}}(g)&\leq& (2b)^d(Q\varphi)^{\mathbf{p}}+ \sum_{m=0}^{\widetilde{\mathbf{m}}} (Qe^{2d(m+1)}\varphi)^{\mathbf{p}}\nu_d\big(\Gamma_m\big)
+\int_{\Gamma_{\mathbf{m}}}\big|U_{\vartheta,p}(x,g)\big|^{\mathbf{p}}\nu_d(\rd x)
\nonumber\\
&=:&(2b)^d(Q\varphi)^{\mathbf{p}}+(Qe^{2d}\varphi)^{\mathbf{p}}\sum_{m=1}^{\widetilde{\mathbf{m}}} e^{2dm\mathbf{p}}
\nu_d\big(\Gamma_m\big)+T_{\widetilde{\mathbf{m}}}.
\end{eqnarray}
Here we have put $\Gamma_m=\big\{x\in (b,b)^d:\; U_{\vartheta,p}(x,g)\geq Qe^{2dm}\varphi\big\}$ and  $Q=2\mathbf{c_2}+e^{d}$, where, recall, $\mathbf{c}_2=\mathbf{a}(2b+1)\mathbf{C}\big(\vec{r}\big)\|w\|_{1,\bR^d}\big(1-e^{-\beta_*}\big)^{-1} $. Moreover, if $\widetilde{\mathbf{m}}=\infty$
we set $T_\infty=0$.

\smallskip

We have in view of the definition of $U_{\vartheta,p}$
\begin{eqnarray*}
\label{eq9:proof-th:adaptive-upper-bound}
\Gamma_m\subset\Big\{x\in (b,b)^d:\; b^*_{\vec{\mh}_{\mathbf{s}(m)}}(\cdot,g)+\varpi_\e V^{-\frac{1}{2}}_{\mathbf{s}(m)}\geq Qe^{2dm}\varphi\Big\}
\end{eqnarray*}
Recall that in view of (\ref{eq5020:proof-key-lemma})
$
e^{-d}V_{\mathbf{s}(m)}\leq e^{-2d}L_\beta^{-1}\varphi^{\frac{1}{\beta}}e^{-4dm}\leq V_{\mathbf{s}(m)}
$
for any $m=0,\ldots,\widetilde{\mathbf{m}}$.

Hence,
$
\varpi_\e V^{-\frac{1}{2}}_{\mathbf{s}(m)}\leq e^d\varphi e^{2dm}
$
and we get
\begin{eqnarray}
\label{eq10:proof-th:adaptive-upper-bound}
&&\Gamma_m\subset\Big\{x\in (b,b)^d:\; b^*_{\vec{\mh}_{\mathbf{s}(m)}}(\cdot,g)\geq 2\mathbf{c_2}e^{2dm}\varphi\Big\}=:\Gamma_m^*
\end{eqnarray}
Note also that
$
 \varpi_\e V^{-\frac{1}{2}}_{\mathbf{s}(m)}\leq e^d\varphi e^{2dm}< e^d\mathbf{c_2}e^{2dm}\varphi<  e^d b^*_{\vec{\mh}_{\mathbf{s}(m)}}(x,g)
$
 for any $x\in\Gamma_m^*$ that yields,
\begin{eqnarray}
\label{eq12:proof-th:adaptive-upper-bound}
&& |U_{\vartheta,p}(x,g)\big|\leq  ( e^d+1)b^*_{\vec{\mh}_{\mathbf{s}(m)}}(x,g), \quad\forall x\in\Gamma_m^*.
\end{eqnarray}
The latter inequality allows us to  bound from above  $T_{\widetilde{\mathbf{m}}}$ if $\widetilde{\mathbf{m}}<\infty$. Indeed, in view of (\ref{eq12:proof-th:adaptive-upper-bound})
\begin{eqnarray}
\label{eq13:proof-th:adaptive-upper-bound}
&& T_{\widetilde{\mathbf{m}}}\leq ( e^d+1)^{\mathbf{p}}\big\|b^*_{\vec{\mh}_{\mathbf{s}(\widetilde{\mathbf{m}})}}(\cdot,g)\big\|^{\mathbf{p}}_{\mathbf{p}}.
\end{eqnarray}

Another bound can be obtained in the case $p^*>\mathbf{p}$. Applying H\"older inequality, the assertion of Lemma \ref{lem:technical-new}
and (\ref{eq9:proof-th:adaptive-upper-bound})
\begin{eqnarray}
\label{eq14:proof-th:adaptive-upper-bound}
&& T_{\widetilde{\mathbf{m}}}\leq (L^*\Upsilon_3)^{\mathbf{p}}\big[\nu_d\big(\Gamma^*_{\widetilde{\mathbf{m}}}\big)\big]^{1-\mathbf{p}/p^*}.
\end{eqnarray}
The definition of $b^*_{\vec{\mh}_{\mathbf{s}(\mathbf{m})}}(\cdot,g)$ implies that for any $m=1,\ldots \widetilde{\mathbf{m}}$
\begin{eqnarray*}
&&\Gamma_m^*\subset \bigcup_{J\in\mJ}\bigcup_{j=1}^d\Big\{x\in (b,b)^d:\; M_{J}\big[b_{\vec{\mh}_{\mathbf{s}(m),j}}\big](x)\geq 2\mathbf{c_2}e^{2dm}\varphi\Big\}.
\end{eqnarray*}
Since in view of (\ref{eq52:proof-key-lemma})
$
\Big\|M_{J}\big[b_{\vec{\mh}_{\mathbf{s}(m)},j}\big]\Big\|_{\infty,\bR^d}\leq \mathbf{c}_2\varphi e^{2dm}
$
for any $j\in \cJ_\infty$ and any $J\in\mJ$ we obtain
\begin{eqnarray}
\label{eq15:proof-th:adaptive-upper-bound}
&&\Gamma_m^*\subset \bigcup_{J\in\mJ}\bigcup_{j\in\cJ_\pm}\Big\{x\in (b,b)^d:\; M_{J}\big[b_{\vec{\mh}_{\mathbf{s}(m),j}}\big](x)\geq \mathbf{c_2}e^{2dm}\varphi\Big\},\quad \forall m=1,\ldots \widetilde{\mathbf{m}}.
\end{eqnarray}
If either $\kappa(p)>0$ or $\kappa(p)\leq 0,\;\tau(p^*)\leq 0$ the following bound is true.
\begin{eqnarray}
\label{eq16:proof-th:adaptive-upper-bound}
\nu_d\big(\Gamma_m^*\big)&\leq& \mathbf{c}_{14}e^{-2dm\omega(2+1/\beta)},\quad \forall m=1,\ldots \widetilde{\mathbf{m}} .
\end{eqnarray}
Indeed, applying the Markov inequality we get for any $m=1,\ldots\widetilde{\mathbf{m}}$ in view of
(\ref{eq:strong-max-partial}) and the bound (\ref{eq:bias-norms-1}) of Lemma \ref{lem:bias-norm-bound}  with $\mathbf{r}=r_j$ and $\vec{M}=\mathbf{a}\vec{L}$
\begin{eqnarray*}
2^{-d}\nu_d\big(\Gamma_m^*\big)&\leq&
\sum_{j\in\cJ_\pm}^d \big[\mathbf{c}_2e^{2dm}\varphi\big]^{-r_j}\big\|b_{\vec{\mh}_{\mathbf{s}(m)},j}\big\|^{r_j}_{r_j,\bR^d}
\leq\sum_{j\in\cJ_\pm}^d \big[e^{2dm}\varphi\big]^{-r_j}\Big(L_j \mh_{s_j(m)}^{\beta_j}\Big)^{r_j}
\nonumber\\
&\leq& \sum_{j\in\cJ_\pm}^d \Big[e^{-2dm}\varphi^{-1} L_j\widetilde{\bleta}^{\beta_j}_j(m)\Big]^{r_j}
\leq de^{-2dm\omega(2+1/\beta)}.
\end{eqnarray*}
If $\kappa(p)\leq 0,\;\tau(p^*)> 0$ we have for any $\widehat{\mathbf{m}}<m\leq \widetilde{\mathbf{m}}$
\begin{eqnarray}
\label{eq17:proof-th:adaptive-upper-bound}
\nu_d\big(\Gamma_m^*\big)
\leq \mathbf{c}_{16}(L_\gamma/L_\beta)^{\upsilon}\varphi^{(1/\beta-1/\gamma)\upsilon}e^{-2dm\upsilon(2+1/\gamma)}.
\end{eqnarray}
Indeed, we obtain for any $m>\widehat{\mathbf{m}} $ in view of
(\ref{eq:strong-max-partial}) and the bound (\ref{eq:bias-norms-2}) of Lemma \ref{lem:bias-norm-bound}  with $\vec{M}=\mathbf{a}\vec{L}$
\begin{eqnarray*}
2^{-d}\nu_d\big(\Gamma_m^*\big)&\leq&
\sum_{j\in\cJ_\pm}^d \big[\mathbf{c}_2e^{2dm}\varphi\big]^{-q_j}\big\|b_{\vec{\mh}_{\mathbf{s}(m)},j}\big\|^{q_j}_{q_j,\bR^d}
\leq \mathbf{c}_{15}\sum_{j\in\cJ_\pm}^d \big[e^{2dm}\varphi\big]^{-q_j}\Big(L_j \mh_{s_j(m)}^{\gamma_j}\Big)^{q_j}
\nonumber\\
&\leq& \mathbf{c}_{15}\sum_{j\in\cJ_\pm}^d \Big[e^{-2dm}\varphi^{-1} L_j\widehat{\bleta}^{\gamma_j}_j(m)\Big]^{q_j}
\leq \mathbf{c}_{16}(L_\gamma/L_\beta)^{\upsilon}\varphi^{(1/\beta-1/\gamma)\upsilon}e^{-2dm\upsilon(2+1/\gamma)}.
\end{eqnarray*}

\subsubsection{\textsf{Proof of the theorem. Derivation of rates}} We will proceed differently in depending on a zone to which the pair $(\vartheta,p)\in\cP^{\text{consist}}$ belongs.

\paragraph{Dense zone: $\kappa(p)>0$. Case $\kappa(p^*)\geq 0$}  Recall that $\widetilde{\mathbf{m}}=\infty$ in this case and, therefore $T_{\widetilde{\mathbf{m}}}=0$. Moreover $\mathbf{p}=\mmp$.  We get from  (\ref{eq8:proof-th:adaptive-upper-bound}) and  (\ref{eq16:proof-th:adaptive-upper-bound})
\begin{eqnarray}
\label{eq19:proof-th:adaptive-upper-bound}
\cR^{\mmp}_{\mmp}(g)\leq \mathbf{c}_{17}\varphi^{\mmp}\sum_{m=0}^{\infty} e^{2dm(\mmp-\omega(2+1/\beta))}\leq \mathbf{c}_{18}\varphi^{\mmp}=
\mathbf{c}_{18}\bde_\e^{\ma\mmp},
\end{eqnarray}
since $\mmp-\omega(2+1/\beta)=-\kappa(\mmp)<0$ in view of the definition of $\mmp$. Taking into account that the right hand side of the latter inequality
is independent of $g$ we obtain in view of (\ref{eq7:proof-th:adaptive-upper-bound})
\begin{eqnarray}
\label{eq20:proof-th:adaptive-upper-bound}
\cR\leq \mathbf{c}_{19}\big(\bde_\e\big)^{\ma}.
\end{eqnarray}
\paragraph{ Dense zone: $\kappa(p)>0$. Case $\kappa(p^*)<0$} Note first that  $\kappa(p^*)<0$ and  $\kappa(\mmp)>0$ implies $\mmp<p^*$ since
$\kappa(\cdot)$ is decreasing. We get in view of (\ref{eq14:proof-th:adaptive-upper-bound}) and (\ref{eq16:proof-th:adaptive-upper-bound})
\begin{eqnarray}
\label{eq21:proof-th:adaptive-upper-bound}
T_{\widetilde{\mathbf{m}}}\leq \mathbf{c}_{20} (L^*)^{\mmp}e^{-2d\widetilde{\mathbf{m}}\omega(2+1/\beta)(1-\mmp/p^*)}\leq
 \mathbf{c}_{21} (L^*)^{\mmp} (\mh_\e)^{\frac{\ell p^*}{\kappa(p^*)}}\varphi^{-\frac{\omega(2+1/\beta)(1-\mmp/p^*)}{\kappa(p^*)/p^*}}
\end{eqnarray}
Noting that
$
\mmp+\frac{\omega(2+1/\beta)(1-\mmp/p^*)}{\kappa(p^*)/p^*}=\frac{p^*\kappa(\mmp)}{\kappa(p^*)}<0
$
and taking into account (\ref{eq2:choice-of-parameters}), we obtain that
$$
\varphi^{-\mmp}(\mh_\e)^{\frac{\ell p^*}{\kappa(p^*)}}T_{\widetilde{\mathbf{m}}}\to 0,\;\e\to 0.
$$
Since (\ref{eq19:proof-th:adaptive-upper-bound}) holds  we assert finally that (\ref{eq20:proof-th:adaptive-upper-bound}) remains true if
$\kappa(p^*)<0$ as well. Thus, the theorem is proved in the case $\kappa(p)>0$.

\paragraph{ New zone: $\kappa(p)\leq 0$, $\tau(p^*)\leq 0$} Recall that $\mathbf{p}=p$ and necessarily $p^*>p$ since we consider $(\vartheta,p)\in \cP^{\text{consist}}$.

Noting that the first inequality in (\ref{eq21:proof-th:adaptive-upper-bound}) remains true we deduce from (\ref{eq8:proof-th:adaptive-upper-bound}) and  (\ref{eq16:proof-th:adaptive-upper-bound})
\begin{eqnarray}
\label{eq22:proof-th:adaptive-upper-bound}
\cR^{p}_{p}(g)\leq \mathbf{c}_{22}\varphi^{p}\sum_{m=0}^{\widetilde{\mathbf{m}}} e^{2dm(p-\omega(2+1/\beta))}+\mathbf{c}_{20} (L^*)^{p}e^{-2d\widetilde{\mathbf{m}}\omega(2+1/\beta)(1-p/p^*)}.
\end{eqnarray}
If $\kappa(p)<0$ we have
\begin{eqnarray}
\label{eq23:proof-th:adaptive-upper-bound}
\cR^{p}_{p}(g)&\leq& \mathbf{c}_{23}\varphi^{p} e^{2d\widetilde{\mathbf{m}}(p-\omega(2+1/\beta))}+\mathbf{c}_{20} (L^*)^{p}e^{-2d\widetilde{\mathbf{m}}\omega(2+1/\beta)(1-p/p^*)}
\nonumber\\
&\leq&\mathbf{c}_{24}(1+L^*)^{p}\varphi^{\frac{pp^*\omega(2+1/\beta)(1/p-1/p^*)}{p^*-\omega(2+1/\beta)}}
=\mathbf{c}_{24}(1+L^*)^{p}\bde_\e^{\frac{pp^*\omega(1/p-1/p^*)}{p^*-\omega(2+1/\beta)}}
\end{eqnarray}
in view of the definition of $\widetilde{\mathbf{m}}$.
Taking into account that the right hand side of the latter inequality
is independent of $g$ we obtain in view of (\ref{eq7:proof-th:adaptive-upper-bound})
\begin{eqnarray}
\label{eq24:proof-th:adaptive-upper-bound}
\cR\leq \mathbf{c}_{19}(1+L^*)\big(\bde_\e\big)^{\ma}.
\end{eqnarray}
If $\kappa(p)=0$ we deduce from (\ref{eq22:proof-th:adaptive-upper-bound}) and the definition of $\widetilde{\mathbf{m}}$
$$
\cR^{p}_{p}(g)\leq \mathbf{c}_{23}\varphi^{p}\widetilde{\mathbf{m}}+(L^*)^{p}\varphi^{\frac{pp^*\omega(2+1/\beta)(1/p-1/p^*)}{p^*-\omega(2+1/\beta)}}\leq
\mathbf{c}_{24}\varphi^{p}|\ln(\e)|+(L^*)^{p}\varphi^p.
$$
Here we have used that   $\frac{p^*\omega(2+1/\beta)(1/p-1/p^*)}{p^*-\omega(2+1/\beta)}=\frac{\beta}{2\beta+1}$ if $\kappa(p)=0$. Thus,  we conclude
\begin{eqnarray}
\label{eq240:proof-th:adaptive-upper-bound}
\cR\leq \mathbf{c}_{25}\big(\bde_\e\big)^{\ma}|\ln(\e)|^{\frac{1}{p}}.
\end{eqnarray}

Thus, the theorem is proved in the case $\kappa(p)\leq 0$.

\paragraph{ Sparse zone: $\kappa(p)\leq 0$, $\tau(p^*)> 0$}

If $p^*=p$ taking into account that $\widetilde{\mathbf{m}}=\widehat{\mathbf{m}}+1$ we deduce from
(\ref{eq8:proof-th:adaptive-upper-bound}), (\ref{eq13:proof-th:adaptive-upper-bound}), (\ref{eq16:proof-th:adaptive-upper-bound}) and (\ref{eq17:proof-th:adaptive-upper-bound})
\begin{eqnarray*}
\label{eq241:proof-th:adaptive-upper-bound}
\cR^{p}_{p}(g)&\leq& \mathbf{c}_{18}\varphi^{p}\sum_{m=0}^{\widehat{\mathbf{m}}} e^{2dm(p-\omega(2+1/\beta))}+( e^d+1)^{p}\big\|b^*_{\vec{\mh}_{\mathbf{s}(\widehat{\mathbf{m}}+1)}}(\cdot,g)\big\|^p_{p}.
\end{eqnarray*}
Using  (\ref{eq13:proof-th:adaptive-upper-bound}), (\ref{eq:strong-max-partial}) and  the triangle inequality we have
\begin{eqnarray*}
\|b^*_{\vec{\mh}_{\mathbf{s}(\widehat{\mathbf{m}}+1)}}(\cdot,g)\big\|_{p}\leq
 \sum_{J\in\mJ}\sum_{j=1}^d \big\|M_{J}\big[b_{\vec{\mh}_{\mathbf{s}(\widehat{\mathbf{m}}+1),j}}\big\|_{p}
\leq \mathbf{c}_{26}\sum_{j=1}^d \big\|b_{\vec{\mh}_{\mathbf{s}(\widehat{\mathbf{m}}+1)},j}\big\|_{p,\bR^d}.
\end{eqnarray*}
Note that $p^*=p$ implies $p_\pm=p$ and, therefore,  $q_j=p$ for any $j=1,\ldots d$, where, recall,  $q_j$ are given in (\ref{eq:gamma-and-q}). Hence we obtain using the bound (\ref{eq:bias-norms-2}) of Lemma \ref{lem:bias-norm-bound}  with $\vec{M}=\mathbf{a}\vec{L}$
\begin{eqnarray}
\label{eq28:proof-th:adaptive-upper-bound}
\|b^*_{\vec{\mh}_{\mathbf{s}(\widehat{\mathbf{m}}+1)}}(\cdot,g)\big\|_{p}&\leq&  \mathbf{c}_{27}\sum_{j=1}^dL_j\mh^{\gamma_j}_{s_j(\widehat{\mathbf{m}}+1)}\leq  \mathbf{c}_{27}\sum_{j=1}^dL_j\widehat{\bleta}^{\gamma_j}_j(\widehat{\mathbf{m}}+1)
\nonumber\\
&\leq& \mathbf{c}_{28}(L_\gamma/L_\beta)^{\upsilon/p}\varphi^{1+(1/\beta-1/\gamma)(\upsilon/p)} e^{2d\widehat{\mathbf{m}}(1-(\upsilon/p)(2+1/\gamma))}
\nonumber\\
&\leq&\mathbf{c}_{28}(L_\gamma/L_\beta)^{\frac{\omega(2+1/\beta)-p}{2p\beta\omega\tau(2)(1/\gamma-1/\beta)}}\varphi^{\frac{(2+1/\beta)\tau(p)}{2\tau(2)}}
\nonumber\\
&=&\mathbf{c}_{28} (L_\gamma/L_\beta)^{\frac{\omega(2+1/\beta)-p}{2p\beta\omega\tau(2)(1/\gamma-1/\beta)}}L^{\frac{\tau(p)}{2\tau(2)}}_\beta
(\e^2|\ln(\e)|)^{\frac{\tau(p)}{2\tau(2)}}=\mathbf{c}_{28} \bde_\e^{\ma}.
\end{eqnarray}
We get in view of the definition of $\widehat{\mathbf{m}}$
\begin{eqnarray}
\label{eq280:proof-th:adaptive-upper-bound}
\varphi^{p}\sum_{m=0}^{\widehat{\mathbf{m}}} e^{2dm(p-\omega(2+1/\beta))}&\leq&\mathbf{c}_{29}\varphi^{p} e^{2d\widehat{\mathbf{m}}(p-\omega(2+1/\beta))}
\leq \mathbf{c}_{30}\bde_\e^{\ma p},\quad\;\;\;\kappa(p)<0;
\\
\label{eq2800:proof-th:adaptive-upper-bound}
\varphi^{p}\sum_{m=0}^{\widehat{\mathbf{m}}} e^{2dm(p-\omega(2+1/\beta))}&\leq&\mathbf{c}_{29}\varphi^{p}(\widehat{\mathbf{m}}+1)
\leq \mathbf{c}_{29}\varphi^{p}|\ln(\e)|, \qquad\quad  \kappa(p)=0.
\end{eqnarray}
Therefore, if $\kappa(p)<0$
\begin{eqnarray}
\label{eq281:proof-th:adaptive-upper-bound}
&&\cR^{p}_{p}(g)\leq \mathbf{c}_{31}\bde_\e^{\ma p}.
\end{eqnarray}
If $\kappa(p)=0$ we can easily check that $\frac{\tau(p)}{2\tau(2)}=\frac{\beta}{2\beta+1}$ that yields in view of the definition
of $\widehat{\mathbf{m}}$
\begin{eqnarray}
\label{eq282:proof-th:adaptive-upper-bound}
&&\cR^{p}_{p}(g)\leq \mathbf{c}_{32}\varphi^{p}|\ln(\e)|=\mathbf{c}_{33}\bde_\e^{\ma p}|\ln(\e)|.
\end{eqnarray}
Taking into account that the right hand sides in (\ref{eq281:proof-th:adaptive-upper-bound}) and (\ref{eq282:proof-th:adaptive-upper-bound})
are independent of $g$ we obtain in view of (\ref{eq7:proof-th:adaptive-upper-bound})
\begin{eqnarray}
\label{eq283:proof-th:adaptive-upper-bound}
\cR\leq \mathbf{c}_{34}\bde_\e^{\ma },\;\; \kappa(p)<0; \qquad \cR\leq\mathbf{c}_{29}\bde_\e^{\ma}|\ln(\e)|,\;\; \kappa(p)=0..
\end{eqnarray}
This completes the proof of the theorem in the case $\kappa(p)\leq 0$, $\tau(p^*)> 0$, $p^*=p$.


If $p^*>p$  we deduce from
(\ref{eq8:proof-th:adaptive-upper-bound}),  (\ref{eq16:proof-th:adaptive-upper-bound}) and (\ref{eq17:proof-th:adaptive-upper-bound})
\begin{eqnarray*}
\label{eq25:proof-th:adaptive-upper-bound}
\cR^{p}_{p}(g)&\leq& \mathbf{c}_{18}\varphi^{p}\sum_{m=0}^{\widehat{\mathbf{m}}} e^{2dm(p-\omega(2+1/\beta))}
\nonumber\\
&&+\mathbf{c}_{16}(L_\gamma/L_\beta)^{\upsilon}\varphi^{p+(1/\beta-1/\gamma)\upsilon}\sum_{m=\widehat{\mathbf{m}}+1}^{\widetilde{\mathbf{m}}} e^{2dm(p-\upsilon(2+1/\gamma))}+T_{\widetilde{\mathbf{m}}}
\nonumber\\
&\leq& \mathbf{c}_{35}\Big[\varphi^{p}\sum_{m=0}^{\widehat{\mathbf{m}}} e^{2dm(p-\omega(2+1/\beta))}+(L_\gamma/L_\beta)^{\upsilon}\varphi^{p+(1/\beta-1/\gamma)\upsilon}
e^{2d\widehat{\mathbf{m}}(p-\upsilon(2+1/\gamma))}\Big]+T_{\widetilde{\mathbf{m}}}.
\end{eqnarray*}
Here we have used that $p\leq p_\pm<\upsilon(2+1/\gamma)<0$ in view of (\ref{eq901:proof-key-lemma}).

Using (\ref{eq00449:proof-key-lemma}) and the definition of $\widehat{\mathbf{m}}$ we compute that
\begin{eqnarray*}
(L_\gamma/L_\beta)^{\upsilon}\varphi^{p+(1/\beta-1/\gamma)\upsilon}
e^{2d\widehat{\mathbf{m}}(p-\upsilon(2+1/\gamma))}&\leq& \mathbf{c}_{36} (L_\gamma/L_\beta)^{\frac{\omega(2+1/\beta)-p}{2\beta\omega\tau(2)(1/\gamma-1/\beta)}} \varphi^{\frac{p(2+1/\beta)\tau(p)}{2\tau(2)}}
\\
&=&\mathbf{c}_{36}V^p_p\big(\vec{L}\big)
(\e^2|\ln(\e)|)^{\frac{p\tau(p)}{2\tau(2)}}=\mathbf{c}_{36} \bde_\e^{\ma p}.
\end{eqnarray*}
It yields together with (\ref{eq280:proof-th:adaptive-upper-bound}) and (\ref{eq2800:proof-th:adaptive-upper-bound})
\begin{eqnarray*}
\label{eq26:proof-th:adaptive-upper-bound}
&&\cR^{p}_{p}(g)\leq \mathbf{c}_{37}\bde_\e^{\ma p}+T_{\widetilde{\mathbf{m}}},\;\; \kappa(p)<0; \qquad \cR^{p}_{p}(g)\leq\mathbf{c}_{29}\bde_\e^{\ma p}|\ln(\e)|+T_{\widetilde{\mathbf{m}}},\;\; \kappa(p)=0.
\end{eqnarray*}
Using  (\ref{eq14:proof-th:adaptive-upper-bound}) and (\ref{eq17:proof-th:adaptive-upper-bound})
\begin{eqnarray*}
\label{eq29:proof-th:adaptive-upper-bound}
T_{\widetilde{\mathbf{m}}}&\leq& \mathbf{c}_{38} A\varphi^{(1/\beta-1/\gamma)\upsilon(1-p/p^*)}e^{-2d\widetilde{\mathbf{m}}\upsilon(2+1/\gamma)(1-p/p^*)}ç
\\
&\leq& \mathbf{c}_{38} A\varphi^{(1/\beta-1/\gamma)\upsilon(1-p/p^*)}e^{-2d\overline{\mathbf{m}}\upsilon(2+1/\gamma)(1-p/p^*)}
=\mathbf{c}_{38} A\varphi^p
\end{eqnarray*}
in view of the definition of $\widetilde{\mathbf{m}}$ and $\overline{\mathbf{m}}$. Here $A=A(L_\beta,L_\gamma,L^*)$ can be easily computed.

Thus we can assert that (\ref{eq283:proof-th:adaptive-upper-bound}) holds in the case $p^*>p$ as well that completes the proof of the theorem.
\epr

\subsection{{Proof of the theorem. Case $p\in\{1,\infty\}$}}
\label{sec:subsec-proof-of-theorem}
Note that $p_\pm=\infty$ if $p=\infty$ and therefore $\vec{\gamma}=\vec{\bga}(\infty)$.
The proof of the theorem in this case is the straightforward  consequence of the Corollary \ref{cor:th:L_p-norm-oracle-inequality}.

Introduce
the vectors $\vec{u}=(u_1,\ldots,u_d)$ and $\vec{v}=(v_1,\ldots,v_d)$ as follows.
\begin{eqnarray*}
u_j&=&L_j^{-1/\beta_j}\big(L_\beta\e^2\big)^{\frac{\beta}{\beta_j(2\beta+1)}},\quad v_j=L_j^{-1/\gamma_j}\big(L_\gamma\e^2|\ln(\e)|\big)^{\frac{\gamma}{\gamma_j(2\gamma+1)}}.
\end{eqnarray*}
It is obvious that both vectors belongs to $\bH^{\text{const}}_\e$ and without loss of generality we can assume that $\vec{u},\vec{v}\in\mH^d$.
Moreover,
$$
\e\Psi^{(\text{const})}_{\e,\infty}\big(\vec{v}\big)\leq \mathbf{c}_{38}\e\sqrt{\big|\ln{(\e)}\big|} V^{-\frac{1}{2}}_{\vec{v}}=\mathbf{c}_{38}\big(L_\gamma\e^2|\ln(\e)|\big)^{\frac{\gamma}{2\gamma+1}}.
$$
Note that
$
1/\gamma =\sum_{j=1}^d\frac{\tau(r_j)}{\beta_j\tau(\infty)}=\frac{1}{\beta\tau(\infty)}\;\;\Rightarrow\;\;\frac{\gamma}{2\gamma+1}=\frac{\tau(\infty)}{2\tau(2)}
$
and, therefore,
$$
\e\Psi^{(\text{const})}_{\e,\infty}\big(\vec{v}\big)=\mathbf{c}_{38}\bde^\ma.
$$
Additionally , we easily compute
$$
\e\Psi^{(\text{const})}_{\e,1}\big(\vec{v}\big)\leq \mathbf{c}_{39}\e V^{-\frac{1}{2}}_{\vec{v}}=\mathbf{c}_{38}\big(L_\beta\e^2\big)^{\frac{\beta}{2\beta+1}}=\mathbf{c}_{39}\bde^\ma.
$$
Applying assertions of Lemma \ref{lem:bias-norm-bound} we obtain
$$
\sum_{j=1}^d\big\|b_{\vec{u},j}\big\|_1\leq\mathbf{c}_{40}\sum_{j=1}^d L_ju^{\beta_j}_j=\mathbf{c}_{40}\bde^\ma;\quad
\sum_{j=1}^d\big\|b_{\vec{v},j}\big\|_\infty\leq\mathbf{c}_{41}\sum_{j=1}^d L_jv^{\gamma_j}_j=\mathbf{c}_{41}\bde^\ma.
$$
The assertion of the theorem follows now from  Corollary \ref{cor:th:L_p-norm-oracle-inequality}.

\epr

\section{Appendix}

\subsection{{Proof of the assertion $(\mathbf{ii})$ of Lemma \ref{lem:about-maximum}}} For any given $\mathbf{s}\in\bN^d$ and any $\vec{\mathrm{h}}\in\mathfrak{S}_d$ define
$$
\Lambda_{s_j}[\mathrm{h}_j]=\big\{x\in (-b,b)^d:\;\;\mathrm{h}_j(x)=s_j\big\},\quad j=1,\ldots, d.
$$
Then $\Lambda_{\mathbf{s}}[\vec{\mathrm{h}}]=\cap_{j=1}^d\Lambda_{s_j}[\mathrm{h}_j]$ and we get putting $\mathbf{s}_j=(s_1,\ldots,s_{j-1},s_{j+1},\ldots,d)$
$$
\nu_d\Big(\Lambda_{s_j}[\mathrm{h}_j]\Big)=\sum_{\mathbf{s}_j\in\bN^{d-1}}\nu_d\Big(\Lambda_{\mathbf{s}}[\vec{\mathrm{h}}]\Big),\quad j=1,\ldots, d.
$$
It yields
for any $\alpha\in (0,1)$ we have for any $\vec{\mathrm{h}}\in\mathfrak{S}_d$
\begin{equation*}
\label{eq1:proof-lem:about-maximum}
\sum_{s_j=1}^\infty\nu^{\frac{\alpha}{d}}_d\Big(\Lambda_{s_j}[\mathrm{h}_j]\Big)\leq
\sum_{\mathbf{s}\in\bN^d}\nu^{\frac{\alpha}{d}}_d\Big(\Lambda_{\mathbf{s}}[\vec{\mathrm{h}}]\Big).
\end{equation*}
Since obviously $\vec{h}\vee\vec{\eta}\in\mathfrak{S}_d$ for any $\vec{h},\vec{\eta}\in\mathfrak{S}_d$, we have
\begin{equation*}
\sum_{s_j=1}^\infty\nu^{\frac{\alpha}{d}}_d\Big(\Lambda_{s_j}[h_j\vee\eta_j]\Big)
\leq\sum_{s_j=1}^\infty\Big\{\nu^{\frac{\alpha}{d}}_d\Big(\Lambda_{s_j}[h_j]\Big)
+\nu^{\frac{\alpha}{d}}_d\Big(\Lambda_{s_j}[\eta_j]\Big)\Big\}.
\end{equation*}
Hence, for any $\alpha\in (0,1)$ and any $\vec{h},\vec{\eta}\in\bH_d\Big(d^{-1}\alpha,2^{-1}L^{\frac{1}{d}}\Big), L>0,$ we get
\begin{equation}
\label{eq2:proof-lem:about-maximum}
\sum_{s_j=1}^\infty\nu^{\frac{\alpha}{d}}_d\Big(\Lambda_{s_j}[h_j\vee\eta_j]\Big)\leq
\sum_{\mathbf{s}\in\bN^d}\nu^{\frac{\alpha}{d}}_d\Big(\Lambda_{\mathbf{s}}[\vec{h}]\Big)
+\sum_{\mathbf{s}\in\bN^d}\nu^{\frac{\alpha}{d}}_d\Big(\Lambda_{\mathbf{s}}[\vec{\eta}]\Big)\leq L^{\frac{1}{d}}.
\end{equation}
Note that $\Lambda_{\mathbf{s}}[\vec{\mathrm{h}}]=\cap_{j=1}^d\Lambda_{s_j}[\mathrm{h}_j]$ implies
 for any $\vec{\mathrm{h}}\in\mathfrak{S}_d$ and $\alpha\in (0,1)$
\begin{equation}
\label{eq3:proof-lem:about-maximum}
\nu^{\alpha}_d\Big(\Lambda_{\mathbf{s}}[\vec{\mathrm{h}}]\Big)\leq \prod_{j=1}^d \nu^{\frac{\alpha}{d}}_d\Big(\Lambda_{s_j}[\mathrm{h}_j]\Big).
\end{equation}
Therefore, we deduce from  (\ref{eq2:proof-lem:about-maximum}) and (\ref{eq3:proof-lem:about-maximum})
$$
\sum_{\mathbf{s}\in\bN^d}\nu^{\alpha}_d\Big(\Lambda_{\mathbf{s}}[\vec{h}\vee\vec{\eta}]
\Big)\leq  \sum_{\mathbf{s}\in\bN^d}\prod_{j=1}^d\nu^{\frac{\alpha}{d}}_d\Big(\Lambda_{s_j}[h_j\vee\eta_j]\Big)
\leq \prod_{j=1}^d \sum_{s_j=1}^\infty\nu^{\frac{\alpha}{d}}_d\Big(\Lambda_{s_j}[h_j\vee\eta_j]\Big)\leq L.
$$
Thus, we obtain that $\vec{h},\vec{\eta}\in\bH_d\Big(d^{-1}\alpha,2^{-1}L^{\frac{1}{d}}\Big)$ implies $\vec{h}\vee\vec{\eta}\in\bH_d(\alpha,L)$.
Putting $\kappa=\alpha/d$ and $\mL=2^{-1}L^{\frac{1}{d}}$  we come to
the assertion of the lemma since $\vec{h},\vec{\eta}\in\cB(\cA)$ implies
$\vec{h}\vee\vec{\eta}\in\cB(\cA)$.
\epr

\subsection{{Proof of Lemma \ref{lem:techniclal-profcor:th:L_p-norm-oracle-inequality}}}
Let $\vec{h},\vec{\eta}\in\mS^{\text{const}}_d$ be fixed. Denote by $\cJ=\big\{j=1,\ldots d:\;\; h_j\vee \eta_j=h_j\big\}$ and suppose first that $\cJ\neq\emptyset$.
Let $\cJ=\left\{j_1<j_2<\cdots<j_k\right\}$,  $k=|\cJ|$, and put $\cJ_l=\left\{j_1<j_2<\cdots<j_l\right\}$, $l=1,\ldots,k$.
Note that for any $x\in\bR^d$
\begin{eqnarray*}
B_{\vec{h},\vec{\eta}}(x,f)&=&
\bigg|\int_{\bR^d}K_{\vec{h}\vee \vec{\eta}}(t-x)\Big[f(t)-f\big(t+\mathbf{E}[\cJ](x-t)\big)\Big]\nu_d(\rd t)
\\
&&-\int_{\bR^d}K_{\vec{\eta}}(t-x)\Big[f(t)-f\big(t+\mathbf{E}[\cJ](x-t)\big)\Big]\nu_d(\rd t)\bigg|.
\end{eqnarray*}
Here we have used Assumption \ref{ass:kernel} ($\mathbf{ii}$) and $\int\cK=1$. Remark also that
\begin{gather}
\label{eq01:proof-cor1}
f(t)-f\big(t+\mathbf{E}[\cJ](x-t)\big)=\sum_{l=1}^k f\big(t+\mathbf{E}[\cJ_{l-1}](x-t)\big)-f\big(t+\mathbf{E}[\cJ_l](x-t)\big),
\end{gather}
where we have put $\cJ_0=\emptyset$.
Thus we have for any $x\in\bR^d$
\begin{eqnarray}
\label{eq06:proof-th:L_p-norm-oracle-inequality}
B_{\vec{h},\vec{\eta}}(x,f)&\leq&
\sum_{l=1}^k \bigg|\int_{\bR^d}K_{\vec{h}\vee \vec{\eta}}(t-x)\Big[f\big(t+\mathbf{E}[\cJ_{l-1}](x-t)\big)-f\big(t+\mathbf{E}[\cJ_{l}](x-t)\big)\Big]\nu_d(\rd t)\bigg|
\nonumber\\
&&+\sum_{l=1}^k\bigg|\int_{\bR^d}K_{\vec{\eta}}(t-x)\Big[f\big(t+\mathbf{E}[\cJ_{l-1}](x-t)\big)-f\big(t+\mathbf{E}[\cJ_{l}](x-t)\big)\Big]\nu_d(\rd t)\bigg|.
\end{eqnarray}
Noting that   $h_{j_l}\vee \eta_{j_l}=h_{j_l}$ for any $l=1,\ldots,k,$ in view of the definition of $\cJ$ we have
\begin{eqnarray}
\label{eq1:proof-cor1}
&&\bigg|\int_{\bR^d}K_{\vec{h}\vee \vec{\eta}}(t-x)\Big[f\big(t+\mathbf{E}[\cJ_{l-1}](x-t)\big)-f\big(t+\mathbf{E}[\cJ_{l}](x-t)\big)\Big]\nu_d(\rd t)\bigg|
\\
&\leq&\int_{\bR^{|\bar{\cJ_l}|}}\Big|K_{\vec{h}\vee\vec{\eta},\bar{\cJ_l}}\big(t_{\bar{\cJ_l}}-x_{\bar{\cJ_l}}\big)
\Big|b_{\vec{h},j_l}\big(x_{\cJ_l},t_{\bar{\cJ_l}}\big)
\nu_{|\bar{\cJ_l}|}\big(\rd t_{\bar{\cJ_l}}\big)
=\Big[\big|K_{\vec{h}\vee\vec{\eta}}\big|\star b_{\vec{h},j_l} \Big]_{\cJ_l}(x),
\nonumber \\*[2mm]
&&\bigg|\int_{\bR^d}K_{\vec{\eta}}(t-x)\Big[f\big(t+\mathbf{E}[\cJ_{l-1}](x-t)\big)-f\big(t+\mathbf{E}[\cJ_{l}](x-t)\big)\Big]\nu_d(\rd t)\bigg|
\nonumber \\
 \label{eq2:proof-cor1}
&\leq&\int_{\bR^{|\bar{\cJ_l}|}}\Big|K_{\vec{\eta},\bar{\cJ_l}}\big(t_{\bar{\cJ_l}}-x_{\bar{\cJ_l}}\big)
\Big|b_{\vec{h},j_l}\big(x_{\cJ_l},t_{\bar{\cJ_l}}\big)
\nu_{|\bar{\cJ_l}|}\big(\rd t_{\bar{\cJ_l}}\big)
=\Big[\big|K_{\vec{\eta}}\big|\star b_{\vec{h},j_l} \Big]_{\cJ_l}(x).
\nonumber
\end{eqnarray}
Here we have used once again Assumption \ref{ass:kernel} ($\mathbf{ii}$) and $\int\cK=1$.

Thus, for any $\vec{h},\vec{\eta}\in\mS^{\text{const}}_d$ for which $\cJ\neq\emptyset$ the first assertion of the lemma follows  from (\ref{eq06:proof-th:L_p-norm-oracle-inequality}), (\ref{eq1:proof-cor1}) and (\ref{eq2:proof-cor1}). It remains to note that $B_{\vec{h},\vec{\eta}}(\cdot,f)\equiv 0$ if $\cJ=\emptyset$ and, therefore, the first assertion is true with an arbitrary choice of $\{j_1,\ldots,j_k\}$.
In particular, one can choose $k=d$ that corresponds to $\{j_1,\ldots,j_k\}=\{1,\ldots,d\}$.

To get the second assertion we choose $\cJ:=\{j_1,\ldots,j_k\}=\{1,\ldots,d\}$ that yields $\cJ_l=\{1,\ldots,l\}$ and note that
(\ref{eq01:proof-cor1})  remains true. Repeating the computations led to (\ref{eq2:proof-cor1}) with $\vec{\eta}$ replaced by $\vec{h}$ we come to the second assertion of the lemma.
\epr

\subsection{{Proof of Lemma \ref{lem:nik-lep}}}

As it was already mentioned if $r^*(s)=s$ the assertion of the lemma is proved in \cite{Nikolski}, Section 6.9.
Thus, it remains to study the case $r^*>s$, where we put $r^*=\max_{j=1,\ldots d}r_j$. Set also $\vec{r}^*=(r^*,\ldots,r^*)$ and
denote $J_+=\{j: r_j\geq s\}$ and $J_-=\{1,\ldots,d\}\setminus J_+$.

The assumption  $\tau\big(r^*(s)\big)=\tau\big(r^*\big)>0$ together with $r_j\leq r^*$ for any $j=1,\ldots, d$ makes possible the application
of the theorem of Section 6.9, \cite{Nikolski} that yields
\begin{equation}
\label{eq1:embedd-nik-proof}
 \bN_{\vec{r},d}\big(\vec{\beta},\vec{L}\big) \subseteq
\bN_{\vec{r}^*,d}\big(\vec{\bga}(r^*),\mathbf{c}\vec{L}\big)
\end{equation}
Note that for any $j\in J_-$ we have $\|f\|_{r_j}\leq L_j$ since $f\in\bN_{\vec{r},d}\big(\vec{\beta},\vec{L}\big)$ and $\|f\|_{r^*}\leq \mathbf{c}L_j$
in view of (\ref{eq1:embedd-nik-proof}). Noting that $r_j<s=r_j(s)< r^*$ we have $\|f\|_{r_j(s)}\leq \mathbf{c}_1L_j$ for any $j\in\cJ_-$ in view of H\"older inequality. It remains to note that $r_j(s)=r_j$ for any $j\in J_+$ and we assert that
\begin{equation}
\label{eq2:embedd-nik-proof}
 \|f\|_{r_j(s)}\leq \mathbf{c}_1L_j, \;\;\forall j=1,\ldots d.
\end{equation}
Since $f\in\bN_{\vec{r},d}\big(\vec{\beta},\vec{L}\big)$ and $\bga_j(s)=\beta_j$,  $r_j(s)=r_j, j\in\cJ_+$ one has
\begin{equation}
\label{eq3:embedd-nik-proof}
  \Big\|\Delta_{u,j}^{k_j} g\Big\|_{r_j(s),\bR^d}= \Big\|\Delta_{u,j}^{k_j} g\Big\|_{r_j,\bR^d} \leq L_j |u|^{\beta_j}=L_j |u|^{\bga_j(s)},\;\;\;\;
\forall u\in \bR,\;\;\;\forall j\in\cJ_+.
\end{equation}
Let now $j\in J_-$. If $r^*=\infty$ we have
\begin{equation}
\label{eq4:embedd-nik-proof}
 \Big\|\Delta_{u,j}^{k_j} g\Big\|_{s,\bR^d}^s\leq \Big\|\Delta_{u,j}^{k_j} g\Big\|^{r_j}_{r_j,\bR^d}\Big\|\Delta_{u,j}^{k_j} g\Big\|^{s-r_j}_{\infty,\bR^d}
 \leq  \mathbf{c}^{s-r_j} L^{s}_j |u|^{r_j\beta_j+(s-r_j)\beta_j\tau(\infty)\tau^{-1}(r_j)},
\end{equation}
in view of (\ref{eq1:embedd-nik-proof}). If $r^*<\infty$, writing
$$
s=\frac{r_j(r^*-s)}{r^*-r_j}+\frac{r^*(s-r_j)}{r^*-r_j}
$$
and applying
the H\"older inequality with exponents $\frac{r^*-r_j}{r^*-s}$ and $\frac{r^*-r_j}{s-r_j}$ we obtain
\begin{eqnarray}
\label{eq5:embedd-nik-proof}
 \Big\|\Delta_{u,j}^{k_j} g\Big\|_{s,\bR^d}^s&\leq& \Big(\Big\|\Delta_{u,j}^{k_j} g\Big\|_{r_j,\bR^d}\Big)^{\frac{(r^*-s)r_j}{r^*-r_j}}
 \Big(\Big\|\Delta_{u,j}^{k_j} g\Big\|_{r^*,\bR^d}\Big)^{\frac{(s-r_j)r^*}{r^*-r_j}}
 \nonumber\\
 &\leq&  \mathbf{c}_1^{\frac{(s-r_j)r^*}{r^*-r_j}} L^{s}_j |u|^{a_j},\;\;\forall u\in \bR,
\end{eqnarray}
in view of (\ref{eq1:embedd-nik-proof}) with
$$
a_j=\frac{(r^*-s)\beta_jr_j}{r^*-r_j}+\frac{\bga_j(r^*)(s-r_j)r^*}{r^*-r_j}=\frac{(r^*-s)\beta_jr_j}{(r^*-r_j)}
+\frac{\tau(r^*)(s-r_j)\beta_jr^*}{\tau(r_j)(r^*-r_j)}.
$$
Note that (\ref{eq4:embedd-nik-proof}) is a particular case of (\ref{eq5:embedd-nik-proof}).

We easily compute that
$
b_j:=\tau(r_j)(r^*-s)\beta_jr_j+\tau(r^*)(s-r_j)\beta_jr^*=s\beta_j\tau(s)(r^*-r_j)
$
and, therefore,
$$
a_j:=\frac{b_j}{\tau(r_j)(r^*-r_j)}=\frac{s\tau(s)\beta_j}{\tau(r_j)}=s\bga_j(s).
$$
Thus, we obtain from (\ref{eq1:embedd-nik-proof})
\begin{eqnarray}
\label{eq6:embedd-nik-proof}
&& \Big\|\Delta_{u,j}^{k_j} g\Big\|_{s,\bR^d}
 \leq  \mathbf{c}_1 L^{s}_j |u|^{\bga_j(s)},\;\;\forall u\in \bR,\;\;\;\forall j\in\cJ_-.
\end{eqnarray}
The required embedding follows now from  (\ref{eq2:embedd-nik-proof}), (\ref{eq3:embedd-nik-proof}) and (\ref{eq6:embedd-nik-proof}).
\epr

\subsection{{Proof of Lemma \ref{lem:bias-norm-bound}}}
We obviously have
\begin{eqnarray*}
b_{\vec{\mathrm{h}},j}(x)&=&\sup_{\eta\leq \mathrm{h}_j,\;\eta\in\mH}
\bigg|\int_{\bR} w_\ell(\mathfrak{z})\big[f\big(x+\mathfrak{z}\eta\mathbf{e}_j\big)-f(x)\big]\nu_1(\rd \mathfrak{z}) \bigg|
\\
&=&\sup_{\eta\leq \mathrm{h}_j,\;\eta\in\mH}
\bigg|\int_{\bR} w_\ell(\mathfrak{z})\big[\Delta_{\mathfrak{z}\eta,j}f(x)\big]\nu_1(\rd \mathfrak{z}) \bigg|.
\end{eqnarray*}
For $j=1,\ldots, d$ we have
\begin{eqnarray*}
&&\int_{\bR} w_\ell(\mathfrak{z}) \Delta_{\mathfrak{z}\eta, j}f(x)\nu_1(\rd \mathfrak{z}) =
\int_{\bR} \sum_{i=1}^\ell \binom{\ell}{i} (-1)^{i+1}\frac{1}{i} w\Big(\frac{\mathfrak{z}}{i}\Big)\big[\Delta_{\eta\mathfrak{z}, j}f(x)\big]
\nu_1(\rd \mathfrak{z})
\nonumber
\\
&&= (-1)^{\ell-1}\int_{\bR} w(z) \sum_{i=1}^\ell \binom{\ell}{i} (-1)^{i+\ell}\big[\Delta_{iz\eta, j}f(x)\big]
\nu_1(\rd z)
= (-1)^{\ell-1} \int_{\bR} w(z) \big[\Delta^\ell_{z\eta, j}\, f(x)\big]\nu_1(\rd z).
\label{eq:int-representation}
\end{eqnarray*}
The last equality follows from the definition of $\ell$-th order difference operator
(\ref{eq:Delta}).
Thus, for any $j=1,\ldots,d$ and any $ x\in(-b,b)^d$
\begin{eqnarray}
\label{eq:B-h-j}
&&b_{\vec{\mathrm{h}},j}(x,f)=\sup_{\eta\leq \mathrm{h}_j,\;\eta\in\mH}\bigg|\int_{\bR} w(z) \big[\Delta^\ell_{z\eta, j}\, f(x)\big] \nu_1(\rd z)\bigg|
\leq\sum_{\eta\leq \mathrm{h}_j}\bigg|\int_{\bR} w(z) \big[\Delta^\ell_{z\eta, j}\, f(x)\big] \nu_1(\rd z)\bigg|,
\end{eqnarray}
since  $\mH$ is a discrete set.
Therefore, by the Minkowski inequality for integrals [see, e.g., \cite[Section~6.3]{folland}] and the triangle inequality, choosing
$\mathfrak{s}$ from the relation $ e^{-\mathfrak{s}-2}=\mathrm{h}_j$ (recall that $\mathrm{h}_j\in\mH$) we obtain
\begin{eqnarray*}
&& \big\|b_{\vec{\mathrm{h}},j}(\cdot,f)\big\|_{\mathbf{r},\bR^d}
\leq\sum_{s=\mathfrak{s}}^\infty
  \int_{-1/(2\ell)}^{1/(2\ell)} |w(z)|\big\| \Delta^\ell_{ze^{-s-2}, j}\, f\big\|_{\mathbf{r},\bR^d}\;\nu_1(\rd z) .
\end{eqnarray*}
Here we have also used that $w$ is compactly supported on $[-1/(2\ell),1/(2\ell)]$.

Note that $ \Delta^\ell_{ze^{-s-2}, j}\, f$ is  supported on $\cY:=(-b-1/2,b+1/2)^d$ for any  $z\in[-1/(2\ell),1/(2\ell)]$. Hence, taking into account that $\mathbf{r}\leq r_j$ we get
\begin{eqnarray*}
\big\| \Delta^\ell_{ze^{-s-2}, j}\, f\big\|_{\mathbf{r},\bR^d}&=&\big\| \Delta^\ell_{ze^{-s-2}, j}\, f\big\|_{\mathbf{r},\cY}\leq (2b+1)^{d(\frac{1}{\mathbf{r}}-\frac{1}{r_j})}
\big\| \Delta^\ell_{ze^{-s-2}, j}\, f\big\|_{r_j,\cY}
\\
&\leq&(2b+1)^d
\big\| \Delta^\ell_{ze^{-s-2}, j}\, f\big\|_{r_j,\bR^d}\leq (2b+1)M_j(ze^{-s-2})^{\beta_j},
\end{eqnarray*}
since $f\in \bN_{\vec{r},d}\big(\vec{\beta},\vec{M}\big)$. Hence, for any $\mathbf{r}\in [1,r_j]$
\begin{eqnarray*}
\big\|b_{\vec{\mathrm{h}},j}(\cdot,f)\big\|_{\mathbf{r},\bR^d}&\leq&(2b+1)^d
M_j   \int_{-1/(2\ell)}^{1/(2\ell)} |w(z)|\,|z|^{\beta_j} \nu_1(\rd z) \sum_{s=\mathfrak{s}}^\infty \big(e^{-s-2}\big)^{\beta_j}
\\
&\leq& (2b+1)^d\|w\|_{1,\bR^d}\big(1-e^{-\beta_j}\big)^{-1} M_j\mathrm{h}_j^{\beta_j}.
\end{eqnarray*}
This proves (\ref{eq:bias-norms-1}).

The inequality in (\ref{eq:bias-norms-2}) follows by the same reasoning with $r_j$ replaced by
$q_j$, $\beta_j$ replaced by $\gamma_j$
and with the use of embedding (\ref{eq:embedd-nik}).

\epr

\subsection{{Proof of Lemma \ref{lem:set-of-band}}}

We will analyze the set $\mH_\e(\vartheta,p)$ separately for different values of $(\vartheta,p)$.

\smallskip

$1^{0}.\; \text{Case}\; \kappa(p)> 0.\;$
If $\kappa(p^*)\geq 0$ we have $r_j\leq p^*\leq \omega(2+1/\beta)$ for all $j=1,\ldots d$. Therefore, for any $m\geq 0$
\begin{eqnarray*}
\label{eq0041-new:proof-key-lemma}
\widetilde{\bleta}_j(m)\leq e^{-2}\big(L_j^{-1}\varphi\big)^{1/\beta_j}, \quad \forall j=1,\ldots d.
\end{eqnarray*}
Thus, for all $\e>0$ small enough $\bar{\bleta}_j(m):=\widetilde{\bleta}_j(m)< \mh_\e$. It yields
\begin{eqnarray}
\label{eq0042-new:proof-key-lemma}
\mh_{s_j(m)}\leq \bar{\bleta}_j(m)<e\mh_{s_j(m)},\quad j=1,\ldots,d,
\end{eqnarray}
\par
If $\kappa(p^*)<0$, that is possible only if $p^*>p$ in view of $\kappa(p)>0$,
 we have for any $0\leq m\leq \widetilde{\mathbf{m}}$ and any $j=1,\ldots d$
\begin{eqnarray*}
\label{eq0043-new:proof-key-lemma}
\widetilde{\bleta}_j(m)\leq e^{-2}\big(L_j^{-1}\varphi\big)^{1/\beta_j}
e^{2dm\left(\frac{1}{\beta_j}-\frac{\omega(2+1/\beta)}{\beta_j p^*}\right)}\leq e^{-2}
\bigg(L_j^{-1}\varphi e^{2d\widetilde{\mathbf{m}}\left(-\frac{\kappa(p^*)}{p^*}\right)}\bigg)^{1/\beta_j}.
\end{eqnarray*}
The definition of $\widetilde{\mathbf{m}}$ implies
$
L_j^{-1}\varphi e^{2d\widetilde{\mathbf{m}}\left(-\frac{\kappa(p^*)}{p^*}\right)}\leq \mh^\ell_\e
$
and we  assert that for all $\e>0$ small enough
$$
\bar{\bleta}_j(m):=\widetilde{\bleta}_j(m)< \mh^{\ell/\beta_j}_\e\leq \mh_\e
$$
since $\beta_j\leq\ell$ and $\mh_\e<1$.
Thus, we conclude that for any $(\vartheta,p)\in\cP^{\text{consist}}$ provided $\kappa(p)>0$ and  for all $\e>0$ small enough
\begin{eqnarray}
\label{eq00440-new:proof-key-lemma}
\mH_\e(\vartheta,p)\subset \mH^d(\mh_\e).
\end{eqnarray}
It implies, in particular, that (\ref{eq0042-new:proof-key-lemma})
takes place when $\kappa(p^*)< 0$ as well. Hence we obtain
\begin{eqnarray}
\label{eq0045-new:proof-key-lemma}
e^{-d}\prod_{j=1}^d\mh_{s_j(m)}\leq \prod_{j=1}^d\bar{\bleta}_j(m)=e^{-2d}L_\beta^{-1}\varphi^{\frac{1}{\beta}}e^{-4dm}\leq\prod_{j=1}^d\mh_{s_j(m)}.
\end{eqnarray}
It yields
\begin{eqnarray}
\label{eq0046-new:proof-key-lemma}
\mathbf{s}(m)\neq \mathbf{s}(n),\quad \forall m\neq n,\; m,n=0,\ldots, \widetilde{\mathbf{m}}.
\end{eqnarray}

$2^{0}.\; \text{Case}\;\kappa(p)\leq 0,\; \tau(p^*)\leq 0.\;$
Since we consider $(\vartheta,p)\in\cP^{\text{consist}}$ the later case is possible only if  $p^*>p$. It implies
$\kappa(p^*)<0$, and
as previously we have
\begin{eqnarray*}
\widetilde{\bleta}_j(m)\leq e^{-2}
\bigg(L_j^{-1}\varphi e^{2d\widetilde{\mathbf{m}}
\left(1-\frac{\omega(2+1/\beta)}{p^*}\right)}\bigg)^{1/\beta_j}.
\end{eqnarray*}
It yields in view of the definition of $\widetilde{\mathbf{m}}$
\begin{eqnarray*}
\widetilde{\bleta}_j(m)
\leq e^{-2}\Big((L_0^{-1}\varphi) e^{-2d\widetilde{\mathbf{m}}
\left(\kappa(p^*)/p^*\right)}\Big)^{1/\beta_j}\leq e^{-2}.
\end{eqnarray*}
We conclude that for any $(\vartheta,p)\in\cP^{\text{consist}}$ such that $\kappa(p)\leq0$  and $\tau(p^*)\leq 0$ for all $\e>0$
\begin{eqnarray}
\label{eq00448-new:proof-key-lemma}
\mH_\e(\vartheta,p)\subset \mH^d,
\end{eqnarray}
and that  (\ref{eq0042-new:proof-key-lemma}) holds.  Hence, in view of  (\ref{eq0045-new:proof-key-lemma}) we assert that (\ref{eq0046-new:proof-key-lemma})
is also fulfilled.


\smallskip

$3^{0}.\; \text{Case}\;\kappa(p)\leq 0,\; \tau(p^*)>0.\;$  Recall that this case is possible only if $p>2$ that implies in particular that $\tau(2)>0$.

We start with presenting some relations between the parameters $\beta,\gamma,\omega$
and $\upsilon$ whose proved are given in Section \ref{sec:subsec-equalities}.
\begin{eqnarray}
\label{eq901:proof-key-lemma}
\gamma<\beta,\qquad \upsilon(2+1/\gamma)>p_\pm;
\\*[2mm]
\label{eq0044901:proof-key-lemma}
1/\omega-1/\upsilon=\beta\big(1/\gamma-1/\beta\big)\big(1-1/\omega\big).
\end{eqnarray}
We deduce from  the equality (\ref{eq0044901:proof-key-lemma})
\begin{eqnarray}
\label{eq00449:proof-key-lemma}
\upsilon(2+1/\gamma)-\omega(2+1/\beta)&=&\omega\upsilon\left[(2+1/\beta)(1/\omega-1/\upsilon)+
(1/\gamma-1/\beta)\omega^{-1}\right]
\nonumber\\
&=&2\beta\tau(2)\omega\upsilon\big(1/\gamma-1/\beta\big).
\end{eqnarray}
Using (\ref{eq00449:proof-key-lemma}) we easily get for any $\mathbf{r}>0$
\begin{eqnarray}
\label{eq00449-useful:proof-key-lemma}
&&1-\frac{\mathbf{r}-\upsilon(2+1/\gamma)}{2\mathbf{r}\beta\omega\tau(2)}-\frac{(1/\gamma-1/\beta)\upsilon}{\mathbf{r}}
=\frac{(2+1/\beta)\tau(\mathbf{r})}{2\tau(2)}.
\end{eqnarray}
Since $\upsilon(2+1/\gamma)\geq p_\pm$, $\widetilde{\mathbf{m}}>\widehat{\mathbf{m}}$  in view of the definition of $\widetilde{\mathbf{m}}$ and
  $q_j\leq p_\pm, j\in\cJ_\pm$ and $q_j=\infty, j\in\cJ_\infty$
  we get
\begin{eqnarray}
\bar{\bleta}_j(m)&=&\widetilde{\bleta}_j(m)\leq e^{-2}\big(L_j^{-1}\varphi\big)^{1/\beta_j}
e^{2d\widehat{\mathbf{m}}\big(\frac{1}{\beta_j}-\frac{\omega(2+1/\beta)}{\beta_jp^*}\big)}, \qquad\qquad\qquad\quad\;\; m\leq\widehat{\mathbf{m}};
\nonumber\\
\label{eq004490:proof-key-lemma}
\bar{\bleta}_j(m)&=&\widehat{\bleta}_j(m)\leq e^{-2}
\big(L_j^{-1}\varphi\big)^{1/\gamma_j}
e^{2d\widehat{\mathbf{m}}\left(\frac{1}{\gamma_j}-\frac{\upsilon(2+1/\gamma)}{\gamma_jq_j}\right)}
\left[\frac{L_\gamma\varphi^{1/\beta}}{L_\beta\varphi^{1/\gamma}}\right]^{\frac{\upsilon}{\gamma_jq_j}}, \quad m>\widehat{\mathbf{m}}.
\end{eqnarray}
 We obtain in view of the definition of $\widehat{\mathbf{m}}$
\begin{eqnarray*}
\big\{\bar{\bleta}_j(m)\big\}^{\beta_j}&\leq& e^{-2\beta_j}L_j^{-1}\varphi^{\frac{(2+1/\beta)\tau(p^*)}{2\tau(2)}},
\qquad\qquad\quad\quad\;\;\quad m\leq\widehat{\mathbf{m}};
\\
\big\{\bar{\bleta}_j(m)\big\}^{\gamma_j}&\leq&  e^{-2\gamma_j}T_1
L_j^{-1}\varphi^{1-\frac{q_j-\upsilon(2+1/\gamma)}{2q_j\beta\omega\tau(2)}-\frac{(1/\gamma-1/\beta)\upsilon}{q_j}}, \quad m>\widehat{\mathbf{m}},\;j\in\cJ_\pm;
\\
\big\{\bar{\bleta}_j(m)\big\}^{\gamma_j}&\leq& e^{-2\gamma_j}
L_j^{-1}T_2\varphi^{\frac{(2+1/\beta)(1-1/\omega)}{2\tau(2)}}, \qquad\qquad\quad\;\; m>\widehat{\mathbf{m}},\;j\in\cJ_\infty,
\end{eqnarray*}
where $T_1=T_1(L_\gamma/L_\beta)$  and $T_2=T_2(L_\gamma/L_\beta)$ can be easily deduced.

Thus, we assert that for any $j=1,\ldots, d$ and any $\e>0$ small enough
\begin{eqnarray}
\label{eq0050-new:proof-key-lemma}
&&\bar{\bleta}_j(m)\leq \mh_\e,\quad \forall m\leq \widehat{\mathbf{m}}.
\end{eqnarray}
Moreover, if $j\in\cJ_\pm$ we have in view of (\ref{eq00449-useful:proof-key-lemma})
$$
\big\{\bar{\bleta}_j(m)\big\}^{\gamma_j}\leq e^{-2\gamma_j} T_1
L_j^{-1}\varphi^{\frac{(2+1/\beta)\tau(q_j)}{2\tau(2)}}\to 0, \e\to 0,
$$
since $\tau(q_j)>0$ for any $j=1,\ldots d$ in view of $\tau(p^*)>0$.

Note also that if $\cJ_\infty\neq\emptyset$ then $p^*=\infty$ and therefore, $\tau(\infty)=1-1/\omega>0$ and, therefore,
$$
e^{-2\gamma_j}T_2
L_j^{-1}\varphi^{\frac{(2+1/\beta)(1-1/\omega)}{2\tau(2)}}\to 0,\;\e\to 0.
$$

Hence, for all $\e>0$ small enough
$
\bar{\bleta}_j(m)\leq \mh_\e,\; \forall m> \widehat{\mathbf{m}}.
$
Taking into account (\ref{eq0050-new:proof-key-lemma}), we conclude that (\ref{eq0042-new:proof-key-lemma}) and (\ref{eq00440-new:proof-key-lemma})   hold in the case $\kappa(p)\leq 0, \tau(p^*)> 0$. Moreover, (\ref{eq0046-new:proof-key-lemma}) is fulfilled if $m\leq \widehat{\mathbf{m}}$ as well in view of (\ref{eq0045-new:proof-key-lemma}).

On the other hand, in view of (\ref{eq0042-new:proof-key-lemma})
\begin{eqnarray}
\label{eq0051-new-new:proof-key-lemma}
e^{-d}\prod_{j=1}^d\mh_{s_j(m)}\leq \prod_{j=1}^d\bar{\bleta}_j(m)=e^{-2d}L_\beta^{-1}\varphi^{\frac{1}{\beta}}e^{-4dm}\leq\prod_{j=1}^d\mh_{s_j(m)},\quad\forall m>\widehat{\mathbf{m}},
\end{eqnarray}
and, therefore, (\ref{eq0046-new:proof-key-lemma}) is fulfilled for any  $m\geq 0$.
\epr

\subsection{{Proof of Lemma \ref{lem:technical}}}
Let $\mu\in (0,1)$ be the number whose choice will be done later and put $\vec{\mu}=(\mu,\ldots,\mu)$. Without loss of generality one can assume that
$\vec{\mu}\in\mS_d^{\text{const}}$. For any $f\in\bN_{\vec{r},d}\big(\vec{\beta},\vec{L}\big)$ introduce
$$
S_{\vec{\mu}}(x,f)=\int_{\bR^d}K_{\vec{\mu}}(t-x)f(t)\nu_d(\rd t),\;\;x\in \bR^d,
$$
where, recall,  $K$ is given in (\ref{eq:w-function}).

\smallskip

$1^{0}.\;$ Let us  prove that for any $\mu\in (0,1)$
\begin{equation}
\label{eq1:proof-lem:technical}
S_{\vec{\mu}}(\cdot,f)\in\bN^*_{\vec{r},d}\big(\vec{\beta},\mathbf{a}\vec{L}\big),\quad\;\forall f\in\bN_{\vec{r},d}\big(\vec{\beta},\vec{L}\big).
\end{equation}
First, we note that $S_{\vec{\mu}}(\cdot,f)$ is compactly supported on $(-b-1,b+1)^d$ in view of the definition of the kernel $K$, since $\mu\in (0,1)$.
Next, taking into account that $K$ is Lipschitz-continuous and compactly supported as well as $f\in \bL_{r^*}\big(\bR^d\big)$, $r^*=\max_{l=1,\ldots,d}r_l$,
since $f\in\bN_{\vec{r},d}\big(\vec{\beta},\vec{L}\big)$, and applying H\"older inequality, we can assert that $S_{\vec{\mu}}(\cdot,f)\in\bC\big(\bR^d\big)$
and moreover
\begin{equation}
\label{eq01:proof-lem:technical}
S_{\vec{\mu}}(\cdot,f)\in\bL_q\big(\bR^d\big),\quad \forall q\geq 1.
\end{equation}
Thus,  $S_{\vec{\mu}}(\cdot,f)\in\bC_{\bK}\big(\bR^d\big)$. It remains to prove that $S_{\vec{\mu}}(\cdot,f)\in\bN_{\vec{r},d}\big(\vec{\beta},\mathbf{a}\vec{L}\big)$. Indeed, applying the Young inequality we obtain for any
$j=1,\ldots, d$
$$
\big\|S_{\vec{\mu}}(\cdot,f)\big\|_{r_j,\bR^d}\leq \|K\|_{1,\bR^d}\|f\|_{r_j,\bR^d}\leq L_j\|K\|_{1,\bR^d},
$$
since $f\in\bN_{\vec{r},d}\big(\vec{\beta},\vec{L}\big)$.
Moreover, for any $j=1,\ldots, d$, $k\in\bN^*$ and any $u\in\bR$
\begin{eqnarray*}
\Delta_{u,j}^k S_{\vec{\mu}}(x,f)&:=&\Delta_{u,j}^k\bigg\{\int_{\bR^d}K_{\vec{\mu}}(z)f(x+z)\nu_d(\rd z)\bigg\}=
\int_{\bR^d}K_{\vec{\mu}}(z)\big\{\Delta_{u,j}^kf(x+z)\big\}\nu_d(\rd z)
\\
&=&\int_{\bR^d}K_{\vec{\mu}}(t-x)\big\{\Delta_{u,j}^kf(t)\big\}\nu_d(\rd t).
\end{eqnarray*}
Thus, applying the Young inequality, we have for any integer $k_j>\beta_j$
$$
\big\|\Delta_{u,j}^{k_j} S_{\vec{\mu}}(\cdot,f)\big\|_{r_j,\bR^d}\leq \|K\|_{1,\bR^d}\big\|\Delta_{u,j}^{k_j}f\big\|_{r_j,\bR^d}\leq L_j\|K\|_{1,\bR^d}|u|^{\beta_j},\quad\forall u\in\bR,
$$
since $f\in\bN_{\vec{r},d}\big(\vec{\beta},\vec{L}\big)$. Moreover,  $S_{\vec{\mu}}(\cdot,f)\in\bL_2\big(\bR^d\big)$ in view of (\ref{eq01:proof-lem:technical}).

We conclude that $S_{\vec{\mu}}(\cdot,f)\in\bN_{\vec{r},d}\big(\vec{\beta},\mathbf{a}\vec{L}\big)$
and, therefore, (\ref{eq1:proof-lem:technical}) is established.

\smallskip

$2^{0}.\;$ We will need the following auxiliary result. For any $(\vartheta,p)\in\cP^{\text{consist}}$ there exists $\mathrm{p}> p$ such that
\begin{equation}
\label{eq2:proof-lem:technical}
f\in\bN_{\vec{r},d}\big(\vec{\beta},\vec{L}\big)\;\Rightarrow\;f\in\bL_{\mathrm{p}}\big(\bR^d\big).
\end{equation}
Indeed, if $p^*>p$ we can choose $\mathrm{p}=p^*$ in view of the definition of an anisotropic Nikol'skii class.
If $p<2$ one can choose $\mathrm{p}=2$ since the definition of $\bN_{\vec{r},d}\big(\vec{\beta},\vec{L}\big)$ implies that $f\in\bL_{2}\big(\bR^d\big)$.

 It remains to consider the case $p^*=p$ and $p\geq 2$. Since $(\vartheta,p)\in\cP^{\text{consist}}$  necessarily  in this case $\tau(p)>0$ and, therefore, one can find
 $\mathrm{p}>p$ such that $\tau(\mathrm{p})>0$. In view of $p^*=p<\mathrm{p}$ and $\tau(\mathrm{p})>0$ the assertion of Lemma \ref{lem:nik-lep} holds with
 $s=\mathrm{p}$ and $\vec{r}(s)=(\mathrm{p},\ldots,\mathrm{p})$ and, therefore, $f\in\bL_{\mathrm{p}}\big(\bR^d\big)$ in view of the definition of an anisotropic Nikol'skii class. Thus, (\ref{eq2:proof-lem:technical}) is established.

\smallskip

$3^{0}.\;$ Let $f,g\in\bL_{p}\big(\bR^d\big)$ be arbitrary functions. We obviously have
$$
\sup_{\vec{h}\in\mS_d}\Big|\cB^{(p)}_{\vec{h}}(g)-\cB^{(p)}_{\vec{h}}(f)\Big|\leq 3\sup_{\vec{\mathrm{h}}\in\mS_d}
\big\|S_{\vec{\mathrm{h}}}(\cdot,g-f)\big\|_p+\|g-f\|_p.
$$
Since $K$ is compactly supported on $[-1/2,1/2]^d$ we obviously have that
$$
\big|S_{\vec{\mathrm{h}}}(x,g-f)\big|\leq \|K\|_{\infty,\bR^d}M[|g-f|](x),\quad x\in \bR^d.
$$
Applying $(p,p)$-strong maximal inequality (\ref{eq:strong-max}) we obtain for any $p>1$
$$
\big\|S_{\vec{\mathrm{h}}}(\cdot,g-f)\big\|_p\leq \bar{C}(p)\|K\|_{\infty,\bR^d}\|g-f\|_{p,\bR^d}.
$$
 Noting that the right hand side of the latter inequality is independent of $\vec{\mathrm{h}}$ we obtain finally
 \begin{equation*}
\sup_{\vec{h}\in\mS_d}\Big|\cB^{(p)}_{\vec{h}}(g)-\cB^{(p)}_{\vec{h}}(f)\Big|\leq (3\bar{C}(p)\|K\|_{\infty,\bR^d}+1)\|g-f\|_{p,\bR^d}.
\end{equation*}
Choosing $g=S_{\vec{\mu}}(\cdot,f)$ and noting that $\big|S_{\vec{\mu}}(\cdot,f)-f(\cdot)\big|=:B_{\vec{\mu}}(\cdot,f)$ we get
\begin{equation}
\label{eq3:proof-lem:technical}
\sup_{\vec{h}\in\mS_d}\Big|\cB^{(p)}_{\vec{h}}\big(S_{\vec{\mu}}(\cdot,f)\big)-\cB^{(p)}_{\vec{h}}(f)\Big|\leq (3\bar{C}(p)\|K\|_{\infty,\bR^d}+1)\big\|B_{\vec{\mu}}(\cdot,f)\big\|_{p,\bR^d}.
\end{equation}

$4^{0}.\;$ Some remarks are in order. First, $B_{\vec{\mu}}(\cdot,f)$ is compactly supported on $\bK$ for any $\mu\in(0,1)$. Next,
$B_{\vec{\mu}}(\cdot,f)\in\bL_{\mathrm{p}}\big(\bR^d\big)$ in view of (\ref{eq01:proof-lem:technical}) and (\ref{eq2:proof-lem:technical}).
At last, in view of (\ref{eq1:proof-cor}) and the first assertion of Lemma \ref{lem:bias-norm-bound} we have
$
\limsup_{\mu\to 0}\big\|B_{\vec{\mu}}(\cdot,f)\big\|_{1,\bR^d}=0.
$

All saying above allows us to assert that
$
\limsup_{\mu\to 0}\big\|B_{\vec{\mu}}(\cdot,f)\big\|_{p,\bR^d}=0.
$
It yields together with (\ref{eq3:proof-lem:technical}), that
for any $\varkappa>0$ and any $f\in\bN_{\vec{r},d}\big(\vec{\beta},\vec{L}\big)$
one can find
$\bmu=\bmu(\varkappa,f)$ such that
\begin{equation*}
\label{eq4:proof-lem:technical}
\sup_{\vec{h}\in\mS_d}\Big|\cB^{(p)}_{\vec{h}}\big(S_{\vec{\bmu}}(\cdot,f)\big)-\cB^{(p)}_{\vec{h}}(f)\Big|\leq \varkappa,
\end{equation*}
where as previously $\vec{\bmu}=(\bmu,\ldots,\bmu)$.

This obviously implies, for any $f\in\bN_{\vec{r},d}\big(\vec{\beta},\vec{L}\big)$ and any
$\bH\subseteq\mS_d$
\begin{eqnarray*}
\inf_{\vec{h}\in\bH}\Big[\cB^{(p)}_{\vec{h}}(f)+\e\Psi_{\e,p}\big(\vec{h}\big)\Big]&\leq&
\inf_{\vec{h}\in\bH}\Big[\cB^{(p)}_{\vec{h}}\big(S_{\vec{\bmu}}(\cdot,f)\big)+\e\Psi_{\e,p}\big(\vec{h}\big)\Big]+\varkappa
\\
&\leq&
\sup_{g\in\bN^*_{\vec{r},d}\big(\vec{\beta},\mathbf{a}\vec{L}\big)}\inf_{\vec{h}\in\bH}\Big[\cB^{(p)}_{\vec{h}}(g)+\e\Psi_{\e,p}\big(\vec{h}\big)\Big]+\varkappa,
\end{eqnarray*}
where to get the last inequality we have used (\ref{eq1:proof-lem:technical}). Since the right hand side of the latter inequality is independent of $f$ one gets
$$
\sup_{f\in\bN_{\vec{r},d}\big(\vec{\beta},\vec{L}\big)}\inf_{\vec{h}\in\bH}\Big[\cB^{(p)}_{\vec{h}}(f)+\e\Psi_{\e,p}\big(\vec{h}\big)\Big]\leq
\sup_{g\in\bN^*_{\vec{r},d}\big(\vec{\beta},\mathbf{a}\vec{L}\big)}\inf_{\vec{h}\in\bH}\Big[\cB^{(p)}_{\vec{h}}(g)+\e\Psi_{\e,p}\big(\vec{h}\big)\Big]+\varkappa,
$$
and the assertion of the lemma follows since $\varkappa$ is an arbitrary number.
\epr

\subsection{{Proof of Lemma \ref{lem:technical-new}}}

We obviously  have
$
U_{\vartheta,p}(x,g)\leq  b^*_{\vec{\mh}_{\mathbf{s}(0)}}(x,g)+\varpi_\e V^{-\frac{1}{2}}_{\mathbf{s}(0)}
$
and, therefore,
$$
U:=\sup_{g\in\bN^*_{\vec{r},d}\big(\vec{\beta},\mathbf{a}\vec{L}\big)}\big\|U_{\vartheta,p}(\cdot,g)\big\|_{p^*}\leq
\sup_{g\in\bN^*_{\vec{r},d}\big(\vec{\beta},\mathbf{a}\vec{L}\big)}\big\| b^*_{\vec{\mh}_{\mathbf{s}(0)}}(x,g)\big\|_{p^*}
+(2b)^{\frac{d}{p}}\varpi_\e V^{-\frac{1}{2}}_{\mathbf{s}(0)}.
$$
Note that  in view of (\ref{eq5020:proof-key-lemma})
$
\varpi_\e V^{-\frac{1}{2}}_{\mathbf{s}(0)}\to 0,\;\e\to 0
$
and, therefore, for all $\e>0$ small enough
\begin{equation}
\label{eq1:proof-lemma:technical-new}
U\leq
\Upsilon_1\sup_{g\in\bN^*_{\vec{r},d}\big(\vec{\beta},\mathbf{a}\vec{L}\big)}\big\| b^*_{\vec{\mh}_{\mathbf{s}(0)}}(x,g)\big\|_{p^*}
+L^*.
\end{equation}
Recall that
$
b^*_{\vec{h}}(x,g)=\sup_{J\in\mJ}\sup_{j=1,\ldots d}M_{J}\big[b_{\vec{h},j}\big](x)
$
and, therefore, we obtain first, applying (\ref{eq:strong-max-partial})
\begin{equation}
\label{eq2:proof-lemma:technical-new}
\big\| b^*_{\vec{\mh}_{\mathbf{s}(0)}}(\cdot,g)\big\|_{p^*}\leq
2^d\mathbf{C}_{p^*}\sum_{j=1}^d\big\| b_{\vec{\mh}_{\mathbf{s}(0)},j}\big\|_{p^*}.
\end{equation}
Next, we have for any $j=1,\ldots,d$ and any $x\in\bR^d$
\begin{eqnarray*}
b_{\vec{\mh}_{\mathbf{s}(0)},j}(x)&:=&\sup_{k:\:\mh_k\leq \mh_{s_j(0)}}
\left|\int_{\bR}w_\ell(u)g\big(x+u\mh_k\mathbf{e}_j\big)\nu_1(\rd u)-g(x) \right|
\\
&\leq& \|w_\ell\|_{\infty,\bR^d}M_{J_j}[g](x) +|g(x)|,
\end{eqnarray*}
where we have denoted  $J_j=\{1,\ldots, d\}\setminus\{j\}$. Thus, applying once again (\ref{eq:strong-max-partial}) we obtain
$$
\big\| b_{\vec{\mh}_{\mathbf{s}(0)},j}\big\|_{p^*}\leq \big(\mathbf{C}_{p^*}\|w_\ell\|_{\infty,\bR^d}+1\big)\|g\|_{p^*}.
$$
Noting that in view of the definition of the Nikolskii class  $\|g\|_{p^*}\leq \mathbf{a} L^*$ for any $g\in\bN^*_{\vec{r},d}\big(\vec{\beta},\mathbf{a}\vec{L}\big)$ and the assertion of the lemma follows from (\ref{eq1:proof-lemma:technical-new}) and
(\ref{eq2:proof-lemma:technical-new}).
\epr

\subsection{{Proof of formulas  (\ref{eq901:proof-key-lemma}) and  (\ref{eq0044901:proof-key-lemma})}}
\label{sec:subsec-equalities}

 In view of the definition of $\vec{\gamma}$ we have
$$
\frac{1}{\gamma}=\sum_{j\in J_\pm}\frac{\tau(r_j)}{\tau(p_\pm)\beta_j}+\sum_{j\in J_\infty}\frac{1}{\beta_j}\geq \frac{1}{\beta}
$$
since $\tau(r_j)\geq\tau(p_\pm)$. In view of the definition $\upsilon$
$$
\frac{p_\pm}{\upsilon}=\sum_{j\in J_\pm}\frac{1}{\gamma_j}\leq \sum_{j\in J_\pm}\frac{1}{\gamma_j}+\sum_{j\in J_\infty}\frac{1}{\beta_j}=\frac{1}{\gamma}
$$
and, therefore, $p_\pm\leq \upsilon/\gamma<\upsilon (2+1/\gamma)$.

\smallskip
\noindent {\bf Proof of (\ref{eq0044901:proof-key-lemma}).}\;
First, we remark that
\begin{eqnarray*}
\label{eq000101:proof-th}
&&p_\pm\Big(\frac{1}{\omega}-\frac{1}{\upsilon}\Big)+
\frac{1}{\gamma}-\frac{1}{\beta}=
\sum_{j\in J_\pm}\bigg(\Big[\frac{p_\pm}{r_j\beta_j}-\frac{1}{\gamma_j}\Big]+\Big[\frac{1}{\gamma_j}-\frac{1}{\beta_j}\Big]\bigg)
\\
&&=p_\pm\sum_{j\in J_\pm}\bigg(\frac{1}{r_j\beta_j}-\frac{1}{p_\pm\beta_j}\bigg)=:Ap_\pm.
\end{eqnarray*}

Next,
\begin{eqnarray*}
\sum_{j\in J_\pm}\frac{1}{\gamma_j}&=&\sum_{j\in J_\pm}\frac{\tau_j}{\tau(p_\pm)\beta_j}
= \frac{1}{\tau(p_\pm)}
\sum_{j\in J_\pm}\frac{1-1/\omega+1/(r_j\beta)}{\beta_j}
\\
&=&
\frac{1-1/\omega}{\tau(p_\pm)}\sum_{j\in J_\pm}\frac{1}{\beta_j}+\frac{1}{\tau(p_\pm)\beta}\sum_{j\in J_\pm}\bigg(\frac{1}{r_j\beta_j}-\frac{1}{p_\pm\beta_j}\bigg)
+\frac{1}{\tau(p_\pm)\beta p_\pm}\sum_{j\in J_\pm}\frac{1}{\beta_j}
\\
&=&\sum_{j\in J_\pm}\frac{1}{\beta_j}+\frac{A}{\tau(p_\pm)\beta}.
\end{eqnarray*}
It yields,
$
\frac{1}{\gamma}-\frac{1}{\beta}=\sum_{j\in J_\pm}\Big(\frac{1}{\gamma_j}-\frac{1}{\beta_j}\Big)=\frac{A}{\tau(p_\pm)\beta}
$
and, therefore,
$$
p_\pm\big(1/\omega-1/\upsilon\Big)=\big(1/\gamma-1/\beta\big)\big(\tau(p_\pm)\beta p_\pm-1\big)=
\big(1/\gamma-1/\beta\big)\beta p_\pm(1-1/\omega).
$$
The relation (\ref{eq0044901:proof-key-lemma}) is proved.
\epr


\bibliographystyle{agsm}

\end{document}